\documentclass[review]{elsarticle}
\makeatletter
\def\ps@pprintTitle{%
 \let\@oddhead\@empty
 \let\@evenhead\@empty
 \def\@oddfoot{\centerline{\thepage}}%
 \let\@evenfoot\@oddfoot}
\makeatother

\pdfoutput=1
\usepackage{graphicx}
\usepackage{subfigure}
\usepackage{amsfonts, amsmath, amssymb}
\usepackage{mathabx}
\usepackage{booktabs,siunitx}
\usepackage[svgnames,table]{xcolor}
\usepackage[tableposition=above]{caption}
\usepackage{pifont}
\usepackage{bm}
\usepackage{cancel}
\usepackage{breqn}
\usepackage{caption}
\usepackage{color}
\usepackage{enumerate}
\usepackage{float}
\usepackage{hyperref}
\usepackage{soul}
\usepackage[english]{babel}
\usepackage[margin=1in]{geometry}
\usepackage{mathtools}
\usepackage{multirow}
\usepackage{setspace}
\usepackage{subfigure}
\usepackage{stmaryrd}
\usepackage{tikz}
\usepackage{lineno}
\DeclareUnicodeCharacter{00B3}{\textsuperscript{3}}
\usepackage[ruled,linesnumbered]{algorithm2e}
\RestyleAlgo{ruled}
\SetKwComment{Comment}{/* }{ */}

\SetCommentSty{mycommfont}
\DeclareMathAlphabet{\mathpzc}{OT1}{pzc}{m}{it}

\newcommand \dd[2] {\frac{{\rm d} #1}{{\rm d} #2}}
\newcommand \dD[2] {\frac{{\rm d^2} #1}{{\rm d} {#2}^2}}

\newcommand \D [2]{\frac{\partial #1}{\partial #2}}
\newcommand \DD[2]{\frac{\partial^2 #1}{\partial #2 ^2}}
\newcommand \DDD [2]{\frac{{\rm D} #1}{{\rm D} #2}}

\renewcommand{\vec}[1]{\bm{\mathrm{#1}}}
\newcommand{\V}[1]{\bm{\mathrm{#1}}}

\def \div{\nabla \cdot \mbox{}}
\def \grad{\nabla}

\def \x{\vec{x}}

\def \n{\vec{n}}
\def \s{\vec{s}}
\def \u{\vec{u}}

\def \T{\vec{T}}

\def \L{\vec{L}}

\def \dt{\Delta t}
\def \dx{\Delta x}
\def \dy{\Delta y}

\def \rhos{\rho^{\rm S}}
\def \rhol{\rho^{\rm L}}

\def \rhog{\rho^{\rm{G}}}

\def \cps{C^{\rm S}}
\def \cpl{C^{\rm L}}
\def \cpg{C^{\rm G}}

\def \cpm{C^{\rm M}}

\def \Ts{T^{\rm S}}
\def \Tl{T^{\rm L}}
\def \Tum{T^{\rm M}}

\def \Tr{T_r}

\def \DeltaT{\Delta T}
\def \DeltaTstar{\Delta T^{*}}
\def \Tsol{T^{\rm sol}}
\def \Tliq{T^{\rm liq}}
\def \Tm{T_{m}}

\def \ks{\kappa^{\rm S}}
\def \kl{\kappa^{\rm L}}
\def \kg{\kappa^{\rm G}}
\def \km{\kappa^{\rm M}}

\def \hs{h^{\rm S}}

\def \hl{h^{\rm L}}

\def \mul{\mu^{\rm{L}}}
\def \mus{\mu^{\rm{S}}}
\def \mug{\mu^{\rm{G}}}

\def \ul{u^{\rm L}}
\def \us{u^{\rm S}}
\def \um{u^{\rm M}}

\def \Omegas{\Omega^{\rm S}}
\def \Omegal{\Omega^{\rm L}}

\def \Omegam{\Omega^{\rm M}}

\def \alphas{\alpha^{\rm S}}

\def \alphal{\alpha^{\rm L}}

\def \half{\frac{1}{2}}
\def \3half{\frac{3}{2}}

\def \hliq{h^{\rm liq}}
\def \hsol{h^{\rm sol}}
\def \hm{h^{\rm M}}

\def \sonet{s_1(t)}

\def \stwot{s_2(t)}

\def \lambdal{\lambda_l}
\def \lambdas{\lambda_s}
\def \lambdam{\lambda_m}
\def \lambdastar{\lambda^*}
\def \lambdalstar{\lambda_l^*}
\def \lambdasstar{\lambda_s^*}
\def \lambdamstar{\lambda_m^*}

\def \us{u^{\rm S}}
\def \ul{u^{\rm L}}

\def \Veta{V(\eta)}
\def \Vl{V_l}

\def \rhot{\rho_T}

\def \T0{T_0}
\def \Ti{T_i}
\def \Tj{T^j}
\def \ThetaS{\Theta^{\rm S}}
\def \ThetaM{\Theta^{\rm M}}
\def \ThetaL{\Theta^{\rm L}}
\def \alpham{\alpha^{\rm M}}
\def \Thetaj{\Theta^{j}}

\def \yone {y_1}
\def \ytwo {y_2}

\def \svarphit{s(t)}

\def \qmax {q_{\max}}

\graphicspath{{figures/}}

\usepackage{caption}  

\newcommand{\upperRomannumeral}[1]{\uppercase\expandafter{\romannumeral#1}}

\begin{document}
\let\today\relax
\let\underbrace\LaTeXunderbrace
\let\overbrace\LaTeXoverbrace

\begin{frontmatter}

\title{Self-similar structure of non-isothermal variable-density mushy Stefan problems and an improved low-Mach enthalpy method}
 \author[SDSU]{Yavkreet Swami}
\author[SDSU]{Amneet Pal Singh Bhalla\corref{mycorrespondingauthor}}
\ead{asbhalla@sdsu.edu}

\address[SDSU]{Department of Mechanical Engineering, San Diego State University, San Diego, CA}

\begin{abstract}

The enthalpy method was introduced in the late 1970s to simulate phase-change problems on fixed grids without explicitly tracking the moving phase-change front. It remains one of the most widely used approaches in academic and commercial software for the simulation of industrial melting and solidification processes.  For pure phase-change materials (PCMs) that melt or solidify at a single temperature, the enthalpy method introduces an artificial\footnote{In contrast, for metal alloys and glassy materials, $\Tsol$ and $\Tliq$ are physical, thermodynamic quantities.} mushy region bounded by the solidus temperature, $\Tsol$, and the liquidus temperature, $\Tliq$. As the numerical parameter $\Delta T=\Tliq-\Tsol$ approaches zero, the solution obtained with the enthalpy method is generally assumed to converge to that of the classical Stefan problem, in which the phase-change front is infinitesimally thin. This assumption is largely based on benchmark studies performed under the simplifying assumption of equal solid and liquid densities.

In this work, we systematically investigate the accuracy and spatio-temporal convergence properties of the enthalpy method for both low- and high-density-ratio phase-change problems. Because the limiting behavior $\Delta T\rightarrow0$ is difficult to realize numerically, owing to the diminishing thickness of the mushy region, we formulate and analyze the finite-$\Delta T$ mushy Stefan problem solved by the enthalpy method. We show that this problem possesses a self-similar structure that reduces the governing equations to a boundary-value problem involving two unknown parameters. The resulting similarity solution yields the temperature fields in the solid, liquid, and mushy regions, as well as the positions of the solidus, liquidus, and phase-change fronts. For unequal-density cases, the similarity solution additionally predicts spatially varying and spatially uniform velocity fields in the mushy and liquid regions, respectively, arising from density-change-induced flow. The analysis reveals that the enthalpy method solves the mushy Stefan problem, which differs fundamentally from the classical Stefan problem. In particular, the width of the mushy region grows proportionally to $\sqrt{t}$, whereas the phase-change front in the classical Stefan problem remains infinitesimally thin for all time. We show that, for large-density-ratio cases, $\Delta T$ must be sufficiently small for the mushy Stefan problem to closely approximate the sharp-interface Stefan problem. In contrast, for equal- and low-density-ratio cases, good agreement between the two problems is obtained over a wide range of $\Delta T$ values. The theoretical analysis also enables improvements to our previously developed low-Mach enthalpy method, enhancing its stability and accuracy as the density ratio between the two phases increases from $\mathcal{O}(1)$ to $\mathcal{O}(3)$.  A series of test cases with equal and unequal phase densities demonstrate excellent agreement between the improved low-Mach enthalpy method and the corresponding similarity solution of the mushy Stefan problem over a wide range of $\Delta T$ values. The improved method is further evaluated for large density contrasts between the solid and liquid metal phases, with particular emphasis on manufacturing defects such as pipe shrinkage and porosity. The results show that the improved low-Mach enthalpy method remains stable and produces physically accurate predictions under large density contrasts, whereas the original enthalpy method yields incorrect results and becomes unstable when predicting manufacturing defects involving large density ratios.

\end{abstract}

\begin{keyword}
\emph{diffuse-interface methods}\sep \emph{phase change} \sep \emph{metal manufacturing} \sep \emph{heat and mass transfer}  
\end{keyword}

\end{frontmatter}

\section{Introduction}\label{sec_introduction}

Modeling phase-change phenomena is central to understanding a wide range of natural and engineering processes, including ice-cap formation~\cite{buffo2021dynamics,buffo2021characterizing}, thermal energy storage~\cite{el2017thermal,allouhi2018optimization,badiei2020performance,hossain2019two}, casting and welding~\cite{king2015overview,khairallah2016laser}, and metal additive manufacturing~\cite{katayama2013handbook,blakey2021metal}. A prototypical mathematical model for melting and solidification is the classical Stefan problem, in which the phase-change front is represented as a sharp moving interface separating the solid and liquid phases~\cite{crank1984,alexiades2018mathematical}. The evolution of the interface is governed by the Stefan condition, which requires evaluating the jump in the temperature gradient across the phase boundary.

Directly tracking a sharp interface becomes challenging in practical melting and solidification problems because interfaces may merge, split, appear, or disappear during the simulation. Consequently, many commercial and research codes instead employ diffuse-interface formulations, in which the phase-change front is represented over a finite thickness and explicit interface tracking is avoided. Among the most widely used diffuse-interface approaches are phase-field methods based on the Allen--Cahn equation and the enthalpy method. Of these two approaches, the enthalpy method has been more widely adopted in commercial software because of its simplicity and its reliance on thermophysical properties that are typically available for engineering materials. In contrast, phase-field methods~\cite{huang2022consistent,zhang2022phasefieldicing} introduce additional length and time scales associated with the diffuse interface and require material parameters, such as the Cahn number and mobility coefficient, that are often unavailable for engineering materials and are therefore selected empirically. These parameters can significantly influence the predicted evolution of the phase-change front.

 Enthalpy methods were developed and refined throughout the late 1970s and 1980s for the simulation of phase-change problems. One of the earliest numerical solutions of the two-phase Stefan problem using an enthalpy formulation was reported by Crowley~\cite{crowley1978numerical}, who demonstrated the applicability of the method to multidimensional Stefan problems. In 1981, Voller and Cross~\cite{vollerCross1981} showed that accurate moving-boundary solutions can be obtained using the enthalpy formulation. The method was subsequently extended to more general convection--diffusion phase-change problems by Voller and co-workers~\cite{voller1987enthalpy}, allowing phase change to occur either at a single melting temperature (isothermal phase change) or over a finite temperature interval (non-isothermal phase change). The enthalpy methods have become a standard approach in commercial software packages such as FLOW-3D AM~\cite{FLOW-3D}, ANSYS Fluent~\cite{ansysFluentSolidificationMelting2025}, and AM PravaH~\cite{ampravahUserManual2024} for modeling industrial metal manufacturing processes (welding, casting, additive manufacturing) involving solidification and melting. However, most existing implementations are based on classical enthalpy formulations that assume equal densities for the solid and liquid phases. This limitation is significant because solid--liquid phase change is generally accompanied by volume change. In many commercial~\cite{FLOW-3D,cadence_simufact_additive,comsol_phase_change_material_63} and academic~\cite{schreterfleischhacker2025meltvapordynamics,li2025adaptive,chen2024meltrelocation,bennani2026monolithic} codes, shrinkage and expansion effects are either neglected or incorporated through empirical or phenomenological correction models~\cite{zanini2024analysing}, rather than being derived directly from the governing conservation laws. When the solid density, $\rhos$, differs from the liquid density, $\rhol$, phase change induces fluid motion, thereby coupling the thermal problem with the conservation of mass and momentum equations. Consequently, there is a need for modern enthalpy formulations that account for density-change effects in a thermodynamically and fluid-mechanically consistent manner.


Recent enthalpy-based formulations have begun to account for density-change effects in melting and solidification simulations. Zhao et al.~\cite{zhao2020generalized} developed an enthalpy-based lattice Boltzmann method for variable-density solid--liquid phase change by introducing a velocity-divergence constraint at the phase interface to account for volume change. The method was validated against the one-dimensional sharp-interface Stefan problem for density ratios up to $\rhol/\rhos=2$. Lyu et al.~\cite{lyu2021hybrid} investigated the influence of liquid-fraction--temperature relations and the phase-change interval $\DeltaT$ on ice formation problems. They reported that a linear liquid-fraction mapping produced accurate phase-front predictions relative to the classical Stefan solution for density ratios up to $\rhol/\rhos=1.67$, while smoother transition functions reduced start-up oscillations. More recently, Guo et al.~\cite{Guo2024DropletSolidificationLBM} studied droplet solidification using an enthalpy-based lattice Boltzmann method over a range of solid--liquid density ratios, $\rhos/\rhol = 0.8-1.3$. They showed that increasing the density ratio reduced the solidification time and altered the final droplet morphology, and validated the numerical predictions against experiments on hexadecane and water droplet solidification. Huang et al.~\cite{Huang2026DiffuseInterfaceFreezing} coupled a modified phase-field model for the gas--liquid interface with an enthalpy formulation for solid--liquid phase change to investigate freezing in fractures and porous media. Their model was validated against the classical Stefan problem for density ratios ranging from $0.6$ to $1.6$. In our previous works, we developed the low-Mach enthalpy method~\cite{thirumalaisamy2023low,thirumalaisamy2025consistent} to account for density-change-induced flow during melting and solidification and validated the method against the sharp-interface Stefan problem for density ratios up to 5.4.

Collectively, these studies demonstrate substantial progress in the development of variable-density enthalpy methods. However, their verification has remained largely based on the sharp-interface Stefan solution, corresponding to the limiting case $\DeltaT\rightarrow0$, and on density ratios close to unity. Consequently, the behavior of enthalpy methods over a broad range of phase-change intervals $\DeltaT$ and density contrasts has not been systematically investigated. 

To address these longstanding gaps in the enthalpy method literature, we first formulate a non-isothermal variable-density mushy Stefan problem for a finite phase-change interval, $\Delta T$. We show that the governing equations possess a self-similar structure that reduces the problem to a boundary-value problem. For the equal-density case, the similarity solution is obtained in closed form, whereas the unequal-density case requires the numerical solution of the resulting nonlinear boundary-value problem. Closed-form similarity solutions for equal-density mushy Stefan problems have previously been reported in the literature~\cite{solomon1982mushy,Ceretani2018SimilarityStefanMushy}. However, these studies assume an isothermal mushy region maintained at the phase-change temperature. As a result, their analytical solutions are not suitable for validating enthalpy methods, in which the mushy region is inherently non-isothermal due to the finite phase-change interval $\Delta T$. In the present work, we derive a closed-form analytical solution for the equal-density, non-isothermal mushy Stefan problem, thereby providing a more appropriate analytical benchmark for the verification of enthalpy methods.

The similarity analysis reveals that, for large-density-ratio cases, $\Delta T$ must be sufficiently small for the mushy Stefan problem to closely approximate the sharp-interface Stefan problem. In contrast, for equal- and low-density-ratio cases, good agreement between the two problems is maintained over a broad range of $\Delta T$ values. The similarity formulation also enables improvements to our previously developed low-Mach enthalpy method, yielding a formulation that remains stable and accurate as the density ratio between the two phases increases from $\mathcal{O}(1)$ to $\mathcal{O}(3)$. A series of verification cases with equal and unequal phase densities demonstrates excellent agreement between the improved low-Mach enthalpy method and the corresponding mushy Stefan similarity solution over a wide range of $\Delta T$ values. We further demonstrate that the improved low-Mach enthalpy method achieves a spatial–temporal convergence rate between first and second order, even at density ratios of $\mathcal{O}(3)$. Finally, metal-manufacturing defects, including pipe shrinkage during casting and porosity caused by gas entrapment during metal 3D printing, are simulated using the improved low-Mach enthalpy method. In these tests, the solid-to-liquid density ratio is systematically increased from \(\mathcal{O}(1)\) to \(\mathcal{O}(2)\) to examine the effects of volumetric shrinkage on the free metal surface and the resulting defects. The results demonstrate that the improved method remains stable for large density ratios and produces physically consistent predictions, whereas the original enthalpy method yields nonphysical results and becomes unstable.


\section{Enthalpy model formulation}\label{sec_specific_enthalpy}

Recent developments in the enthalpy method literature~\cite{lyu2021hybrid,zhao2020generalized,Guo2024DropletSolidificationLBM,Huang2026DiffuseInterfaceFreezing}, including our previous work on the low-Mach enthalpy method~\cite{thirumalaisamy2023low,thirumalaisamy2025consistent}, have incorporated variable thermophysical properties and density-change-induced flow directly into the governing equations. Nevertheless, the definition of the specific enthalpy within the mushy region remains a modeling choice. This choice determines the relationship among liquid fraction, temperature, and specific enthalpy, and therefore influences the evolution of the phase-change front. In this section, we examine two plausible definitions of the mushy region specific enthalpy and show that, although both satisfy standard thermodynamic requirements, they lead to markedly different predictions when the density ratio between the solid and liquid phases is large, e.g., $\mathcal{O}(2)$ or $\mathcal{O}(3)$.

\subsection{Mushy region specific enthalpy}
\label{sec_enthalpy_justification}

In the enthalpy method, phase change is assumed to occur over a finite temperature interval bounded by the temperatures $\Tsol \le T(\x,t) \le \Tliq$. The region $\Omegam(t)$ satisfying $\Tsol \le T \le \Tliq$ is referred to as the mushy region and consists of a mixture of solid and liquid phases. The liquid fraction is represented by the variable $\varphi(\x,t)$, such that $\varphi(\x \in \Omegas,t)=0$, $\varphi(\x \in \Omegal,t)=1$, and $0<\varphi(\x \in \Omegam,t)<1$. The computational domain $\Omega$ is partitioned into non-overlapping solid, mushy, and liquid regions, i.e., $\Omega=\Omegas(t)\cup\Omegam(t)\cup\Omegal(t)$. The assumption of a mushy region is physically justified for alloys and glassy materials, in which the phase transition occurs over a range of temperatures rather than at a single melting temperature, $\Tm$. In contrast, phase change in a pure material occurs at a single temperature.
 
 The specific enthalpies (enthalpy per unit mass) at the solidus and liquidus temperatures are defined as
\begin{subequations} \label{eq:hsol_hliq_general} 
\begin{alignat}{2} &\hsol = \cps(\Tsol-\Tr),\label{eq:hsol_general}\\ &\hliq = \hsol + L + \widebar{C}(\Tliq-\Tsol),\label{eq:hliq_general} 
\end{alignat} 
\end{subequations}
in which $L$ denotes the latent heat of melting, $\cps$, $\cpl$, and $\widebar{C}$ are the specific heat capacities of the solid, liquid, and mushy regions, respectively, and $\Tr$ is a reference temperature, typically chosen as the phase-change temperature in the bulk at the reference pressure. A simple choice for the effective heat capacity of the mushy region is the arithmetic average of the solid and liquid values, $\widebar{C}=\frac{\cps+\cpl}{2}$. The remaining modeling question is how to define the specific enthalpy in the mushy region, $\hm$.

To define a physically reasonable model, we require that the mushy region specific enthalpy, $\hm$,
\begin{enumerate}
    \item be continuous between $\hsol$ and $\hliq$;
    \item account for both sensible and latent heat contributions; and
    \item be such that the total enthalpy of the mushy region (per unit bulk volume) interpolates linearly between the total enthalpies of the solid and liquid phases.
\end{enumerate}
While the first two (thermodynamic) requirements are straightforward, the third requires further justification. The enthalpy transport equation governs the evolution of the specific enthalpy $h$ (or, equivalently, the temperature $T$). Because the enthalpy method identifies the mushy region through the liquid-fraction field $\varphi$, an additional relation is required to determine $\varphi$ and close the system. The third requirement yields the simplest closure, allowing the liquid fraction to be expressed as either $\varphi=\varphi(T)$ or $\varphi=\varphi(h)$. With the density of the mushy region defined by
\begin{equation}
\rho=\varphi\rhol + (1-\varphi)\rhos,
\label{eqn_rho_mixture}
\end{equation}
the third requirement can be written as
\begin{equation}
\rho \hm=(1-\varphi)\rhos\hsol+\varphi\rhol\hliq,
\label{eqn_mushy_enthalpy_closure}
\end{equation}
in which $\rhos$ and $\rhol$ denote the constant densities of the solid and liquid phases, respectively.

In the mushy region, latent heat is released in proportion to the mass of liquid present. Consequently, the latent-heat contribution to $\hm$ naturally takes the form $\varphi\rhol L/\rho$. The remaining modeling choice concerns the sensible-heat contribution to $\hm$. One possibility is to model this contribution as $\widebar{C}(T-\Tsol)$. Alternatively, the sensible heat may be weighted by the liquid mass fraction, analogous to the treatment of latent heat. This leads to the sensible-heat contribution as $\rhol \widebar{C}(T-\Tsol)/\rho$.

Both models satisfy the second requirement introduced above. The difference between them is negligible when $\rhos$ and $\rhol$ are comparable, but becomes significant when the density ratio $\rhol/\rhos$ deviates substantially from unity. We therefore consider the following two definitions for the mushy-region specific enthalpy:
\begin{subequations}\label{eqn_mushyh1_h2}
    \begin{alignat}{2}
        &\hm_1 = \widebar{C}(T-\Tsol)+\hsol+\varphi\frac{\rhol}{\rho}L, \qquad \Tsol \le T<\Tliq, \label{eqn_old_mushy_definition}\\
        &\hm_2 = \frac{\rhol}{\rho}\widebar{C}(T-\Tsol)+\hsol+\varphi\frac{\rhol}{\rho}L, \qquad \Tsol \le T<\Tliq.\label{eqn_new_mushy_definition}
    \end{alignat}
\end{subequations}

The corresponding piecewise definitions of the specific enthalpy in the PCM are
\begin{subequations}
    \begin{alignat}{2}
&h_1 = \begin{cases}
 \cps (T - \Tr) ,  &T<\Tsol,\\ 
 \widebar{C}(T - \Tsol) + \hsol + \varphi \frac{\displaystyle \rhol}{\displaystyle \rho}L,&\Tsol \le T < \Tliq, \\ 
 \cpl(T-\Tliq)+\hliq,& T \ge \Tliq ,
 \end{cases} \label{eqn_h_pcm_old} \\
&h_2 = \begin{cases}
 \cps (T - \Tr) ,  &T<\Tsol,\\ 
 \frac{\rhol}{\rho}\widebar{C}(T-\Tsol)+\hsol+\varphi\frac{\displaystyle \rhol}{\displaystyle \rho}L,&\Tsol \le T < \Tliq, \\ 
 \cpl(T-\Tliq)+\hliq,& T \ge \Tliq.
\end{cases}
\label{eqn_h_pcm_new}
\end{alignat}
\end{subequations}

To determine $\varphi(T)$, substitute $\hm$ from Eqs.~\eqref{eqn_mushyh1_h2} and $\rho$ from Eq.~\eqref{eqn_rho_mixture} into the closure relation Eq.~\eqref{eqn_mushy_enthalpy_closure}. This yields $\varphi_1 = \varphi_1(T)$ and $\varphi_2 = \varphi_2(T)$, corresponding to mushy specific enthalpy $\hm_1$ and $\hm_2$, respectively as 
\begin{subequations}
    \begin{alignat}{2}
        &\varphi_1(T)=\frac{\rhos(T-\Tsol)}{T(\rhos-\rhol)+\rhol\Tliq-\rhos\Tsol} ,\label{eqn_varphi_temp_old}\\
        &\varphi_2(T)=\frac{T-\Tsol}{\Tliq-\Tsol}.\label{eqn_varphi_temp_new}
    \end{alignat}
\end{subequations}

From Eq.~\eqref{eqn_varphi_temp_old}, it follows that $\varphi_1(T)$ is nonlinear whenever $\rhos\neq\rhol$. Moreover, the phase-change front, defined by the isotherm $T=\Tm$, corresponds to an iso-contour of $\varphi_1$ whose value depends on the densities of the solid and liquid phases. In contrast, $\varphi_2(T)$ is linear and independent of $\rhos$ and $\rhol$. Consequently, the value of $\varphi_2$ associated with the isotherm $T=\Tm$ remains fixed. In particular, when $\Tsol$ and $\Tliq$ are chosen symmetrically about $\Tm$, the phase-change front is given by $\varphi_2=0.5$ regardless of the values of $\rhos$ and $\rhol$. Note that when $\rhos \approx \rhol$, $\varphi_1(T) \approx \varphi_2(T)$.

The formulation defined by Eqs.~\eqref{eqn_old_mushy_definition} and~\eqref{eqn_h_pcm_old} was introduced and analyzed in our previous works \cite{thirumalaisamy2023low,thirumalaisamy2025consistent} and will hereafter be referred to as the original low-Mach enthalpy formulation. The formulation defined by Eqs.~\eqref{eqn_new_mushy_definition} and~\eqref{eqn_h_pcm_new} will be referred to as the improved low-Mach enthalpy formulation. In our previous studies, the original formulation was tested for density ratios between the solid and liquid phases of $\mathcal{O}(1)$ and was found to produce accurate solutions for a range of benchmark problems. However, when the density ratio is increased to $\mathcal{O}(2)$ or $\mathcal{O}(3)$, the original formulation exhibits spurious instabilities in the evolution of the phase-change front; see Fig.~\ref{fig_Ratio540Diff}. In contrast, the improved formulation predicts stable and accurate phase-front evolution, as demonstrated in Sec.~\ref{subsec_ratio540DiffDeltaT}.

To assess the accuracy of the enthalpy method, a reference solution is required. Existing studies typically benchmark enthalpy methods against the classical sharp-interface Stefan problem. Such comparisons are appropriate for assessing the limiting behavior as $\DeltaT\rightarrow0$. However, practical enthalpy simulations are performed with a finite $\DeltaT$. In the next section, we show that the enthalpy method solves a fundamentally different problem, namely the mushy Stefan problem. By exploiting its self-similar structure, we derive a similarity solution that is obtained numerically from a boundary-value problem. For the equal-density case, the similarity solution reduces to a closed-form analytical solution.

\section{Mushy Stefan problem: self-similar structure and similarity solution}\label{sec_mushy_stefan}

\begin{figure}[]
\centering{\includegraphics[scale=0.42]{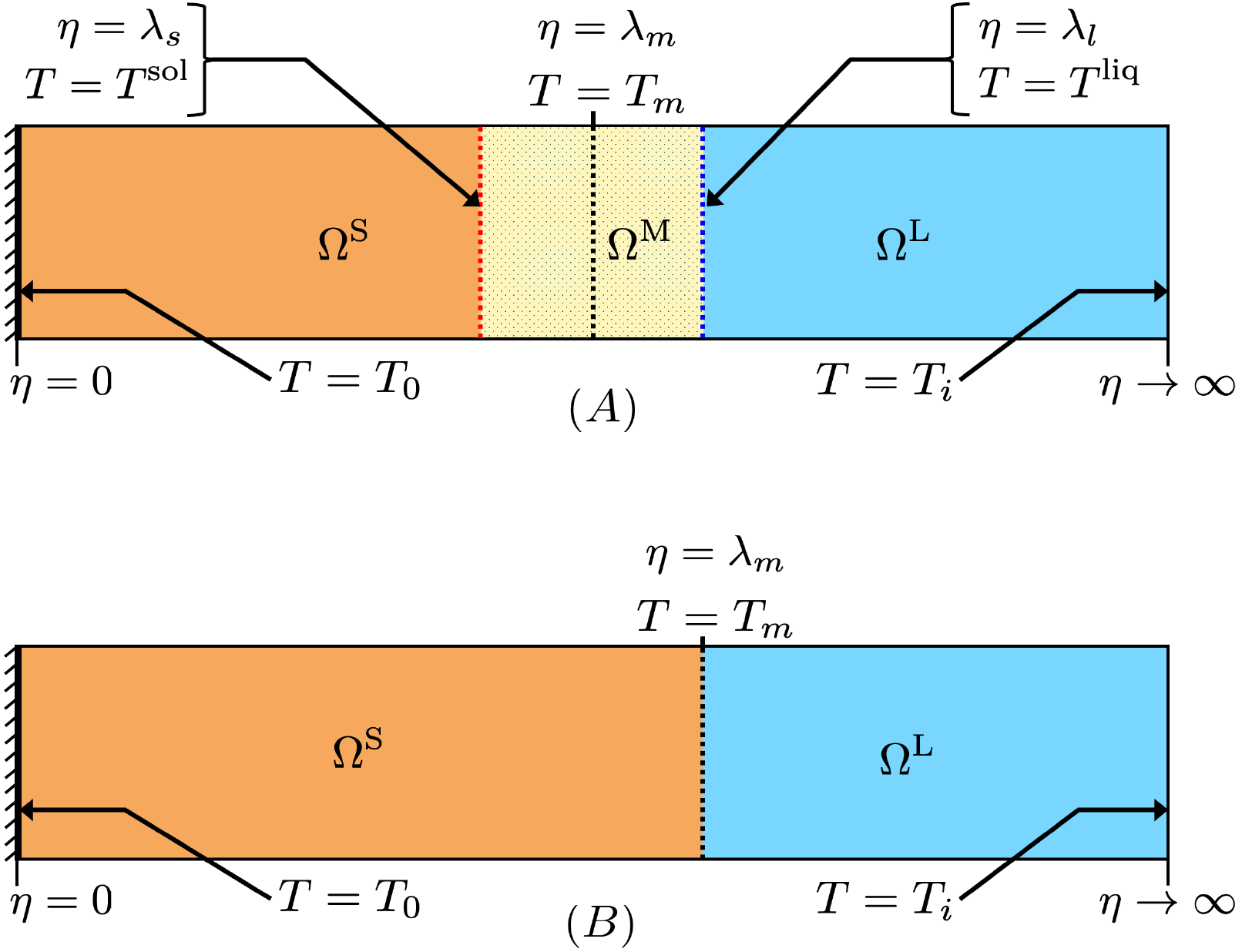}}
\caption{Schematic illustration of the (A) mushy Stefan problem and (B) sharp-interface Stefan problem for an initially liquid phase-change material (PCM) undergoing solidification from the left boundary at $t=0^{+}$.}\label{fig_Schematic}
\end{figure}

We consider two versions of the mushy Stefan problem in this section: (i) the equal-density case, in which the solid, liquid, and mushy regions all have the same density; and (ii) the unequal-density case, in which the solid and liquid phases have different but spatially uniform densities, while the density in the mushy region varies spatially. The former case admits a closed-form analytical solution, which serves as the initial guess for the boundary-value solver used to compute the latter.

\subsection{Equal-density mushy Stefan problem}\label{subsec_matched_density_mushy_stefan}

Consider a one-dimensional, semi-infinite solidification problem on the domain $\Omega:=\{x\ge 0\}$. Initially, the domain consists entirely of liquid at a uniform temperature $T_i$. At $t=0^{+}$, the temperature at the left boundary is suddenly reduced to $\T0<\Tsol$ and maintained at that value thereafter. As a result, solidification commences, and the solidus and liquidus fronts are assumed to emerge simultaneously at $t=0^{+}$ and propagate in the positive $x$-direction; see Fig.~\ref{fig_Schematic}(A). The positions of the solidus and liquidus fronts are denoted by $x=\sonet$ and $x=\stwot$, respectively, such that the spatial regions
\begin{subequations}\label{eqn_mushy_stefan_regions}
    \begin{alignat}{2}
        0 <\; & x < \sonet \quad && \in \Omegas,\\
        \sonet <\; & x < \stwot \quad && \in \Omegam,\\
        & x > \stwot \quad && \in \Omegal.
    \end{alignat}
\end{subequations}
The solidus, liquidus, and phase-change fronts are defined by the isotherms
\begin{subequations}
    \begin{alignat}{2}
       & T(\sonet,t) &&= \Tsol, \\
       & T(\stwot,t) &&= \Tliq, \\
       & T(s(t),t) &&= \Tm.
    \end{alignat}
\end{subequations}

In phase $j\in\{\rm{S,L,M}\}$, the one-dimensional mass-balance and enthalpy equations take the form
\begin{subequations} \label{eqn_mass_enthalpy}
    \begin{alignat}{2}
        &\D{\rho^j}{t} + \D{(\rho^j u^j)}{x} = 0,\label{eqn_mass}\\
        & \D{(\rho^j h^j)}{t} + \D{(\rho^j u^j h^j)}{x}
        = \kappa^j \DD{T^j}{x}. \label{eqn_enthalpy}
   \end{alignat}       
\end{subequations}

For the equal-density case, there is no density-change-induced flow, i.e.,
$u^j(x\in\Omega^j,t)=0$. This follows directly from the mass-conservation equation~\eqref{eqn_mass}, which reduces to
$\partial{u^j}/\partial{x} = 0$ in each region. Together with the zero-velocity condition at $x=0$ and continuity of velocity across the solidus and liquidus interfaces, the velocity vanishes everywhere in the domain, i.e., $u(x \in \Omega, t) = 0$.  Substituting the solid- and liquid-phase specific enthalpies,
\begin{subequations}\label{eqn_solid_liquid_enthalpy_matched}
    \begin{alignat}{2}
        &\hs = \cps(T-\Tr),\\
        &\hl = \hliq + \cpl(T-\Tliq),
    \end{alignat}
\end{subequations}
into the enthalpy equation~\eqref{eqn_enthalpy} yields
\begin{subequations}
    \begin{alignat}{2}
         &\rho\cps\D{\Ts}{t}
         = \ks\DD{\Ts}{x},
         &&\qquad 0<x<\sonet,
         \label{eqn_heat_solid_matched}\\
         &\rho\cpl\D{\Tl}{t}
         = \kl\DD{\Tl}{x},
         &&\qquad x>\stwot.
         \label{eqn_heat_liquid_matched}
    \end{alignat}
\end{subequations}

Next, we derive the heat equation in the mushy region using the improved enthalpy definition given by Eq.~\eqref{eqn_new_mushy_definition}. Since $\rhos=\rhol=\rho$, the mushy region specific enthalpy reduces to
\begin{equation}
    \hm = \hsol + \widebar{C}(T-\Tsol) + \varphi L.
    \label{eqn_hm_matched_density}
\end{equation}
Using the linear liquid-fraction relation
$\varphi(T)=(T-\Tsol)/(\Tliq-\Tsol)$,
Eq.~\eqref{eqn_hm_matched_density} can be rewritten as
\begin{equation}
    \hm = \hsol + \left(\widebar{C} + \frac{L}{\Tliq-\Tsol}\right)(T-\Tsol).
    \label{eqn_hm_matched_density_simplified}
\end{equation}
Therefore, the derivative of the mushy region specific enthalpy with respect to temperature, which defines the effective specific heat capacity of the mushy region, is constant and given by
\begin{equation}
    \cpm = \dd{\hm}{T}
    = \widebar{C} + \frac{L}{\Tliq-\Tsol}.
    \label{eqn_dhmdT_matched_density}
\end{equation}

Substituting Eq.~\eqref{eqn_hm_matched_density_simplified} into the enthalpy equation~\eqref{eqn_enthalpy}, and applying the chain rule,
$\D{\hm}{t}=\dd{\hm}{T}\D{\Tum}{t}$,
yields the heat equation in the mushy region:
\begin{equation}
    \rho \cpm \D{\Tum}{t}
    = \km\DD{\Tum}{x},
    \qquad \sonet < x < \stwot.
    \label{eqn_heat_mush_matched}
\end{equation}

Combining the solid-, mushy-, and liquid-region heat equations, the equal-density mushy Stefan problem is governed by
\begin{subequations}
    \begin{alignat}{2}
        &\rho\cps\D{\Ts}{t}
        =\ks\DD{\Ts}{x},
        &&\qquad 0<x<\sonet,
        \label{eqn_heat_solid_matched_summary}\\
        &\rho\left(\widebar{C}+\frac{L}{\Tliq-\Tsol}\right)\D{\Tum}{t}
        =\km\DD{\Tum}{x},
        &&\qquad \sonet < x < \stwot,
        \label{eqn_heat_mush_matched_summary}\\
        &\rho\cpl\D{\Tl}{t}
        =\kl\DD{\Tl}{x},
        &&\qquad x>\stwot,
    \label{eqn_heat_liquid_matched_summary}
    \end{alignat}
\label{eqn_heat_equations_matched_summary}
\end{subequations}
in which $\ks$, $\km$, and $\kl$ denote the thermal conductivities of the solid, mushy, and liquid regions, respectively. The thermal conductivity of the mushy region is taken as the arithmetic average of the solid and liquid conductivities,
$ \km=\frac{\ks+\kl}{2}$.

To obtain a self-similar solution, we seek temperature fields of the form
\begin{equation}
    \Ts(x,t)=\ThetaS(\eta), \qquad
    \Tum(x,t)=\ThetaM(\eta), \qquad
    \Tl(x,t)=\ThetaL(\eta),
    \label{eqn_temperature_similarity_fields}
\end{equation}
in which $\eta$ is the similarity variable defined by
\begin{equation}
    \eta=\frac{x}{2\sqrt{t}}.
    \label{eqn_eta_definition}
\end{equation}
The interface locations are parameterized by the similarity constants $\lambdas$ and $\lambdal$ as
\begin{equation}
    \sonet=2\lambdas\sqrt{t}, \qquad
    \stwot=2\lambdal\sqrt{t}.
    \label{eqn_interface_similarity}
\end{equation}
Accordingly, the solid, mushy, and liquid regions in similarity space are given by
\begin{equation}
    0<\eta<\lambdas, \qquad
    \lambdas<\eta<\lambdal, \qquad
    \eta>\lambdal.
    \label{eqn_eta_regions}
\end{equation}

Defining $\Tj(x,t)=\Thetaj(\eta)$, the required derivatives are
\begin{subequations}
    \begin{alignat}{2}
        &\D{\Tj}{t}
        =-\frac{\eta}{2t}\dd{\Thetaj}{\eta},
        \label{eqn_Tt_similarity}\\
        &\D{\Tj}{x}
        =\frac{1}{2\sqrt{t}}\dd{\Thetaj}{\eta},
        \label{eqn_Tx_similarity}\\
        &\DD{\Tj}{x}
        =\frac{1}{4t}\dD{\Thetaj}{\eta}.
        \label{eqn_Txx_similarity}
    \end{alignat}
    \label{eqn_similarity_derivatives}
\end{subequations}

Substituting these expressions into the governing heat equations yields
\begin{subequations}
    \begin{alignat}{2}
        &\alphas\dD{\ThetaS}{\eta}
        +2\eta\dd{\ThetaS}{\eta}=0,
        &&\qquad 0<\eta<\lambdas,
        \label{eqn_solid_similarity_ode_matched}\\
        &\alpham\dD{\ThetaM}{\eta}
        +2\eta\dd{\ThetaM}{\eta}=0,
        &&\qquad \lambdas<\eta<\lambdal,
        \label{eqn_mush_similarity_ode_matched}\\
        &\alphal\dD{\ThetaL}{\eta}
        +2\eta\dd{\ThetaL}{\eta}=0,
        &&\qquad \eta>\lambdal,
        \label{eqn_liquid_similarity_ode_matched}
    \end{alignat}
    \label{eqn_similarity_odes_matched}
\end{subequations}
in which $\alpha^j=\kappa^j/(\rho C^j)$ denotes the thermal diffusivity of region $j$.
The temperature boundary and interfacial conditions are
\begin{subequations}\label{eqn_temperature_bcs_matched}
    \begin{alignat}{2}
    &\ThetaS(0)=\T0,
    \label{eqn_theta_s_left_bc}\\
    &\ThetaS(\lambdas)=\ThetaM(\lambdas)=\Tsol,
    \label{eqn_theta_solidus_bc}\\
    &\ThetaM(\lambdal)=\ThetaL(\lambdal)=\Tliq,
    \label{eqn_theta_liquidus_bc}\\
    &\ThetaL(\infty)=\Ti.
    \label{eqn_theta_farfield_bc}
    \end{alignat}
\end{subequations}

The solutions of Eqs.~\eqref{eqn_similarity_odes_matched} satisfying the boundary and interfacial conditions~\eqref{eqn_temperature_bcs_matched} are
\begin{subequations}\label{eqn_theta_solutions_matched}
    \begin{alignat}{2}
        &\ThetaS(\eta)
        =
        \T0+(\Tsol-\T0)
        \frac{\operatorname{erf}\left(\eta/\sqrt{\alphas}\right)}
        {\operatorname{erf}\left(\lambdas/\sqrt{\alphas}\right)},
        &&\qquad 0<\eta<\lambdas,
        \label{eqn_theta_s_solution_matched}\\
        &\ThetaM(\eta)
        =
        \Tsol+(\Tliq-\Tsol)
        \frac{
        \operatorname{erf}\left(\eta/\sqrt{\alpham}\right)
        -
        \operatorname{erf}\left(\lambdas/\sqrt{\alpham}\right)}
        {
        \operatorname{erf}\left(\lambdal/\sqrt{\alpham}\right)
        -
        \operatorname{erf}\left(\lambdas/\sqrt{\alpham}\right)},
        &&\qquad \lambdas<\eta<\lambdal,
        \label{eqn_theta_m_solution_matched}\\
        &\ThetaL(\eta)
        =
        \Ti+(\Tliq-\Ti)
        \frac{\operatorname{erfc}\left(\eta/\sqrt{\alphal}\right)}
        {\operatorname{erfc}\left(\lambdal/\sqrt{\alphal}\right)},
        &&\qquad \eta>\lambdal.
        \label{eqn_theta_l_solution_matched}
    \end{alignat}
\end{subequations}

The similarity constants $\lambdas$ and $\lambdal$ are determined by enforcing continuity of heat flux at the solidus and liquidus interfaces. Unlike the classical Stefan problem, latent heat in the mushy Stefan problem is released continuously over a finite region, resulting in continuous heat flux across both interfaces. Accordingly,
\begin{subequations}\label{eqn_flux_conditions_matched}
    \begin{alignat}{2}
        &\ks\left.\dd{\ThetaS}{\eta}\right|_{\eta=\lambdas}
        =
        \km\left.\dd{\ThetaM}{\eta}\right|_{\eta=\lambdas},
        \label{eqn_flux_solidus_matched}\\
        &\km\left.\dd{\ThetaM}{\eta}\right|_{\eta=\lambdal}
        =
        \kl\left.\dd{\ThetaL}{\eta}\right|_{\eta=\lambdal}.
        \label{eqn_flux_liquidus_matched}
    \end{alignat}
\end{subequations}

Substituting the similarity solutions~\eqref{eqn_theta_solutions_matched} into the flux-continuity conditions~\eqref{eqn_flux_conditions_matched} yields the coupled transcendental system
\begin{subequations}\label{eqn_transcedental_match_density} 
\begin{align} 
\frac{\ks(\Tsol-\T0)}{\sqrt{\alphas}}\frac{\exp(-\lambdas^2/\alphas)}{\operatorname{erf}\left(\lambdas/\sqrt{\alphas}\right)}
&=\frac{\km(\Tliq-\Tsol)}{\sqrt{\alpham}}\frac{\exp(-\lambdas^2/\alpham)}{\operatorname{erf}\left(\lambdal/\sqrt{\alpham}\right)-\operatorname{erf}\left(\lambdas/\sqrt{\alpham}\right)},\label{eqn_lambda_s_matched}\\
\frac{\km(\Tliq-\Tsol) }{\sqrt{\alpham}}\frac{\exp(-\lambdal^2/\alpham)}{\operatorname{erf}\left(\lambdal/\sqrt{\alpham} \right) - \operatorname{erf} \left(\lambdas/\sqrt{\alpham}\right)}
&=\frac{\kl(\Ti-\Tliq)}{\sqrt{\alphal}}\frac{\exp(-\lambdal^2/\alphal) }{\operatorname{erfc}\left( \lambdal/\sqrt{\alphal} \right)}.
\label{eqn_lambda_l_matched} 
\end{align} 
\end{subequations}

The transcendental system~\eqref{eqn_transcedental_match_density} is solved numerically using MATLAB's \texttt{fsolve} function to determine $\lambdas$ and $\lambdal$. Once these similarity constants are known, the positions and velocities of the solidus and liquidus fronts follow directly from Eq.~\eqref{eqn_interface_similarity}. The width of the mushy region is given by
\begin{equation}
    w(t)
    =
    \stwot-\sonet
    =
    2(\lambdal-\lambdas)\sqrt{t}.
\end{equation}
Thus, the width of the mushy region grows proportionally to $\sqrt{t}$, in contrast to the classical Stefan problem, in which the phase-change front remains infinitesimally thin for all time.

To compare the phase-front location predicted by the mushy Stefan problem with that of the classical Stefan problem, we solve Eq.~\eqref{eqn_theta_m_solution_matched} for the similarity coordinate $\eta=\lambdam$ satisfying $\ThetaM(\lambdam)=\Tm$. Fig.~\ref{fig_Schematic}(B) shows the schematic representation of the sharp-interface Stefan problem. The parameter $\lambdam$, which lies within the interval $\lambdas<\lambdam<\lambdal$, identifies the location of the phase-change front within the mushy region. The corresponding physical position of the phase-change front is
\begin{equation}
    s(t)=2\lambdam\sqrt{t}.
    \label{eqn_st_pc}
\end{equation}

\subsection{Unequal-density mushy Stefan problem}\label{subsec_unmatched_density_mushy_stefan}

We now consider the more general case in which the solid and liquid densities are different, i.e., $\rhos \neq \rhol$. In this case, phase change induces fluid motion due to density variations. The velocity in the solid phase remains zero, $\us = 0$, whereas the liquid and mushy regions exhibit non-zero velocities that contribute to convective heat transport.

The density in the mushy region, given by Eq.~\eqref{eqn_rho_mixture}, can be rewritten as
\begin{equation}
\rho(\Tum)=\rhos + \rhot(\Tum-\Tsol),\label{eqn_mushy_density_unmatched_T}
\end{equation}
after substituting $\varphi(T)$ from Eq.~\eqref{eqn_varphi_temp_new}. Here, $\rhot = \dd{\rho}{T} = (\rhol-\rhos)/(\Tliq-\Tsol)$. Rewriting the improved specific enthalpy, Eq.~\eqref{eqn_new_mushy_definition}, as a function of temperature yields
\begin{equation}
\hm=\hsol+ \frac{\rhol \widebar{C}(\Tum-\Tsol)}{\rho(\Tum)} + \frac{\rhol L (\Tum-\Tsol)}{\rho(\Tum)(\Tliq-\Tsol)},\label{eqn_hm}
\end{equation}
Substituting Eq.~\eqref{eqn_hm} into the non-conservative form of the enthalpy equation~\eqref{eqn_enthalpy} gives
\begin{subequations}
\begin{alignat}{2}
\rho(\Tum) \DDD{\hm}{t}= \rho(\Tum) \dd{\hm}{\Tum}\DDD{\Tum}{t} = \Gamma(\Tum)\DDD{\Tum}{t} &=  \km\DD{\Tum}{x}, \label{eqn_non-cons-enthalpy} \\
\hookrightarrow \; \Gamma(\Tum) \left(\D{\Tum}{t} + \um\D{\Tum}{x}\right) &=\km\DD{\Tum}{x}.
\label{eqn_mush_heat_unmatched_physical}
\end{alignat}
\end{subequations}
Substituting the expressions for the density and specific enthalpy into Eq.~\eqref{eqn_non-cons-enthalpy} yields the effective volumetric heat capacity of the mushy region,
\begin{equation}
\Gamma(\Tum) = \rho(\Tum) \dd{\hm}{\Tum} = \frac{\rhol \rhos A}{\rho(\Tum)}.
\end{equation}
Here, $A = \widebar{C} + L/(\Tliq - \Tsol)$. Note that the effective heat capacity of the mushy region, $\cpm = \dd{\hm}{T} = \Gamma(T)/\rho(T) = \rhol \rhos A/\rho(T)^2$, is now temperature dependent. This contrasts with the equal-density case, for which $\cpm = A$ is constant; see Eq.~\eqref{eqn_dhmdT_matched_density}.

In the solid region, the velocity vanishes, $\us=0$, and therefore the heat equation remains unchanged from Eq.~\eqref{eqn_heat_solid_matched}. In the liquid region, the density is constant. The mass-conservation equation therefore reduces to $\partial \ul/\partial x = 0$, implying that the liquid velocity is spatially uniform. The heat equation in the liquid region then takes the form
\begin{equation}
\rhol\cpl \left(\D{\Tl}{t}+\ul\D{\Tl}{x}\right)=\kl\DD{\Tl}{x} .
\label{eqn_liquid_heat_unmatched}
\end{equation}

Introducing the similarity variable $\eta=x/(2\sqrt{t})$ and parameterizing the interface locations by the similarity constants $\lambdas$ and $\lambdal$, we seek velocity fields in the mushy and liquid regions of the form
\begin{subequations}
\begin{alignat}{2}
& \um(x,t)=\frac{\Veta}{\sqrt{t}}, &&\qquad \lambdas<\eta<\lambdal, \\
& \ul(x,t)=\frac{\Vl}{\sqrt{t}}, &&\qquad \eta > \lambdal.
\label{eqn_velocity_similarity_unmatched}
\end{alignat}
\end{subequations}
Here, $\Veta$ denotes the mushy region velocity function, which satisfies the conditions $V(\lambdas)=0$ and $V(\lambdal)=\Vl$. We again seek temperature fields of the form $\Ts(x,t)=\ThetaS(\eta)$, $\Tum(x,t)=\ThetaM(\eta)$, and $\Tl(x,t)=\ThetaL(\eta)$ in the three regions.

Following steps analogous to those of the previous section, the solid-region similarity equation and its solution remain unchanged and are given by Eqs.~\eqref{eqn_solid_similarity_ode_matched} and~\eqref{eqn_theta_s_solution_matched}, respectively. Substituting Eq.~\eqref{eqn_velocity_similarity_unmatched} into Eq.~\eqref{eqn_liquid_heat_unmatched} yields the liquid-region similarity equation,
\begin{equation}
\alphal\dD{\ThetaL}{\eta} +2(\eta-\Vl)\dd{\ThetaL}{\eta}=0.
\label{eqn_liquid_similarity_unmatched}
\end{equation}
Subject to the boundary conditions $\ThetaL(\lambdal)=\Tliq$ and $\ThetaL(\infty)=\Ti$, the liquid-region temperature admits the closed-form solution
\begin{equation}
\ThetaL(\eta)=\Ti+(\Tliq-\Ti)\frac{\operatorname{erfc}\left((\eta-\Vl)/\sqrt{\alphal}\right) }{\operatorname{erfc}\left((\lambdal-\Vl)/\sqrt{\alphal}\right)}.
\label{eqn_theta_l_unmatched}
\end{equation}

The mushy region similarity equations are obtained by applying the similarity transformation to the mass and enthalpy equations~\eqref{eqn_mass_enthalpy}. Substituting $\rho=\rho(\ThetaM(\eta))$ and $\um=\Veta/\sqrt{t}$ into Eq.~\eqref{eqn_mass} yields
\begin{align}
& \rho_t + (\rho \um)_x = 0, \nonumber \\
\hookrightarrow & -\frac{\eta}{2t}\rho_\eta + \frac{1}{2t}(\rho \Veta)_\eta = 0,  \nonumber \\
\hookrightarrow    & \dd{\Veta}{\eta}=\frac{\rhot}{\rho(\ThetaM)}\dd{\ThetaM}{\eta}\left(\eta-\Veta\right).
\label{eqn_velocity_ode_unmatched}
\end{align}
In deriving Eq.~\eqref{eqn_velocity_ode_unmatched}, we have used
$\rho_t=-(\eta\rho_\eta)/(2t)$,
$(\rho \um)_x=(\rho \Veta)_\eta/(2t)$,
and
$\rho_\eta=\rhot\ThetaM_\eta$. Eq.~\eqref{eqn_velocity_ode_unmatched} governs the similarity velocity field within the mushy region. 

Similarly, applying the similarity transformation to the mushy region enthalpy equation~\eqref{eqn_mush_heat_unmatched_physical} yields
\begin{align}
& \Gamma(\ThetaM)\left(-\frac{\eta}{2t}\ThetaM_\eta + \frac{V(\eta)}{\sqrt{t}}\frac{1}{2\sqrt{t}}\ThetaM_\eta\right)
 = \km \frac{1}{4t}\ThetaM_{\eta \eta}, \nonumber \\
\hookrightarrow \; & \km \dD{\ThetaM}{\eta} + 2\Gamma(\ThetaM)(\eta-V)\dd{\ThetaM}{\eta} = 0.
\label{eqn_mush_energy_ode_unmatched}
\end{align}

Equations~\eqref{eqn_velocity_ode_unmatched} and~\eqref{eqn_mush_energy_ode_unmatched} form a coupled nonlinear system governing the mushy region temperature $\ThetaM(\eta)$ and velocity $V(\eta)$. It is, however, possible to eliminate $V$ from Eq.~\eqref{eqn_mush_energy_ode_unmatched}, thereby reducing the problem to a single nonlinear differential equation for $\ThetaM$. This reduction is presented next.

Solving the mushy region temperature equation~\eqref{eqn_mush_energy_ode_unmatched} for $\eta-V$ gives
\[
\eta-V = -\frac{\km\ThetaM_{\eta \eta}}{2\Gamma(\ThetaM)\ThetaM_\eta}.
\]
Substitution into Eq.~\eqref{eqn_velocity_ode_unmatched} yields
\begin{subequations}
\begin{alignat}{2}
&\rho(\ThetaM)V_\eta + \rhot\frac{\km\ThetaM_{\eta \eta}}{2\Gamma(\ThetaM)} = 0, \\
\hookrightarrow \; &\rho(\ThetaM)V_\eta + \rhot\frac{\km\rho(\ThetaM)}{2\rhol\rhos A}\ThetaM_{\eta \eta} = 0, \\
\hookrightarrow \; &   V_\eta = -\frac{\rhot \km}{2\rhol\rhos A}\ThetaM_{\eta \eta}. \label{eqn_VM_eta}
\end{alignat}
\end{subequations}

Integrating Eq.~\eqref{eqn_VM_eta} and applying the velocity boundary condition at the solidus front, $V(\lambdas)=0$, gives
\begin{equation}
V(\eta) = \frac{\rhot \km}{2\rhol\rhos A}
\left[\ThetaM_\eta(\lambdas)-\ThetaM_\eta(\eta)\right].
\label{eqn_VM}
\end{equation}

Substituting this expression into the mushy region temperature equation~\eqref{eqn_mush_energy_ode_unmatched} gives
\begin{equation}
\km\ThetaM_{\eta \eta} + 2\Gamma(\ThetaM)
\left(
\eta - \frac{\rhot \km}{2\rhol\rhos A}
\left[\ThetaM_\eta(\lambdas)-\ThetaM_\eta(\eta)\right]
\right)\ThetaM_\eta = 0. \label{eqn_mushy_reduced_ode}
\end{equation}
This is the reduced mushy region similarity ODE with the velocity eliminated. Eq.~\eqref{eqn_mushy_reduced_ode} is a second-order nonlinear ODE in $\eta$ with two unknown parameters, $\lambdas$ and $\lambdal$. Therefore, four boundary conditions are required to close the problem:
\begin{subequations}
\begin{alignat}{2}
\ThetaM(\lambdas) &=\Tsol, \\
\ThetaM(\lambdal) &=\Tliq, \\
\ks\left.\dd{\ThetaS}{\eta}\right|_{\eta=\lambdas} &=\km\left.\dd{\ThetaM}{\eta}\right|_{\eta=\lambdas},\label{eqn_flux_solidus_unmatched}\\
\kl\left.\dd{\ThetaL}{\eta}\right|_{\eta=\lambdal} &= \km\left.\dd{\ThetaM}{\eta}\right|_{\eta=\lambdal}.\label{eqn_flux_liquidus_unmatched}
\end{alignat}
\label{eqn_unmatched_mush_bcs}
\end{subequations}
The left-hand sides of the flux conditions~\eqref{eqn_flux_solidus_unmatched} and~\eqref{eqn_flux_liquidus_unmatched} can be evaluated analytically using Eqs.~\eqref{eqn_theta_s_solution_matched} and~\eqref{eqn_theta_l_unmatched}, respectively.

\subsection{Fixed-domain numerical solution to the similarity equation} \label{subsec_fixed_domain_unmatched}

To numerically solve the nonlinear mushy region ODE~\eqref{eqn_mushy_reduced_ode}, we map the interval $\lambdas\le \eta \le \lambdal$ to the fixed computational domain $z\in[0,1]$ using the transformation
\begin{equation}
z = \frac{\eta-\lambdas}{\lambdal-\lambdas}, \qquad\Delta_\lambda=\lambdal-\lambdas .
\label{eqn_fixed_domain_map_unmatched}
\end{equation}
We then introduce the state variables
\begin{equation}
\yone=\ThetaM,\qquad \ytwo=\dd{\ThetaM}{\eta}.
\label{eqn_state_variables_unmatched}
\end{equation}
Since
\begin{equation}
\dd{}{z}=\Delta_\lambda\dd{}{\eta},
\label{eqn_z_eta_derivative}
\end{equation}
the governing equation can be written as the first-order system
\begin{subequations}
\begin{alignat}{2}
&\dd{\yone}{z}=\Delta_\lambda \ytwo,
\label{eqn_bvp_y1_unmatched}\\
&\dd{\ytwo}{z}=-\Delta_\lambda\frac{2 \Gamma(\yone)}{\km}\left( \eta-V(\eta) \right)\ytwo.
\label{eqn_bvp_y2_unmatched}
\end{alignat}
\label{eqn_bvp_system_unmatched}
\end{subequations}
The constitutive relations appearing in Eqs.~\eqref{eqn_bvp_system_unmatched} are
\begin{subequations}
\begin{alignat}{2}
\eta(z) &= \lambdas + \Delta_\lambda z, \\
\Gamma(\yone) &= \frac{\rhol\rhos A}{\rhos + \rhot(\yone - \Tsol)}, \\
V(\eta) &= \frac{\rhot \km}{2\rhol\rhos A}
\left[ \frac{\ks}{\km}\ThetaS_\eta(\lambdas)-\ytwo\right].
\end{alignat}
\label{eqn_bvp_constitutive_unmatched}
\end{subequations}

The system of first-order nonlinear equations~\eqref{eqn_bvp_system_unmatched}, together with the boundary conditions~\eqref{eqn_unmatched_mush_bcs}, is solved using MATLAB's boundary-value solver \texttt{bvp5c}. The solver is based on a collocation method with adaptive mesh refinement and supports the inclusion of unknown parameters, allowing the similarity constants $\lambdas$ and $\lambdal$ to be determined as part of the solution procedure. As an initial guess for $\lambdas$ and $\lambdal$, we solve the transcendental system~\eqref{eqn_transcedental_match_density} corresponding to the equal-density case, taking $\rho=\min(\rhos,\rhol)$. The initial guess for $\ThetaM(\eta)$ is chosen as a linear profile connecting $\Tsol$ and $\Tliq$. Once $\ThetaM(\eta)$ has been obtained, the velocity field in the mushy region is recovered from Eq.~\eqref{eqn_VM}. The numerical solution for $\ThetaM(\eta)$ also allows the similarity coordinate of the phase-change front, $\eta=\lambdam$, to be determined by solving $\ThetaM(\lambdam)=\Tm$. The position of the phase-change front is $s(t) = 2\lambdam\sqrt{t}$ (see  Eq.~\eqref{eqn_st_pc}).

\noindent\textbf{\underline{Continuation strategy for the boundary value solver}}: For sufficiently challenging parameter combinations, MATLAB's boundary-value solver \texttt{bvp5c} may fail to converge from a given initial guess. This typically occurs when the density ratio between the two phases exceeds approximately $1000$ or when the phase-change interval $\Delta T$ becomes very small. In such cases, we employ a continuation strategy. Specifically, we first solve a nearby problem with either a smaller density ratio or a larger value of $\Delta T$, both of which are easier to converge. The converged solution is then used as the initial guess for the target problem. By repeatedly applying this continuation procedure, we are able to obtain similarity solutions for density ratios as large as $2700$ and phase-change intervals as small as $0.1$ K.

\subsection{Comparison of similarity solutions for the mushy and sharp-interface Stefan problems}\label{subsec_sharpvsmushy_results}

\begin{table}[]
\centering
\caption{Thermophysical properties used in the one-dimensional Stefan problem. Except for the solid density, which is artificially reduced to investigate large-density-ratio cases, all material properties correspond to those of aluminum and are taken from Refs.~\cite{doble2007perry,hatch1984aluminium,desai1987thermodynamic}.}

\label{tab_thermophys_properties_Al}
\begin{tabular}{ll l}
Property & Value & Units\\
\midrule
Thermal conductivity of solid   $\ks$     &   211  & W/m.K  \\
Thermal conductivity of liquid  $\kl$     & 91   & W/m.K    \\
Specific heat of solid  $\cps$            & 910 & J/kg.K      \\
Specific heat of liquid  $\cpl$           & 1042.4 & J/kg.K   \\
Melting temperature  $T_m$                & 933.6  & K     \\
Reference temperature  $T_r$              & 933.6  & K     \\
Latent heat of melting $L$              &  383840 & J/kg    \\
\bottomrule
 \end{tabular}
\end{table}

Here, we compare the similarity solutions of the mushy and sharp-interface Stefan problems for a one-dimensional solidification problem with both equal and unequal phase densities. The thermophysical properties of the phase-change material (PCM) are based primarily on aluminum and are summarized in Table~\ref{tab_thermophys_properties_Al}. The similarity solution of the sharp-interface Stefan problem is available from our previous works~\cite{thirumalaisamy2023low,thirumalaisamy2025consistent}. For both models, the similarity coordinate corresponding to the phase-change isotherm, $T=T_m$, is denoted by $\lambdam$. For the mushy Stefan problem, the width of the mushy region in similarity space is given by $\lambdal-\lambdas$. Figs.~\ref{fig_Schematic}(A) and (B) schematically illustrate the mushy and sharp-interface Stefan problems, respectively. 

To compare results at different density ratios, we nondimensionalize the phase-change interval, temperature, similarity coordinate, and similarity constants as 
\begin{subequations}
\begin{alignat}{2}
\Delta T^* &= \frac{\Delta T}{T_m}, \\
T^* &= \frac{T-T_i}{T_m-T_i}, \\
\eta^* &= \frac{\eta}{\sqrt{\alpha_{\mathrm{ref}}}}, \\
\lambda_j^* &= \frac{\lambda_j}{\sqrt{\alpha_{\mathrm{ref}}}},
\quad j \in \{s,m,l\}.
\end{alignat}
\end{subequations}
Here, $\alpha_{\mathrm{ref}} = \frac{\alphas+\alphal}{2}$ is the reference (taken to be the average) thermal diffusivity. 


\subsubsection{Equal-density case} \label{subsec_no_volume_change}

\begin{figure}[]
	\centering{\includegraphics[width=\textwidth, height = 0.91\textheight,  trim={40 65 30 0},clip,keepaspectratio]{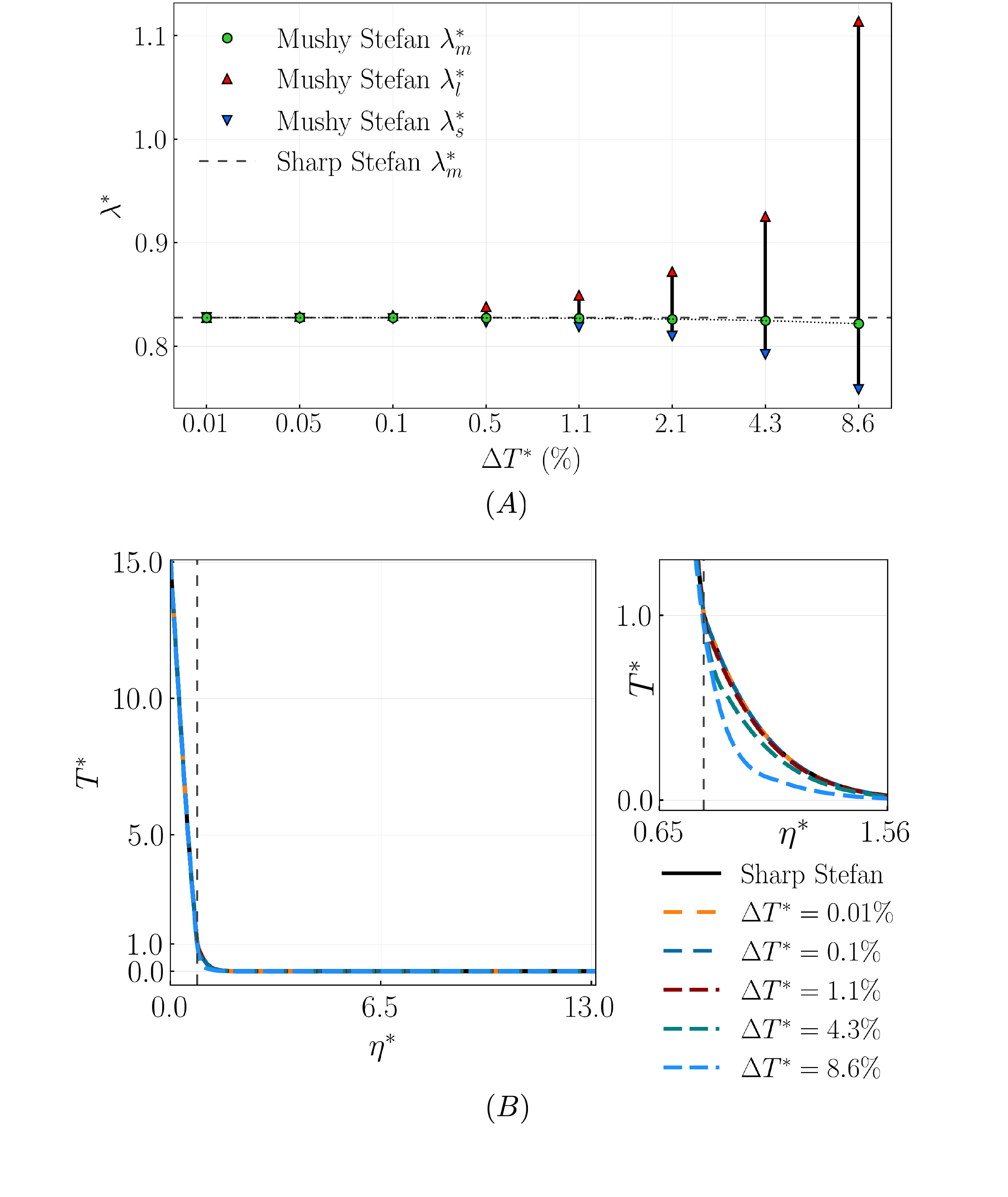}}
\caption{Comparison of the similarity solutions of the mushy and sharp-interface Stefan problems for the equal-density case ($\rhol/\rhos=1$) for different values of the nondimensional phase-change interval $\DeltaTstar(\%)=100 \frac{\DeltaT}{\Tm}$, in which $\Tm=933.6$ K, reported as a percentage of the melting temperature: (A) nondimensional phase-front location; and (B) nondimensional temperature profile. The green filled circles and the horizontal dashed line in (A) denote the nondimensional phase-front similarity coordinate $\lambdamstar$ predicted by the mushy and sharp-interface Stefan problems, respectively. The vertical bars terminated by triangular markers indicate the nondimensional mushy region interval $[\lambdasstar,\lambdalstar]$. The vertical dashed line in (B) denotes the sharp-interface phase-front location. }\label{fig_MushyvsSharpLowRatio}
\end{figure}

We first consider the equal-density case, $\rhos=\rhol=2700$ kg/m$^3$, for which there is no density-change-induced flow. The domain is initially liquid at temperature $\Ti=978.6$ K, with boundary temperatures $\T0=298.6$ K and $\Ti=978.6$ K imposed at the left and right boundaries, respectively. Similarity solutions for the phase-front position and temperature are presented in Fig.~\ref{fig_MushyvsSharpLowRatio}.

The results show that the phase-front location obtained from the mushy Stefan problem is essentially independent of $\DeltaTstar$ over the range considered. The corresponding value of $\lambdamstar$ remains in close agreement with that predicted by the sharp-interface Stefan problem. Although the phase-front location is largely insensitive to $\DeltaTstar$, the width of the mushy region increases with increasing $\DeltaTstar$, as indicated by the larger vertical bars in Fig.~\ref{fig_MushyvsSharpLowRatio}(A). The vertical bar represents the extent of the mushy region in the (nondimensional) similarity space, with its lower and upper ends corresponding to $\lambdastar =\lambdasstar$ and $ \lambdastar =\lambdalstar$, respectively. Thus, for the equal-density case, reducing $\DeltaTstar$ decreases the mushy region width without significantly affecting the phase-front location. The temperature similarity profiles shown in Fig.~\ref{fig_MushyvsSharpLowRatio}(B) further demonstrate that the mushy and sharp-interface solutions remain in close agreement over the entire range considered, even when the phase-change interval reaches approximately $8.5\%$ of the melting temperature.

In the literature, enthalpy methods are most commonly validated against the equal-density Stefan problem~\cite{hahn2012heat}. The results presented here explain why such validations generally exhibit excellent agreement: for the equal-density case, the mushy and sharp-interface Stefan problems yield nearly identical phase-front locations and temperature fields, even when the phase-change interval is several percent of the melting temperature. As a result, the equal-density Stefan problem does not provide a particularly stringent validation case for modern enthalpy methods designed to model density-change-induced flow.


\subsubsection{Low-density-ratio volume-expansion case}\label{subsec_Low_volume_expansion_case}

\begin{figure}[]
	\centering{\includegraphics[width=\textwidth, height = 0.91\textheight,  trim={40 65 30 0},clip,keepaspectratio]{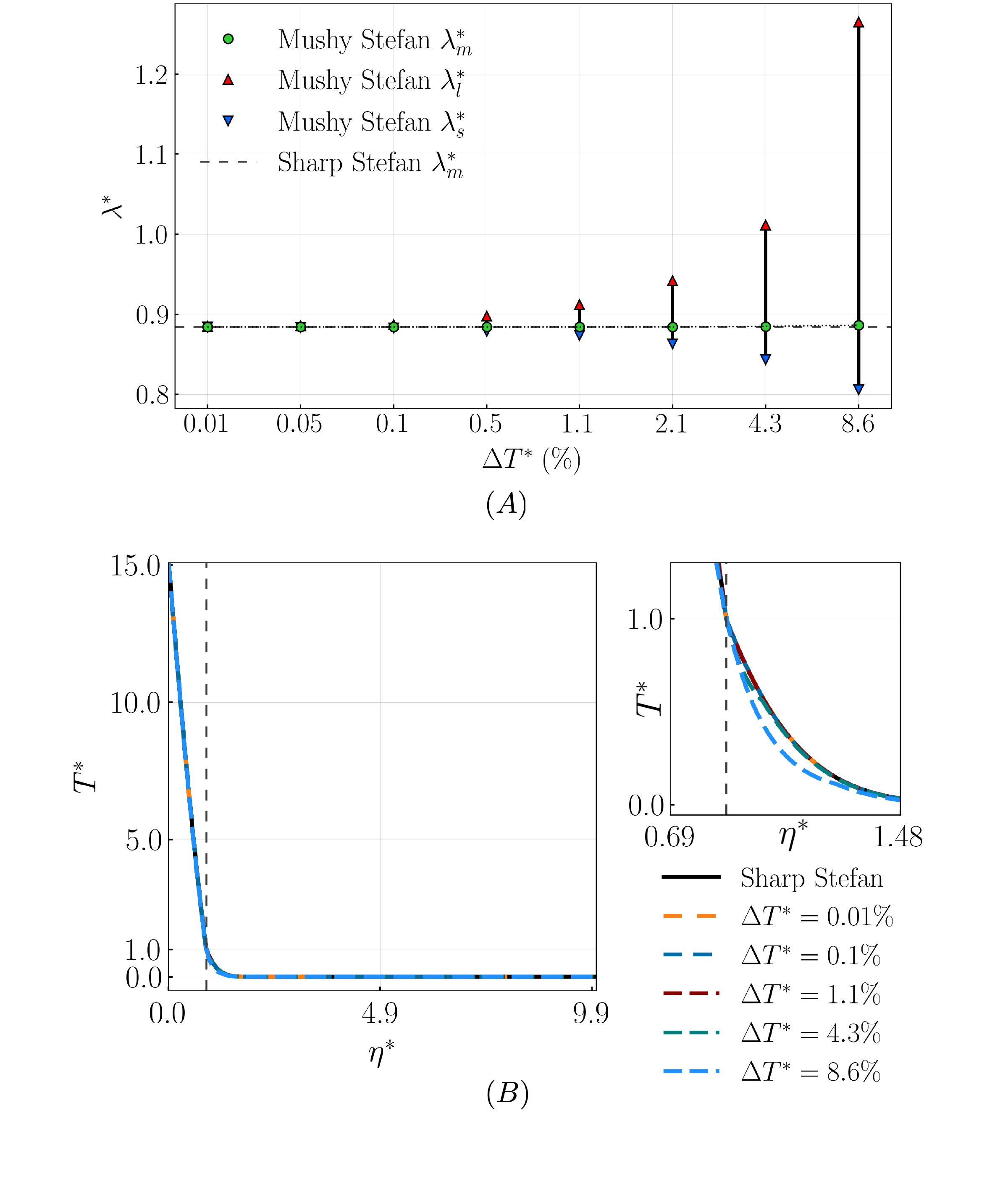}}
\caption{Comparison of the similarity solutions of the mushy and sharp-interface Stefan problems for the low-density-ratio case ($\rhol/\rhos=2$) for different values of the nondimensional phase-change interval $\DeltaTstar(\%)= 100 \frac{\DeltaT}{\Tm}$, in which $\Tm=933.6$ K, reported as a percentage of the melting temperature: (A) nondimensional phase-front location; and (B) nondimensional temperature profile. The green filled circles and the horizontal dashed line in (A) denote the nondimensional phase-front similarity coordinate $\lambdamstar$ predicted by the mushy and sharp-interface Stefan problems, respectively. The vertical bars terminated by triangular markers indicate the nondimensional mushy region interval $[\lambdasstar,\lambdalstar]$. The vertical dashed line in (B) denotes the sharp-interface phase-front location. }\label{fig_MushyvsSharpLowRatioTemperatureComp}
\end{figure}

Next, we consider a moderate density-ratio case with $\rhos=1350$ kg/m$^3$ and $\rhol=2700$ kg/m$^3$, corresponding to a density ratio of two. Since the solid density is lower than the liquid density, the material expands upon solidification. The domain is initially liquid at temperature $\Ti=978.6$ K, with boundary temperatures $\T0=298.6$ K and $\Ti=978.6$ K imposed at the left and right boundaries, respectively.

Fig.~\ref{fig_MushyvsSharpLowRatioTemperatureComp}(A) shows that the nondimensional phase-front location predicted by the mushy Stefan problem remains in close agreement with the sharp-interface Stefan solution over the entire range of $\DeltaTstar$ considered. As in the equal-density case, increasing $\DeltaTstar$ increases the width of the mushy region, as indicated by the larger vertical bars. However, this increase has little effect on the phase-front location, indicating that for density ratios of order unity the mushy Stefan problem reproduces the sharp-interface solution even when the phase-change interval reaches approximately $8.5\%$ of the melting temperature. The temperature similarity profiles shown in Fig.~\ref{fig_MushyvsSharpLowRatioTemperatureComp}(B) exhibit a similar trend. The mushy and sharp-interface Stefan solutions remain in close agreement over the entire range of $\DeltaTstar$ values considered. Thus, even low-density-ratio cases do not provide a stringent validation benchmark for enthalpy methods aimed at resolving volume-change effects during phase change.

\subsubsection{Large-density-ratio volume-expansion case}\label{subsec_Large_volume_expansion_case}

\begin{figure}[]
	\centering{\includegraphics[width=\textwidth, height = 0.91\textheight,  trim={40 65 30 0},clip,keepaspectratio]{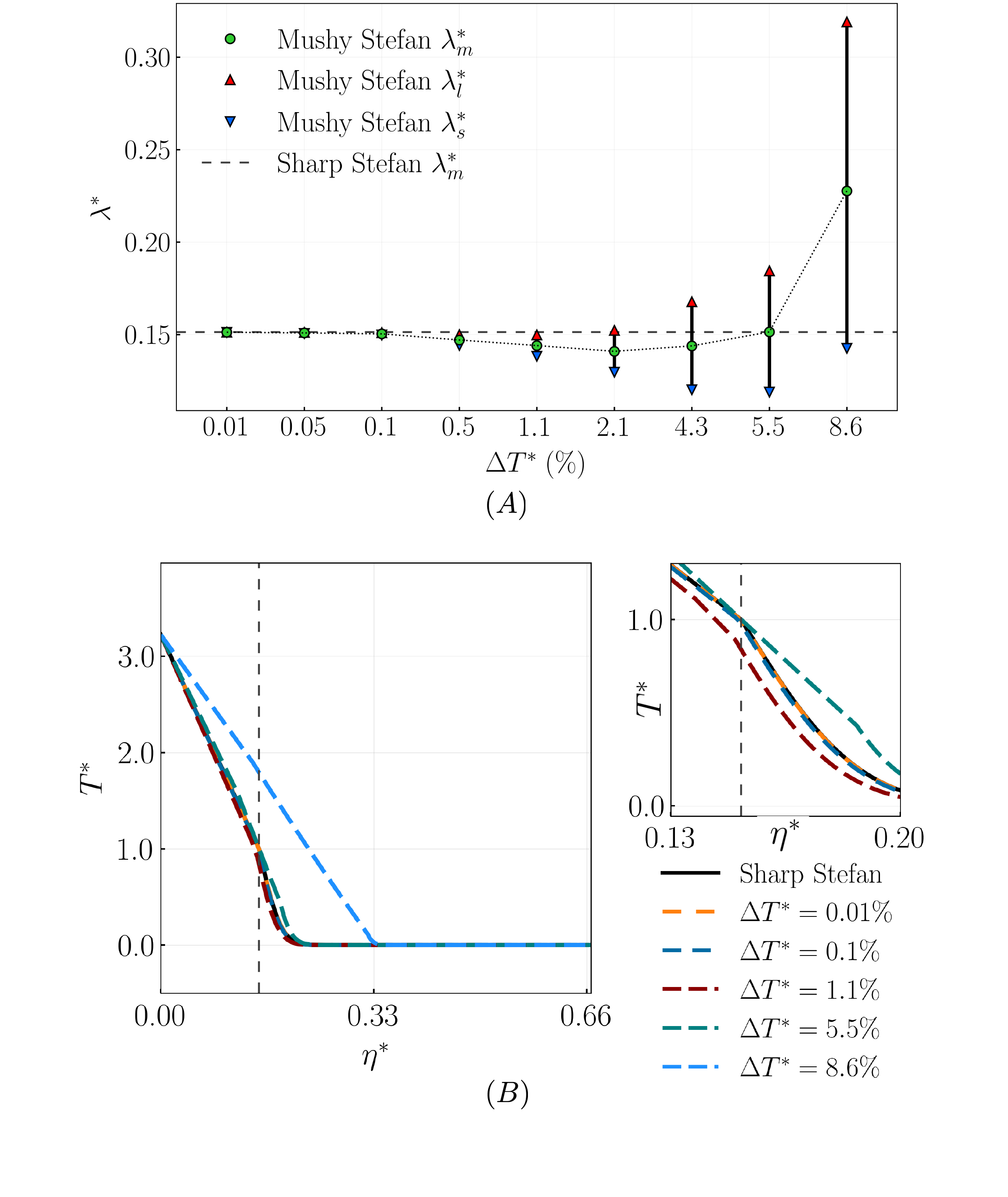}}
\caption{Comparison of the similarity solutions of the mushy and sharp-interface Stefan problems for the large-density-ratio case ($\rhol/\rhos=540$) for different values of the nondimensional phase-change interval $\DeltaTstar (\%)= 100 \frac{\DeltaT}{\Tm}$, in which $\Tm=933.6$ K, reported as a percentage of the melting temperature: (A) nondimensional phase-front location; and (B) nondimensional temperature profile. The green filled circles and the horizontal dashed line in (A) denote the nondimensional phase-front similarity coordinate $\lambdamstar$ predicted by the mushy and sharp-interface Stefan problems, respectively. The vertical bars terminated by triangular markers indicate the nondimensional mushy region interval $[\lambdasstar,\lambdalstar]$. The vertical dashed line in (B) denotes the sharp-interface phase-front location.}\label{fig_MushyvsSharpLargeRatio}
\end{figure}

We next consider a large-density-ratio case with $\rhos=5$ kg/m$^3$ and $\rhol=2700$ kg/m$^3$, corresponding to $\rhol/\rhos=540$. The domain is initially liquid at $\Ti=978.6$ K, with the left boundary maintained at $\T0=833.6$ K and the right boundary at $\Ti$. Relative to the previous cases, the large density contrast generates substantially stronger density-change-induced flow. Consequently, the phase-front location becomes much more sensitive to the choice of $\DeltaTstar$, as shown in Fig.~\ref{fig_MushyvsSharpLargeRatio}(A).

The results reveal a non-monotonic dependence of the mushy Stefan phase-front location on $\DeltaTstar$. As $\Delta T\rightarrow0$, the mushy Stefan solution approaches the sharp-interface limit, as expected. As $\DeltaTstar$ increases, the phase-front location initially deviates from the sharp-interface solution, then moves back toward it, and eventually deviates again at larger values of $\DeltaTstar$. For the present case, a phase-change interval of $\DeltaT = 51.2$ K, corresponding to approximately $5.5\%$ of $\Tm$, yields close agreement between the mushy and sharp-interface phase-front locations despite the relatively large mushy region width. This agreement, however, should not be interpreted as convergence to the sharp-interface limit. Rather, it corresponds to a finite-width mushy Stefan solution whose phase-front location happens to coincide with that of the sharp-interface Stefan problem. The value of $\DeltaTstar$ for which this occurs depends on the density ratio, thermophysical properties, and imposed initial and boundary conditions.

Fig.~\ref{fig_MushyvsSharpLargeRatio}(B) compares the nondimensional temperature profiles of the mushy and sharp-interface Stefan problems for different values of $\DeltaTstar$. When the phase-change interval is approximately $1\%$ of $\Tm$ or smaller, the temperature profiles are nearly indistinguishable, confirming the expected sharp-interface limit. As the phase-change interval increases to one percent or more, noticeable differences emerge between the two solutions. This remains true even for the case of $\DeltaT=51.2$ K, corresponding to approximately $5.5\%$ of $\Tm$, despite the close agreement in the phase-front location.

These results demonstrate that large-density-ratio cases provide a significantly more stringent benchmark for validating modern enthalpy methods than equal- or low-density-ratio problems. In the next section, we present an improved low-Mach enthalpy method designed to accurately resolve density-change-induced flow in the presence of large density contrasts between the two phases. The method is validated against the mushy Stefan problem for a range of $\Delta T$ values, and its spatio-temporal convergence properties are also examined.

\section{Improved low-Mach enthalpy method}\label{sec_improved_low_mach}

In this section, we briefly describe the main aspects of the low-Mach enthalpy method introduced in~\cite{thirumalaisamy2023low,thirumalaisamy2025consistent} and the improvements proposed in the present work. The method also enables the coupling of a passive gas phase with the PCM domain using a level-set-based interface-capturing method. The low-Mach enthalpy method extends classical enthalpy formulations by accounting for density-change-induced flow directly from the governing conservation laws, thereby eliminating the need for reduced-order models to represent shrinkage and expansion during phase change. The method was previously validated on several benchmark problems with density ratios between the solid and liquid phases as large as five. For most engineering metals and alloys, however, the density change during solidification (or melting) is relatively modest, with density ratios typically satisfying $\rhos/\rhol=\mathcal{O}(1)$ (e.g., $\rhos/\rhol\approx 1$--$1.08$). Although enthalpy methods have traditionally been used to model melting and solidification processes, they can also be applied to evaporation, condensation, and boiling problems. In such applications, density ratios between the liquid and vapor phases can readily reach $\mathcal{O}(3)$, making the robustness and accuracy of variable-density enthalpy formulations particularly important. The improved low-Mach enthalpy method proposed in this work is capable of accurately handling such density ratios.

\subsection{Gas--PCM interface representation}

In some of the cases considered in this study, the computational domain 
contains a gas region $\Omega^{\rm G}(t)$ along with the PCM such that
$\Omega=\Omega^{\rm PCM}(t)\cup \Omega^{\rm G}(t), \; \forall t.$ Here, $\Omega^{\rm PCM}(t) = \Omega^{\rm S}(t) \cup \Omega^{\rm M}(t) \cup \Omega^{\rm L}(t)$. The gas--PCM interface $\Gamma(t)$ is represented by the zero level set of the
function $\Phi(\x,t)$. The level-set function is defined as a signed-distance
function satisfying
\begin{align}
\Phi(\x,t)=
\begin{cases}
>0, & \text{in the PCM},\\
<0, & \text{in the gas},\\
=0, & \text{at the gas--PCM interface}.
\end{cases}
\label{eqn_ls_def}
\end{align}
It is advected using the generally non-divergence-free velocity field
$\u(\x,t)$ according to
\begin{equation}
\DDD{\Phi}{t}
=
\D{\Phi}{t}
+\u\cdot\grad\Phi
=0.
\end{equation}
To preserve the signed-distance property of $\Phi$, the level-set function is
reinitialized at the end of each time step. Further details regarding the
level-set advection and reinitialization procedures can be found in our
previous work~\cite{nangia2019robust}.

We define the Heaviside function as
$H(\x,t)=H(\Phi(\x,t))$ to distinguish between the gas and PCM regions. In
practice, a smoothed Heaviside function is used to provide a continuous
transition of the material properties across the gas--PCM interface:
\begin{align}
H(\Phi)=
\begin{cases}
0,
& \Phi<-n_{\text{cells}}\Delta,\\[2mm]
\displaystyle
\frac{1}{2}\left[
1+\frac{\Phi}{n_{\text{cells}}\Delta}
+\frac{1}{\pi}
\sin\left(
\frac{\pi\Phi}{n_{\text{cells}}\Delta}
\right)
\right],
& |\Phi|\leq n_{\text{cells}}\Delta,\\[3mm]
1,
& \Phi>n_{\text{cells}}\Delta.
\end{cases}
\label{eqn_heaviside}
\end{align}
Here, $n_{\text{cells}}\Delta$ defines the width of the interfacial region on
each side of the gas--PCM interface. The parameter
$n_{\text{cells}}\in\mathbb{R}$ denotes the prescribed number of grid cells
over which the transition is distributed, and $\Delta$ is the grid-cell size.
A value of $n_{\text{cells}}=2$ is used in all simulations involving a gas
phase.

The thermophysical properties, including density, thermal conductivity,
viscosity, and specific heat capacity, are denoted by $\rho$, $\kappa$, $\mu$,
and $C$, respectively, and are defined continuously across the gas--PCM
interface. Using $\beta$ to denote a generic thermophysical property, the
mixture model is written as
\begin{equation}
\beta
=
\beta^{\rm G}
+
\left(\beta^{\rm S}-\beta^{\rm G}\right)H
+
\left(\beta^{\rm L}-\beta^{\rm S}\right)H\varphi,
\label{eqn_BetaInterpolation}
\end{equation}
in which $\varphi$ is the liquid fraction defined previously in
Eq.~\eqref{eqn_varphi_temp_new}.

\subsection{Mass-balance and enthalpy equations}\label{subsec_MassBalance_enthalpy_eqns}

The one-fluid model is employed to solve the mass and enthalpy equations throughout the domain,
\begin{subequations} \label{eqn_mass_enthalpy_iLM}
    \begin{alignat}{2}
        &\D{\rho}{t} + \div{(\rho \u)} = 0,\label{eqn_mass_iLM}\\
        & \D{(\rho h)}{t} + \div{(\rho \u h)}
        = \div(\kappa \nabla T). \label{eqn_enthalpy_iLM}
   \end{alignat}
\end{subequations}
Here, $\rho(\x,t)$ and $\kappa(\x,t)$ denote the spatially and temporally varying density and thermal conductivity, respectively, and $\u(\x,t)$ is the Eulerian velocity field. 



\subsection{Improved liquid fraction, specific enthalpy, and temperature relations}\label{subsec_improved_lf_h_T_relations}

The improved liquid-fraction relation as a function of temperature in the three regions of the PCM is given by (see also Eq.~\eqref{eqn_varphi_temp_new})
\begin{align}
    \varphi(T) =
\begin{cases}
0, & T<\Tsol, \\
\dfrac{T-\Tsol}{\Tliq-\Tsol}, & \Tsol\le T<\Tliq, \\
1, & T\ge \Tliq, 
\end{cases} 
\label{eqn_phi_T_improved}
\end{align}
Similarly, using the improved mushy region specific enthalpy definition in Eq.~\eqref{eqn_h_pcm_new}, the liquid-fraction--enthalpy relation in the PCM region reads as
\begin{align}
   \varphi(h) &=
\begin{cases}
0, & h<\hsol,\\
\displaystyle\frac{\rhos(h-\hsol)}{\rhol(\hliq-\hsol)-(\rhol-\rhos)(h-\hsol)},& \hsol\le h<\hliq,\\
1, & h\ge \hliq,
\end{cases} 
\label{eqn_phi_h_improved}
\end{align}
In the gas phase, $\varphi$ is assigned the arbitrary value zero; that is,
$\varphi$ is independent of $T$ and $h$ in $\Omega^{\rm G}$.

In the gas region, specific enthalpy is defined as 
\begin{equation}
    h = \cpg(T - \Tr).
\end{equation} 
The temperature--enthalpy relations in PCM and gas regions can be inverted analytically to obtain
\begin{subequations}
\begin{alignat}{2}
    \left.T(h)\right|_{\rm PCM} &=
\begin{cases}
\displaystyle \frac{h}{\cps}+\Tr, & h<\hsol,\\
\displaystyle\Tsol+\frac{\rhos(h-\hsol)}{\rhol A-\rhot(h-\hsol)},& \hsol\le h<\hliq,\\
\displaystyle\Tliq+\frac{h-\hliq}{\cpl}, & h\ge \hliq, 
\end{cases}\\
\left.T(h)\right|_{\rm gas} &= \frac{h}{\cpg} + \Tr.
\end{alignat}
\label{eqn_T_h_improved}
\end{subequations}
 The derivative of the specific enthalpy with respect to temperature, which corresponds to the effective specific heat capacity, is
\begin{subequations}
\begin{alignat}{2}
    \left.\dd{h}{T}\right|_{\rm PCM} &=
\begin{cases}
\cps, & T<\Tsol,\\
\displaystyle\frac{\rhol\rhos A}{\rho(T)^2},& \Tsol\le T<\Tliq,\\
\cpl, & T\ge \Tliq, 
\end{cases} \\
\left.\dd{h}{T}\right|_{\rm gas} &= \cpg.
\end{alignat}
\label{eqn_dhdT_improved}
\end{subequations}
Equations~\eqref{eqn_phi_T_improved}--\eqref{eqn_dhdT_improved} are used in the Newton iteration employed to solve the nonlinear enthalpy equation~\eqref{eqn_enthalpy_iLM}; details are provided in Sec.~\ref{sec_time_stepping}.

\subsection{Low-Mach velocity-divergence constraint}

We now derive the kinematic constraint on the velocity field that accounts for density-change-induced flow during phase change. Using the mass-conservation equation~\eqref{eqn_mass_enthalpy_iLM} together with the density mixture relation~\eqref{eqn_rho_mixture}, we obtain
\begin{subequations}
\begin{alignat}{2}
&\D{\rho}{t}+\grad\cdot(\rho\u)=0,  \\
\hookrightarrow& \D{\rho}{t} + \rho\grad\cdot\u + \u\cdot\grad\rho =0,  \\
\hookrightarrow& \grad\cdot\u = -\frac{1}{\rho} \left( \D{\rho}{t} + \u\cdot\grad\rho \right) = -\frac{1}{\rho}\DDD{\rho}{t}, \\
\hookrightarrow&\grad\cdot\u =  \frac{\rhos-\rhol}{\rho(\varphi)} \DDD{\varphi}{t}. \label{eqn_divu_general}
\end{alignat}
\end{subequations}

The material derivative of $\varphi$, which is nonzero only within the mushy region $\Omegam$, can be expressed in terms of the material derivative of the specific enthalpy as follows:
\begin{align}
\DDD{\varphi}{t}
=
\dd{\varphi}{h}\,\DDD{h}{t}
=
\frac{\grad\cdot(\kappa\grad T)}{\rho(\varphi)}
\dd{\varphi}{h}.
\end{align}

Using the $\varphi$--$h$ relation given by Eq.~\eqref{eqn_phi_h_improved}, the velocity-divergence constraint can be written as
\begin{subequations}
\begin{alignat}{2}
    \left.\grad\cdot\u\right|_{\rm PCM} &=
\begin{cases}
0,&  h<\hsol,\\
\displaystyle \frac{(\rhos-\rhol)}{\rhol\rhos A\DeltaT}\left[\grad\cdot(\kappa\grad T)\right],&  \hsol\le h<\hliq,\\
0,&  h\ge\hliq,
\end{cases}\\
\left.\grad\cdot\u\right|_{\rm gas} &= 0.
\end{alignat}
\label{eqn_divu_casewise}
\end{subequations}
From Eqs.~\eqref{eq:hsol_hliq_general}, we have the identity $A\Delta T=\hliq-\hsol$.
The velocity-divergence constraint in the mushy region, Eq.~\eqref{eqn_divu_casewise}, is equivalent to the similarity ODE~\eqref{eqn_VM_eta} derived in Sec.~\ref{subsec_unmatched_density_mushy_stefan}. Although the velocity-divergence condition presented in our previous works~\cite{thirumalaisamy2023low,thirumalaisamy2025consistent} appears different from Eq.~\eqref{eqn_divu_casewise}, the two formulations are mathematically equivalent. In the earlier derivation, the coefficient multiplying the term $\div(\kappa\grad T)$ was left in an unsimplified form, making it appear spatially varying. A straightforward algebraic simplification shows that this coefficient is, in fact, constant, yielding the compact expression given in Eq.~\eqref{eqn_divu_casewise}.   

\noindent\textbf{\underline{Sharp-interface limit of the low-Mach constraint}}:
In our previous works~\cite{thirumalaisamy2023low,thirumalaisamy2025consistent}, the sharp-interface limit of the low-Mach velocity-divergence constraint was not investigated. The derivation presented below demonstrates that Eq.~\eqref{eqn_divu_casewise} reduces to the classical sharp-interface constraint in the limit of a vanishing mushy region thickness.

For the classical sharp-interface Stefan problem with a density jump across the solid--liquid interface, the velocity-divergence constraint reads as~\cite{tryggvason2011direct},
\begin{equation}
\grad\cdot\u = \left(1 - \frac{\rhos}{\rhol}\right) u_n \delta(\x-\s).\label{eqn_divu_sharp}
\end{equation}
Here, $u_n=\u^{*}\cdot\n$ denotes the normal velocity of the interface, $\delta$ is the Dirac-delta function, and $\x=\s$ denotes the interface location. At first glance, Eq.~\eqref{eqn_divu_sharp} appears very different from Eq.~\eqref{eqn_divu_casewise}. However, applying the Rankine--Hugoniot jump conditions to the mass and enthalpy equations~\eqref{eqn_mass_enthalpy_iLM} across an infinitesimally thin mushy region moving with velocity $\u^{*}$ yields
\begin{align}
&\left[\kl \grad T - \ks \grad T\right] \cdot \n = \rhos(\hsol - \hliq) u_n, \nonumber \\
\hookrightarrow \; & \div (\kappa \grad T) = \rhos(\hsol - \hliq) u_n \delta(\x-\s).
\label{eqn_heat_flux_jump}
\end{align}
Substituting Eq.~\eqref{eqn_heat_flux_jump} into Eq.~\eqref{eqn_divu_casewise} recovers Eq.~\eqref{eqn_divu_sharp}. 

\subsection{Momentum equation}

The low-Mach constraint in Eq.~\eqref{eqn_divu_casewise} is solved together with the momentum equation
\begin{equation}
\D{(\rho\u)}{t}+\grad\cdot(\rho\u\otimes\u)=-\grad p+\grad\cdot\left[\mu\left(\grad\u+(\grad\u)^\intercal\right)\right] - B_d\u + \rho\vec{g} + {\vec{f}_{\rm st}},
\label{eqn_momentum}
\end{equation}
to obtain the velocity field $\u(\x,t)$ and pressure field $p(\x,t)$. Here, $\mu(\x,t)$ denotes the spatially and temporally varying dynamic viscosity, $\rho(\vec{x},t)\vec{g}$ represents the gravitational body force, with $\vec{g}$ denoting the gravitational acceleration, and ${\vec{f}_{\rm st}}$ represents the surface tension force acting on the gas-liquid interface, which is present in some of the cases considered in Sec.~\ref{sec_results}. Details on the implementation of the surface tension force can be found in~\cite{thirumalaisamy2025consistent}. $B_d$ is the Carman--Kozeny drag coefficient defined as
\begin{equation}
B_d=C_d\frac{(\varphi^{\rm S})^2}{(1-\varphi^{\rm S})^3+\epsilon}.
\label{eqn_carman_kozeny}
\end{equation}
In Eq.~\eqref{eqn_carman_kozeny}, $\epsilon$ is a small regularization parameter introduced to avoid division by zero; throughout this work, $\epsilon=10^{-3}$. The drag term suppresses fluid motion within the solid region. Here, $\varphi^{\rm S}$ denotes the solid fraction, defined as $\varphi^{\rm S}=H(1-\varphi)$. The coefficient $C_d$ is assigned a large value, $C_d=\rho/\Delta t$, where $\Delta t$ is the time-step size.
The Carman--Kozeny drag force strongly influences the pressure distribution across the mushy region and acts analogously to the Darcy--Brinkman resistance used to model flow through porous media~\cite{durlofsky1987analysis}.

\subsection{Spatial discretization}

\begin{figure}
  \centering
  \subfigure[A 2D staggered Cartesian grid]{
    \includegraphics[scale = 0.45]{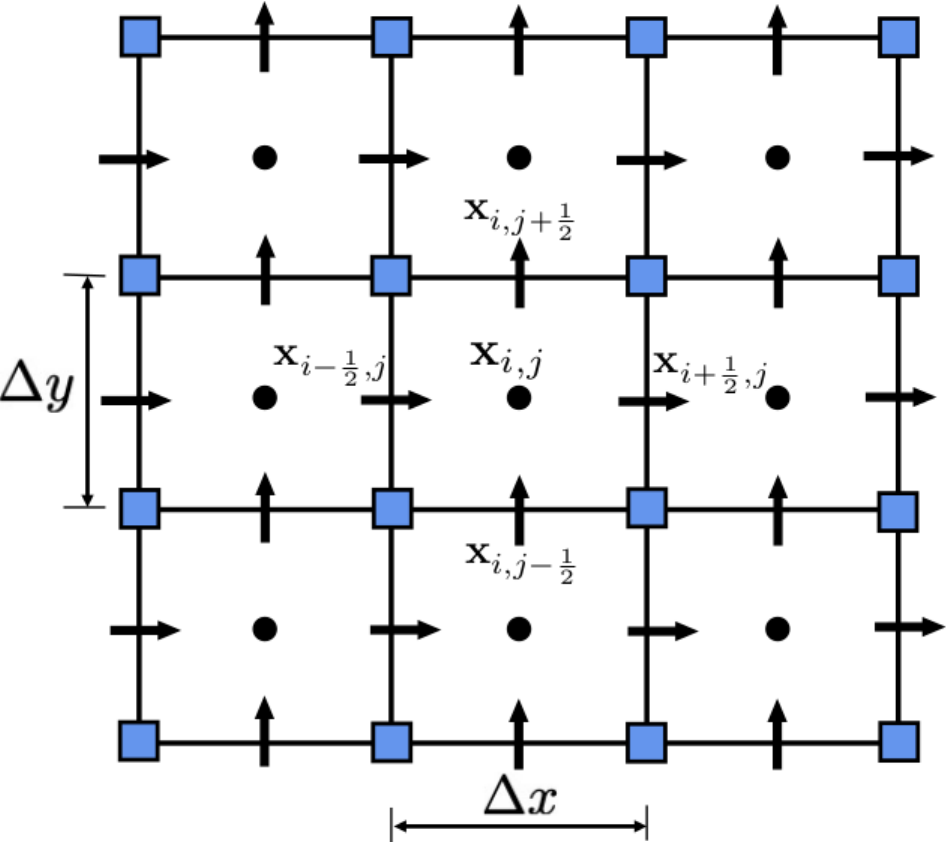} 
    \label{fig_discretized_staggered_grid}
  }
   \subfigure[A single Cartesian grid cell]{
    \includegraphics[scale = 0.35]{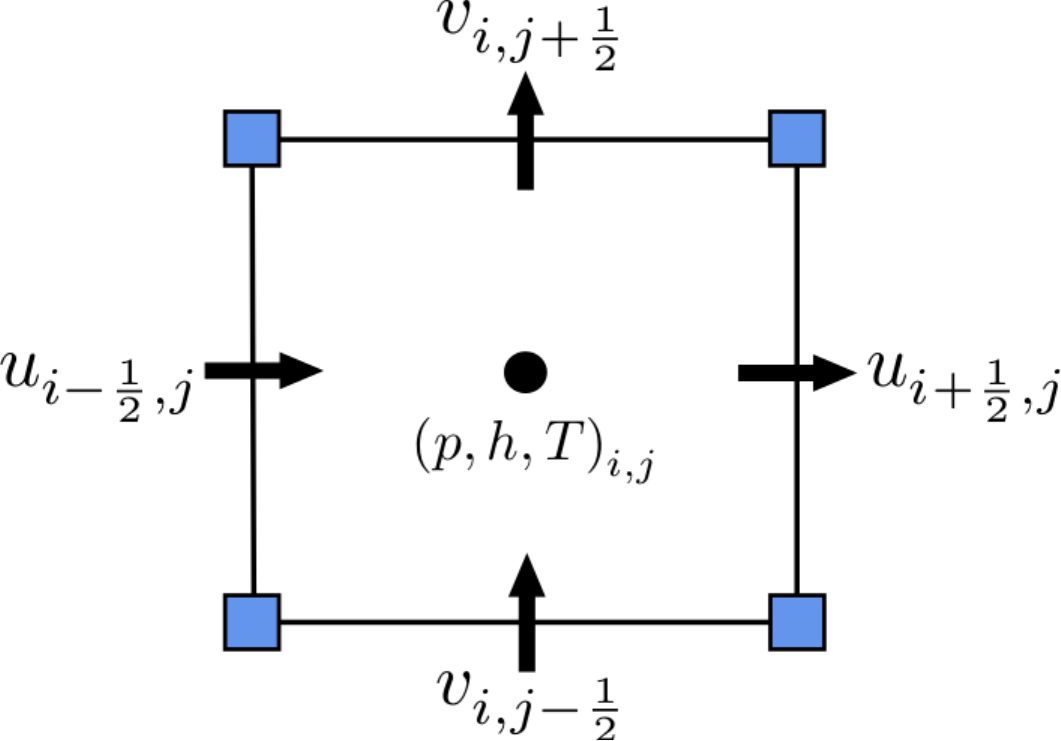}
    \label{fig_single_cell}
  }
  \caption{Schematic of the 2D staggered Cartesian grid. \subref{fig_discretized_staggered_grid} shows the coordinate system for the staggered grid. \subref{fig_single_cell} shows a single grid cell with velocity components $u$ and $v$ approximated at the cell faces (${\bf{\rightarrow}}$) and scalar variables such as pressure $p$, specific enthalpy $h$ and temperature $T$ approximated at the cell center ($\bullet$).}
\label{fig_cfd_grid}
\end{figure}

The governing equations are discretized on a staggered Cartesian grid. Scalar quantities, including pressure $p$, temperature $T$, specific enthalpy $h$, density $\rho$, liquid fraction $\varphi$, and the material properties $\kappa$ and $C$, are stored at cell centers, whereas velocity components are stored at cell faces. Specifically, the $x$- and $y$-components of velocity are stored on the $x$- and $y$-faces, respectively. For a two-dimensional Cartesian grid, the cell-center location is given by
$\x_{i,j} = \left[\left(i+\frac{1}{2}\right)\dx,\left(j+\frac{1}{2}\right)\dy\right]$, where $i=0,\ldots,N_x-1$ and $j=0,\ldots,N_y-1$ for a uniform grid consisting of $N_x\times N_y$ cells. The $x$- and $y$-velocity components are stored at $\x_{i+\frac{1}{2},j}$ and $\x_{i,j+\frac{1}{2}}$, respectively; see Fig.~\ref{fig_cfd_grid}. The spatial discretizations of the key continuous operators are as follows:

\begin{itemize}

\item The discrete divergence operator is evaluated at cell centers as
\begin{equation}
(\div\u)=\frac{u_{i+\frac{1}{2},j}-u_{i-\frac{1}{2},j}}{\dx}+\frac{v_{i,j+\frac{1}{2}}-v_{i,j-\frac{1}{2}}}{\dy}.
\label{eqn_discrete_divergence}
\end{equation}

\item The continuous form of the divergence of the viscous strain rate tensor, which couples the velocity components through spatially variable shear viscosity, is given by
\begin{equation}
\label{eq_visc_cont}
\nabla \cdot \left[\mu \left(\nabla \u + (\nabla \u)^\intercal\right) \right] = \left[
\begin{array}{c}
 (\L_{\mu} \u)^x_{i-\half,j} \\
 (\L_{\mu} \u)^y_{i,j-\half}  \\
\end{array}
\right] = 
\left[
\begin{array}{c}
 2 \D{}{x}\left(\mu \D{u}{x}\right) + \D{}{y}\left(\mu\D{u}{y}+\mu\D{v}{x}\right) \\
 2 \D{}{y}\left(\mu \D{v}{y}\right) + \D{}{x}\left(\mu\D{v}{x}+\mu\D{u}{y}\right) \\
\end{array}
\right].
\end{equation}
The shear viscosity operator is discretized using standard second-order, centered finite differences

\begin{subequations}
\begin{alignat}{2}
 (\L_{\mu} \u)^x_{i-\half,j} &= \frac{2}{\dx}\left[\mu_{i,j}\frac{u_{i+\half,j} - u_{i-\half,j}}{\dx} -
					        \mu_{i-1,j}\frac{u_{i-\half,j} - u_{i-\3half,j}}{\dx}\right] \nonumber \\ 
                    &+ \frac{1}{\dy}\left[\mu_{i-\half, j+\half}\frac{u_{i-\half,j+1} - u_{i-\half,j}}{\dy} - 
					         \mu_{i-\half, j-\half}\frac{u_{i-\half,j} - u_{i-\half,j-1}}{\dy}\right] \nonumber \\
	            &+ \frac{1}{\dy}\left[\mu_{i-\half, j+\half}\frac{v_{i,j+\half} - v_{i-1,j+\half}}{\dx} - 
					         \mu_{i-\half, j-\half}\frac{v_{i,j-\half} - v_{i-1,j-\half}}{\dx}\right] \label{eq_viscx_fd} \\
 (\L_{\mu} \u)^y_{i,j-\half} &= \frac{2}{\dy}\left[\mu_{i,j}\frac{v_{i,j+\half} - v_{i,j-\half}}{\dy} -
					        \mu_{i,j-1}\frac{v_{i,j-\half} - v_{i,j-\3half}}{\dy}\right] \nonumber \\ 
                    &+ \frac{1}{\dx}\left[\mu_{i+\half, j-\half}\frac{v_{i+1,j-\half} - v_{i,j-\half}}{\dx} - 
					         \mu_{i-\half, j-\half}\frac{v_{i,j-\half} - v_{i-1,j-\half}}{\dx}\right] \nonumber \\
	            &+ \frac{1}{\dx}\left[\mu_{i+\half, j-\half}\frac{u_{i+\half,j} - u_{i+\half,j-1}}{\dy} - 
					         \mu_{i-\half, j-\half}\frac{u_{i-\half,j} - u_{i-\half,j-1}}{\dy}\right] \label{eq_viscy_fd},
\end{alignat}
\end{subequations}
in which the shear viscosity is required at both cell centers and nodes of the staggered grid (i.e., $ \mu_{i\pm\half, j\pm\half}$).
Node-centered quantities are obtained via interpolation by harmonically averaging the neighboring cell-centered quantities. 

\item The diffusion term appearing in the enthalpy equation, Eq.~\eqref{eqn_enthalpy_iLM}, is discretized as
\begin{align}
\div(\kappa\grad T)=&\frac{1}{\dx}\left[\kappa_{i+\frac{1}{2},j}\frac{T_{i+1,j}-T_{i,j}}{\dx}-\kappa_{i-\frac{1}{2},j}\frac{T_{i,j}-T_{i-1,j}}{\dx}\right] \nonumber \\
&+\frac{1}{\dy}\left[\kappa_{i,j+\frac{1}{2}}\frac{T_{i,j+1}-T_{i,j}}{\dy}-\kappa_{i,j-\frac{1}{2}}\frac{T_{i,j}-T_{i,j-1}}{\dy}\right].
\label{eqn_discrete_diffusion}
\end{align}

\item The gradient of a cell-centered scalar quantity, such as pressure, is evaluated at cell faces according to
\begin{equation}
(\grad p)_{i+\frac{1}{2},j}^{x}
=
\frac{p_{i+1,j}-p_{i,j}}{\dx},
\qquad
(\grad p)_{i,j+\frac{1}{2}}^{y} =
\frac{p_{i,j+1}-p_{i,j}}{\dy}.
\label{eqn_discrete_gradient}
\end{equation}

\end{itemize}

Thermal conductivity is interpolated from cell centers to cell faces using harmonic averaging. Convective fluxes are discretized using the third-order accurate cubic upwind interpolation (CUI) scheme. The CUI scheme satisfies both the convection-boundedness criterion (CBC) and the total-variation-diminishing (TVD) property~\cite{moukalled2016finite}. In monotone regions, the scheme achieves third-order spatial accuracy, whereas in non-monotone regions it reverts to a first-order upwind discretization to preserve stability and boundedness~\cite{nangia2019robust}.

\subsection{Solution algorithm and time stepping scheme} \label{sec_time_stepping}

We summarize the complete solution algorithm employed in the improved low-Mach enthalpy formulation to advance the solution from time $t^n$ to $t^{n+1}=t^n+\dt$. Time advancement is performed using a fixed-point iteration indexed by $k=0,\ldots,p-1$, where $p$ denotes the prescribed number of fixed-point iterations. Within each fixed-point iteration, the nonlinear enthalpy equation is solved using Newton iterations indexed by $m=0,\ldots,\qmax-1$. Unless otherwise specified, $p=2$ and $\qmax=5$ are used throughout this work. Consequently, the nonlinear enthalpy solve may require up to $p\times\qmax$ iterations per timestep. At the beginning of each timestep, the iterative solution is initialized by setting $\theta^{n+1,0}=\theta^n$ for any generic variable $\theta$. The main steps of the algorithm within each timestep are summarized below:

\begin{enumerate}
\item First, advect the level-set function ${\Phi}^n$ using non-divergence-free velocity $\u$ to obtain $\Phi^{n+1, k+1}$,
\begin{equation}
    \frac{\Phi^{n+1, k+1} - {\Phi}^n}{\Delta t} + \left(\u \cdot \nabla \Phi \right)^{n+1,k} = 0.
\end{equation}
An updated value of Heaviside $H^{n+1,k+1}$ is then computed from $\Phi^{n+1,k+1}$ using Eq.~\ref{eqn_heaviside}.

\item Next, the mass-conservation equation is solved at both cell faces and cell centers. A discrete mass flux $\vec{m}_\rho$ is computed to obtain the updated density field $\breve{\rho}^{n+1,k+1}$,
\begin{equation}
\frac{\breve{\rho}^{n+1,k+1}-\rho^n}{\dt}+\grad\cdot\vec{m}_\rho^{n+1,k}=0.
\label{eqn_discrete_mass_algorithm}
\end{equation}
The resulting mass flux is used in the convective operator of the momentum equation, $\div(\rho\u\otimes\u)=\div(\vec{m}_\rho\otimes\u)$, and in the advective operator of the enthalpy equation, $\div(\rho\u h)=\div(\vec{m}_\rho h)$. Eq.~\eqref{eqn_discrete_mass_algorithm} is integrated using a second-order Runge--Kutta scheme. The solution of Eq.~\eqref{eqn_discrete_mass_algorithm}, together with the use of the same discrete mass flux $\vec{m}_\rho$ in the momentum and enthalpy equations, ensures consistent mass, momentum, and enthalpy transport. Further details on the consistency requirements for mass--momentum--enthalpy transport can be found in our previous work~\cite{thirumalaisamy2025consistent}.

\item Using the updated density field $\breve{\rho}^{n+1,k+1}$ and mass flux $\vec{m}_\rho^{n+1,k}$, the nonlinear enthalpy equation is solved to update the enthalpy $h$, temperature $T$, and liquid fraction $\varphi$. The discrete enthalpy equation reads
\begin{equation}
\frac{\breve{\rho}^{n+1,k+1}h^{n+1,k+1}-\rho^n h^n}{\dt}+\grad\cdot\left(\V{m}_{\rho}^{n+1,k} h^{n+1,k}\right)=\left[\grad\cdot(\kappa\grad T)\right]^{n+1,k+1}.
\label{eqn_discrete_enthalpy}
\end{equation}
The cyclic nonlinear dependence of $h$, $\varphi$, $\rho$, and $\kappa$ on $T$ is resolved using Newton's method. At Newton iteration $m$, the specific enthalpy is linearized as
\begin{equation}
h^{n+1,k+1,m+1}=h^{n+1,k+1,m}+\left(\dd{h}{T}\right)^{n+1,k+1,m}\left(T^{n+1,k+1,m+1}-T^{n+1,k+1,m}\right),
\label{eqn_newton_h_linearization}
\end{equation}
in which $\dd{h}{T}$ is given by Eq.~\eqref{eqn_dhdT_improved}. After solving for $T^{n+1,k+1,m+1}$, the enthalpy field is updated using Eq.~\eqref{eqn_newton_h_linearization}. The liquid fraction and temperature are then synchronized with the updated enthalpy field using the $\varphi$--$h$ and $T$--$h$ relations given by Eqs.~\eqref{eqn_phi_h_improved} and~\eqref{eqn_T_h_improved}, respectively. The thermophysical properties $(\kappa,\mu,\rho,C)$ are subsequently updated using the latest value of $\varphi$ and $H$.

The Newton iteration is considered converged when the relative change in liquid fraction,
\begin{equation}
\Delta_\varphi = \frac{||\varphi^{n+1,k+1,m+1} - \varphi^{n+1,k+1,m}||_2}{1 + || \varphi^{n+1,k+1,m}||_2},
\end{equation}
satisfies $\Delta\varphi \le 10^{-8}$, or when the maximum number of Newton iterations, $\qmax=5$, has been reached.

\item Finally, the momentum equation and the low-Mach velocity-divergence constraint are solved simultaneously for $\u^{n+1,k+1}$ and $p^{n+\frac{1}{2},k+1}$:
\begin{subequations} 
\label{eqn_discrete_momentum} 
\begin{alignat}{2} &\frac{\breve{\rho}^{n+1,k+1}\u^{n+1,k+1}-\rho^n\u^n}{\dt}+\grad\cdot\left(\vec{m}_{\rho}^{n+1,k}\otimes\u^{n+1,k}\right)=-\grad p^{n+\frac{1}{2},k+1} \nonumber \\ &+\grad\cdot\left[\mu\left(\grad\u+ (\grad\u)^\intercal\right)\right]^{n+\frac{1}{2},k+1} -B_d^{n+1,k+1}\u^{n+1,k+1} + \breve{\rho}^{n+1,k+1}\vec{g} + {\vec{f}_{\rm st}}^{n+\frac{1}{2},k+1}, \\ 
& \div (\u^{n+1,k+1}) = -\frac{1}{\rho}\DDD{\rho}{t}. 
\end{alignat} 
\end{subequations}
The viscous term is discretized using the Crank-Nicolson scheme. The Carman--Kozeny term renders the resulting linear system stiff. To efficiently solve Eqs.~\eqref{eqn_discrete_momentum}, we employ a projection-method-based preconditioner that consistently accounts for the stiff drag force term within a flexible GMRES (FGMRES) framework. Details of the linear solver and preconditioning strategy can be found in our previous work~\cite{thirumalaisamy2023pre}.
\end{enumerate}

\section{Software implementation}\label{sec_software_implementation}

The improved low-Mach enthalpy method presented in this work has been implemented within the IBAMR library~\cite{IBAMR-web-page}, an open-source C++ framework for the simulation of computational fluid dynamics (CFD) and fluid--structure interaction problems. IBAMR provides adaptive mesh refinement (AMR) and Cartesian grid management through the SAMRAI framework~\cite{HornungKohn02,samrai-web-page}, while scalable linear algebra and iterative solvers are provided by the PETSc library~\cite{petsc-user-ref,petsc-web-page}. Simulation results are post-processed using in-house MATLAB scripts and the open-source VisIt visualization software~\cite{childs2012visit}.

To facilitate reproducibility, all IBAMR simulations presented in this work are publicly available at \url{https://github.com/IBAMR/IBAMR}. The MATLAB implementation used to compute the similarity solutions of the mushy Stefan problem is available at \url{https://github.com/amneetb/MushyStefan}.

\section{Results and discussion}\label{sec_results}

In this section, we verify the improved low-Mach enthalpy method against the
one-dimensional mushy Stefan solidification problem described in
Sec.~\ref{sec_mushy_stefan}. We then apply the verified method to the
simulation of metal-manufacturing defects. First, the numerical results
obtained with the improved formulation are compared directly with the
corresponding mushy Stefan similarity solutions for the cases considered in
Sec.~\ref{subsec_sharpvsmushy_results}. Next, we assess the accuracy and
spatiotemporal convergence properties of the method. Finally, we evaluate its
robustness for practical manufacturing problems involving substantially
larger density contrasts between the solid and liquid phases.

For the first two parts of this section, a quasi-one-dimensional computational domain is employed to simulate the solidification problem. The phase-change front propagates in the positive $x$-direction. At the left boundary, $x=0$, zero-velocity and fixed-temperature boundary conditions are imposed, i.e., $\u=\mathbf{0}$ and $T=\T0$. At the right boundary, $x=L$, an outflow (zero-pressure) boundary condition is prescribed together with the fixed temperature $T=\Ti$. In the transverse ($y$) direction, periodic boundary conditions are used. The aluminum-based thermophysical properties of the phase-change material are summarized in Table~\ref{tab_thermophys_properties_Al}. For the quasi-one-dimensional problems, shear stresses vanish identically.
Accordingly, the viscosity is set to zero in all three regions of the PCM, thereby
eliminating any stabilizing influence of shear viscosity on the numerical
discretization.

\subsection{Verification of the improved low-Mach enthalpy method against the mushy Stefan solution}\label{subsec_compare_EM_with_MushyStefan}

Here, we compare the improved low-Mach enthalpy method against the mushy Stefan similarity solution for different values of the phase-change temperature interval, $\DeltaT$. The objective is to verify that, for a prescribed value of $\DeltaT$, the numerical solution obtained with the improved enthalpy method converges to the corresponding mushy Stefan solution. The domain is initially liquid at temperature $\Ti$. At $t=0^{+}$, the temperature at the left boundary is suddenly reduced to $\T0<\Tsol$, thereby initiating solidification from the left boundary. As a result, the solid region appears immediately in the domain, avoiding numerical difficulties associated with phase-front nucleation. Simulations are performed on the domain $\Omega \in [0,1] \times [0,0.05]$ for matched and low-density-ratio cases and $\Omega \in [0,0.25] \times [0,0.0125]$ for the large-density-ratio case, using a grid of size $N_x\times N_y=1280\times64$ for all cases. 


\subsubsection{Equal-density case}\label{subsec_ratio1DiffDeltaT}

\begin{figure}[]
	\centering{\includegraphics[scale=0.29]{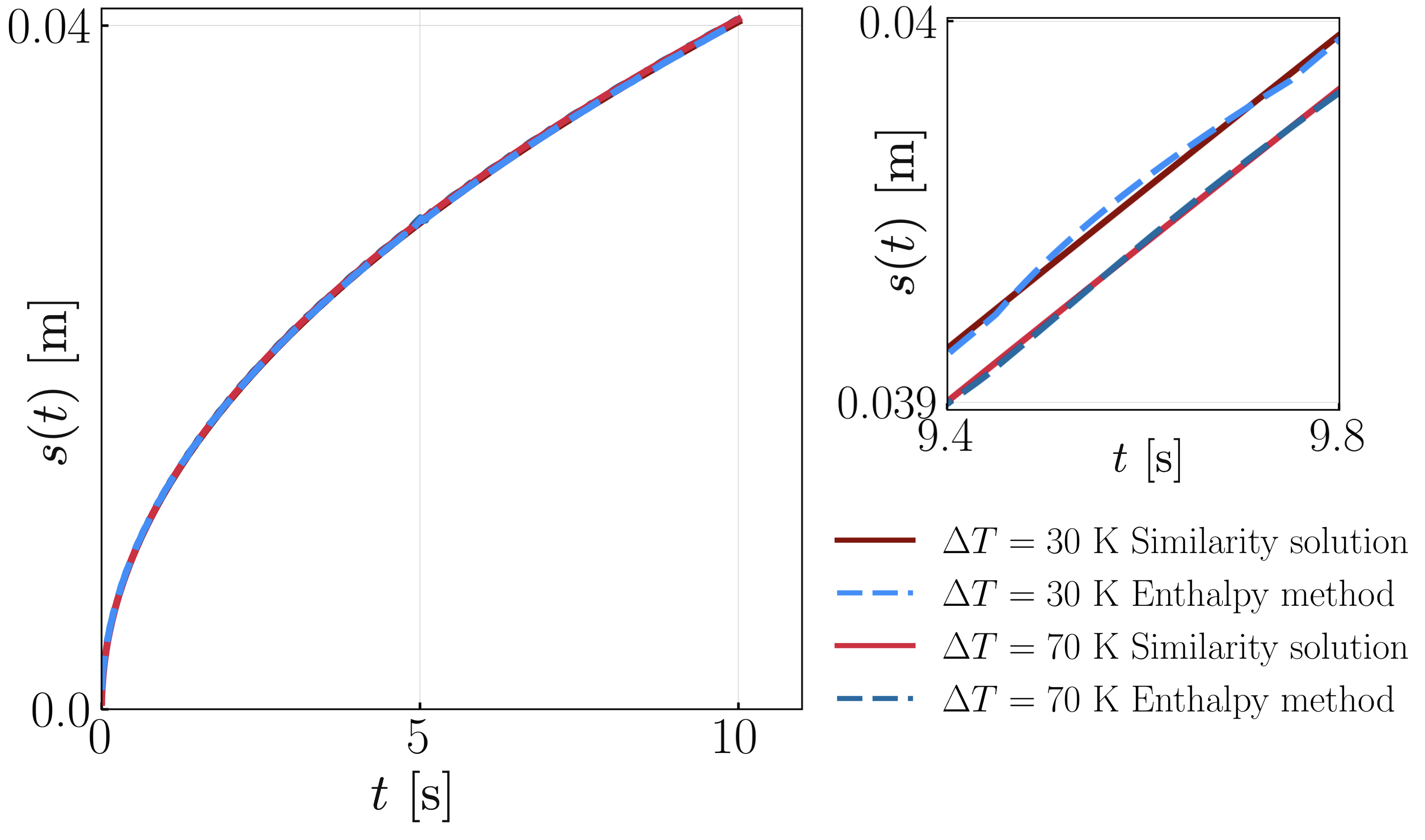}}
\caption{Comparison of the phase-front location, $x=\svarphit$, predicted by the improved low-Mach enthalpy method and the corresponding mushy Stefan similarity solution for the equal-density case ($\rhol/\rhos=1$) for different values of the phase-change interval $\DeltaT$. }\label{fig_Ratio1DiffDeltaT}
\end{figure}

This case corresponds to the equal-density configuration considered in Sec.~\ref{subsec_no_volume_change}, with $\rhos=\rhol=2700$ kg/m$^3$. The left boundary temperature is $\T0=298.6$ K, while the right boundary temperature and initial liquid temperature are both set to $\Ti=978.6$ K. Fig.~\ref{fig_Ratio1DiffDeltaT} compares the phase-front location predicted by the improved enthalpy method with the corresponding equal-density mushy Stefan solution for $\DeltaT=30$ K and $70$ K.

A constant timestep size of $\Delta t=1\times10^{-4}$ s is used for the simulations. For both values of $\DeltaT$, the numerically predicted phase-front location remains in close agreement with the mushy Stefan solution. These results confirm that, for the equal-density case, the enthalpy method is largely insensitive to the choice of the phase-change interval $\DeltaT$.

\subsubsection{Low-density-ratio volume-expansion case}\label{subsec_ratio2DiffDeltaT}

\begin{figure}[]
	\centering{\includegraphics[scale=0.29]{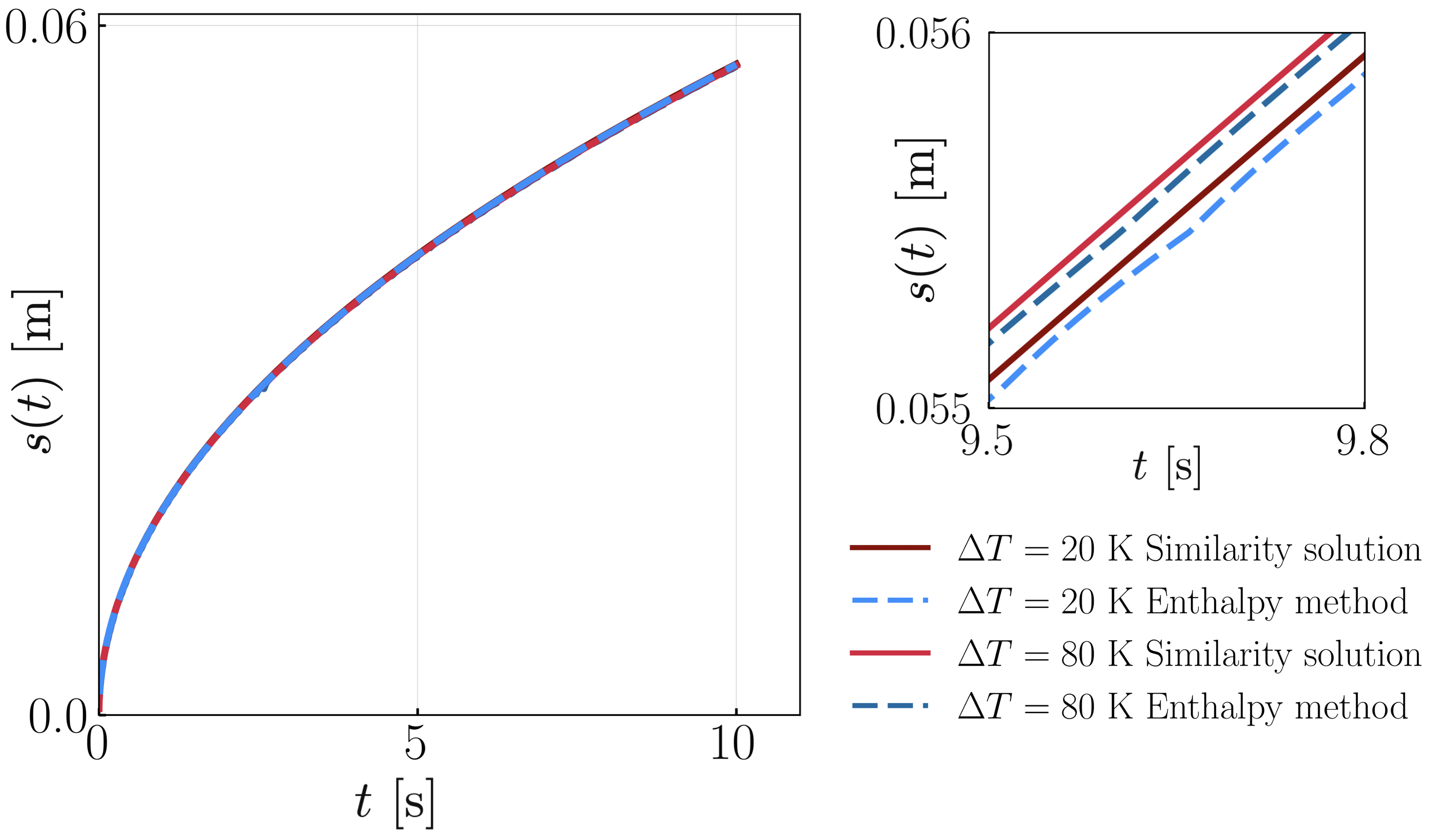}}
\caption{Comparison of the phase-front location, $x=\svarphit$, predicted by the improved low-Mach enthalpy method and the corresponding mushy Stefan similarity solution for the low-density-ratio case ($\rhol/\rhos=2$) for different values of the phase-change interval $\DeltaT$.}\label{fig_Ratio2DiffDeltaT}
\end{figure}

This case corresponds to Sec.~\ref{subsec_Low_volume_expansion_case}, with $\rhos=1350$ kg/m$^3$ and $\rhol=2700$ kg/m$^3$. Numerical solutions obtained using the improved enthalpy method with $\DeltaT=20$ K and $80$ K are compared against the corresponding mushy Stefan similarity solution. The simulations are performed using a constant timestep of $\Delta t=1\times10^{-4}$ s. Fig.~\ref{fig_Ratio2DiffDeltaT} shows that the numerically predicted phase-front location remains in close agreement with the mushy Stefan solution for both values of $\DeltaT$. Consistent with the observations of Sec.~\ref{subsec_Low_volume_expansion_case}, the low-density-ratio case exhibits only a weak dependence on the phase-change interval, and excellent agreement is obtained for both small and large values of $\DeltaT$.

\subsubsection{Large-density-ratio volume-expansion
case}\label{subsec_ratio540DiffDeltaT}

\begin{figure}[]
	\centering{\includegraphics[scale=0.30]{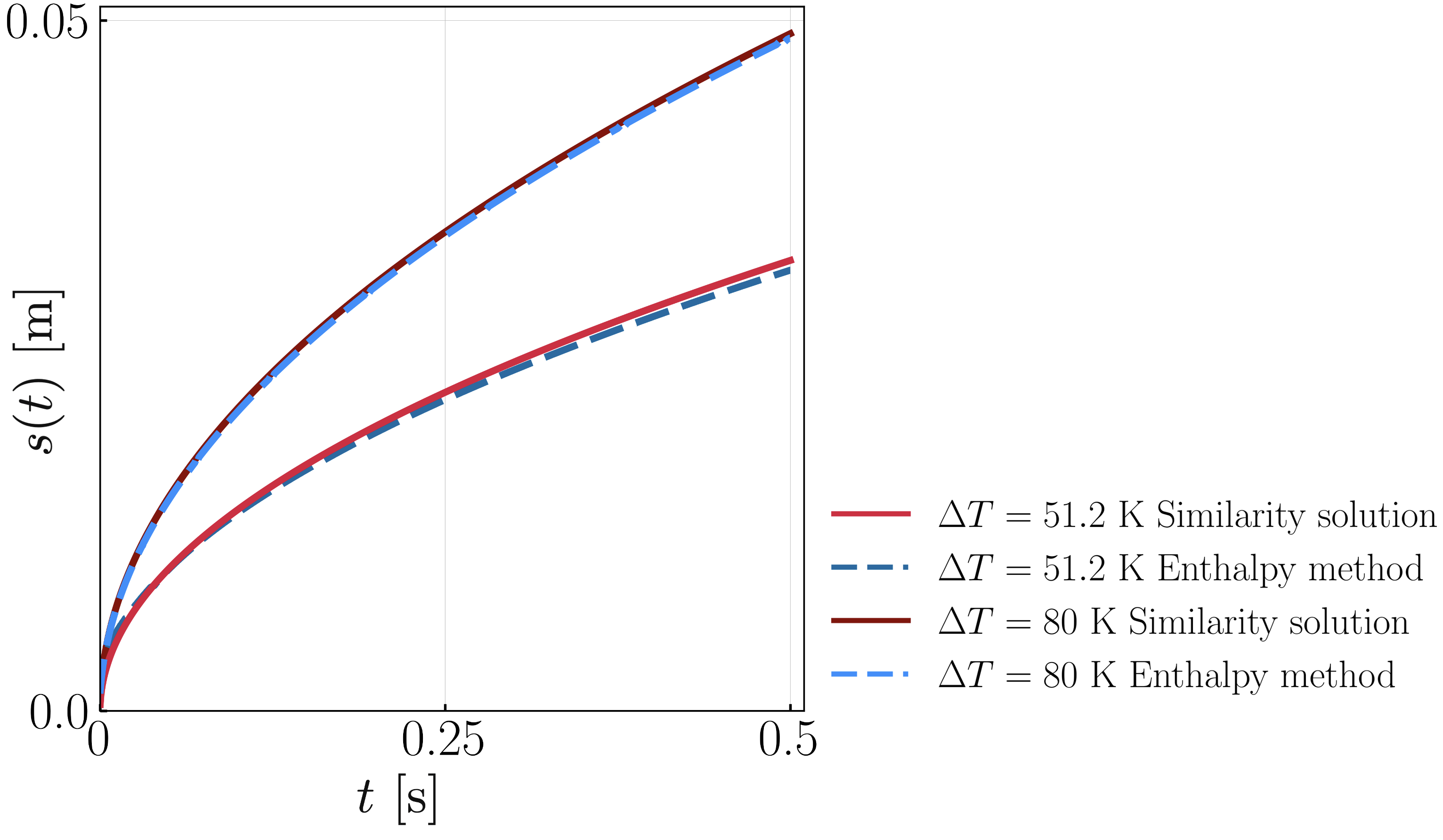}}
\caption{Comparison of the phase-front location, $x=\svarphit$, predicted by the improved low-Mach enthalpy method and the corresponding mushy Stefan similarity solution for the high-density-ratio case ($\rhol/\rhos=540$) for different values of the phase-change interval $\DeltaT$. }\label{fig_Ratio540DiffDeltaT}
\end{figure}

\begin{figure}[]
	\raggedleft\includegraphics[scale=0.52,  trim={40 402 20 300},clip,keepaspectratio]{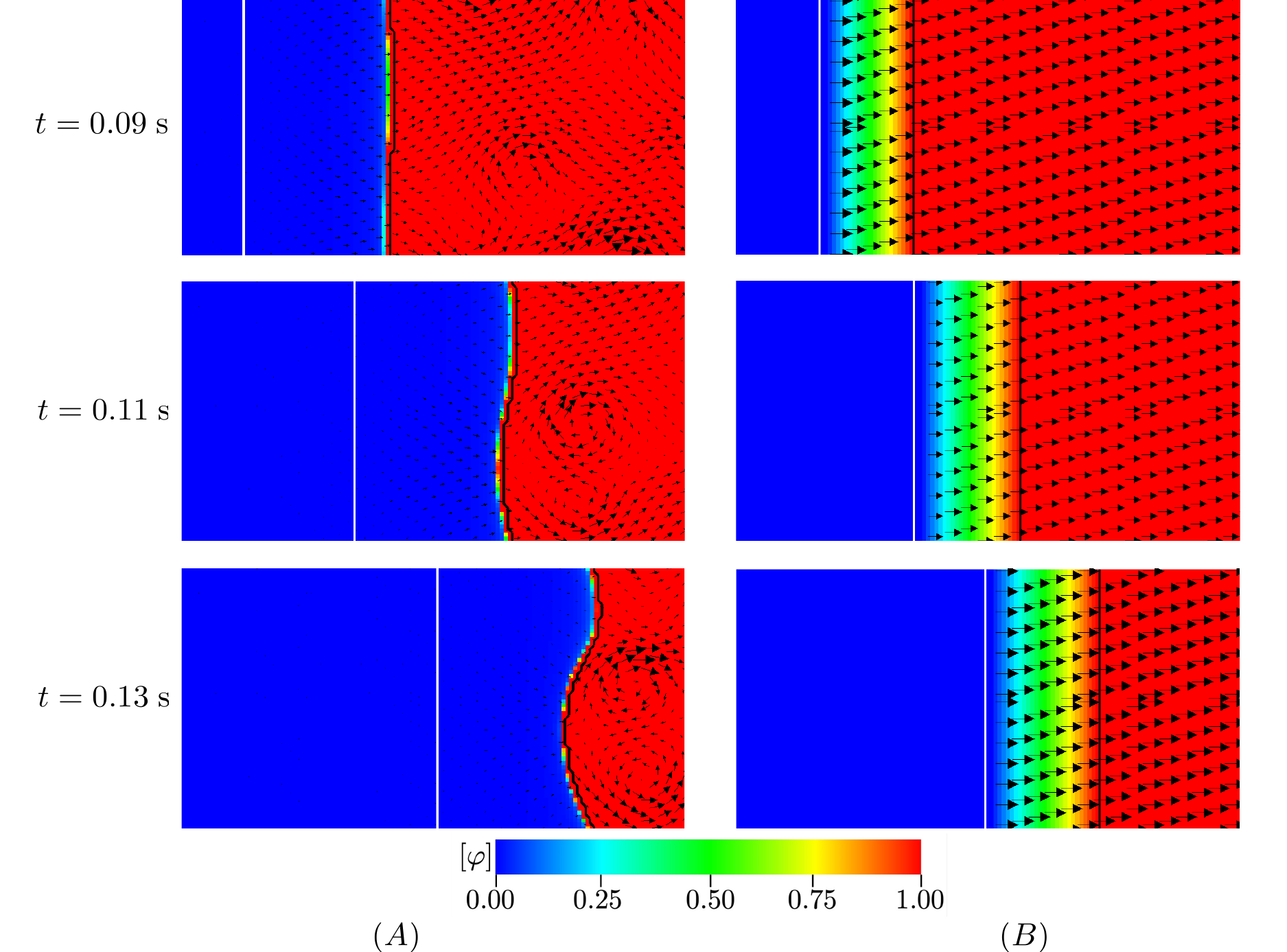}
\caption{Large-density-ratio solidification problem simulated using the (A) original and (B) improved low-Mach enthalpy methods. The time evolution of the liquid-fraction field, $\varphi$, is shown for $\rhol/\rhos=540$, illustrating the onset of interfacial distortions in the mushy region predicted by the original formulation. Blue and red denote the solid and liquid regions, respectively, with the mushy region located between them. White and black contour lines represent the solidus ($T = \Tsol$) and liquidus ($T = \Tliq$) fronts, respectively. Velocity vectors are superimposed on the solution. Simulations were performed with a phase-change interval of $\DeltaT=20$~K.}
\label{fig_Ratio540Diff}
\end{figure}

This case corresponds to Sec.~\ref{subsec_Large_volume_expansion_case}, with $\rhos=5.0$ kg/m$^3$ and $\rhol=2700$ kg/m$^3$. Fig.~\ref{fig_Ratio540DiffDeltaT} compares the phase-front location predicted by the improved enthalpy method with the corresponding large-density-ratio mushy Stefan solution for $\DeltaT=51.2$ K and $80$ K. 
The simulations are performed using a constant timestep of $\Delta t=1\times10^{-6}$ s. For both values of $\DeltaT$, the numerical solution remains in close agreement with the mushy Stefan solution presented earlier in Fig.~\ref{fig_MushyvsSharpLargeRatio}.

Fig.~\ref{fig_Ratio540Diff} compares the evolution of the liquid-fraction field, $\varphi$, obtained using the original low-Mach enthalpy method (panel A) and the improved low-Mach enthalpy method (panel B). In the original formulation, interfacial instabilities develop within the mushy region, leading to an inaccurate prediction of the phase-front location. In contrast, the improved formulation predicts a stable one-dimensional evolution of the phase front. Furthermore, the velocity in the liquid region remains spatially uniform, consistent with the theoretical prediction. In contrast, the original low-Mach enthalpy method produces unphysical vortical structures in the liquid region. These results demonstrate that the improved low-Mach enthalpy method remains accurate for phase-change problems involving large density contrasts between the two phases.

\subsection{Spatio-temporal convergence rate of the improved low-Mach enthalpy method}\label{subsec_convergence_tests}

Having verified that the improved low-Mach enthalpy method solves the one-dimensional mushy Stefan problem, we now assess its spatio-temporal convergence rate. Grid-refinement studies are performed for the three mushy Stefan problems described in Secs.~\ref{subsec_ratio1DiffDeltaT}--\ref{subsec_ratio540DiffDeltaT}. Errors are computed using two reference solutions: (i) the similarity solution of the mushy Stefan problem and (ii) the numerical solution obtained on the finest computational grid.

Although the enthalpy method is mathematically equivalent to the mushy Stefan model in one dimension, the two are not identical in the present quasi-one-dimensional setting. The difference arises because the enthalpy method solves the momentum equation together with the Carman--Kozeny drag force, whereas the mushy Stefan model does not include momentum transport. Consequently, the drag force influences the velocity and pressure fields and introduces a small modeling difference between the two formulations. In the present quasi-one-dimensional configuration, the drag force is active within the solid region and gradually diminishes across the mushy region. Our numerical experiments indicate that removing the drag force leads to small spurious velocities in the solid region, which distort the phase front at large density ratios. These distortions are distinct from the interfacial instabilities shown in Fig.~\ref{fig_Ratio540Diff} for the original low-Mach enthalpy method, which persist even when the drag coefficient is increased by several orders of magnitude. Consequently, the enthalpy method is expected to exhibit a higher observed convergence rate when the finest-grid solution is used as the reference than when compared against the mushy Stefan similarity solution.

For the equal- and low-density-ratio cases, we take $\DeltaT=10$ K because the phase-front location is only weakly sensitive to the choice of $\DeltaT$, while the mushy region remains well resolved by the computational grid. For the large-density-ratio case, we use $\DeltaT=51.2$ K, for which the mushy Stefan solution closely reproduces the sharp-interface phase-front location while maintaining adequate spatial resolution of the mushy region. The computational domain is taken as $\Omega\in[0,1]\times[0,0.05]$ for the equal- and low-density-ratio cases and $\Omega\in[0,0.25]\times[0,0.0125]$ for the large-density-ratio case. The grid-refinement study employs four meshes,
$ N_x\times N_y= \{320\times16, 640\times32, 1280\times64, 2560\times128\}$.  A timestep of $\Delta t=4\times10^{-4}$ s is used on the coarsest grid and is halved with each successive grid refinement to maintain a CFL number below 0.5 throughout the simulation. For the large-density-ratio case, the timestep on the coarsest grid is reduced to $\Delta t=4\times10^{-6}$ s. 

We report errors in both the phase-front location, $x = \svarphit$, over the entire simulation interval and the temperature field, $T(x,t)$, at a fixed time. For a generic quantity $\psi$, the $L^2$ error is computed as the root-mean-square error (RMSE),
$ \|\mathcal{E}_\psi\|_{\rm RMSE} = (\|\psi_{\rm reference}-\psi_{\rm numerical}\|_2)/\sqrt{\mathcal{S}}$, 
in which $\mathcal{S}$ denotes the number of entries in the error vector $\mathcal{E}_\psi$.

\subsubsection{Equal-density case}\label{subsubsec_match_conv_tests}

\begin{figure}[]
	\centering{\includegraphics[width=0.92\textwidth]{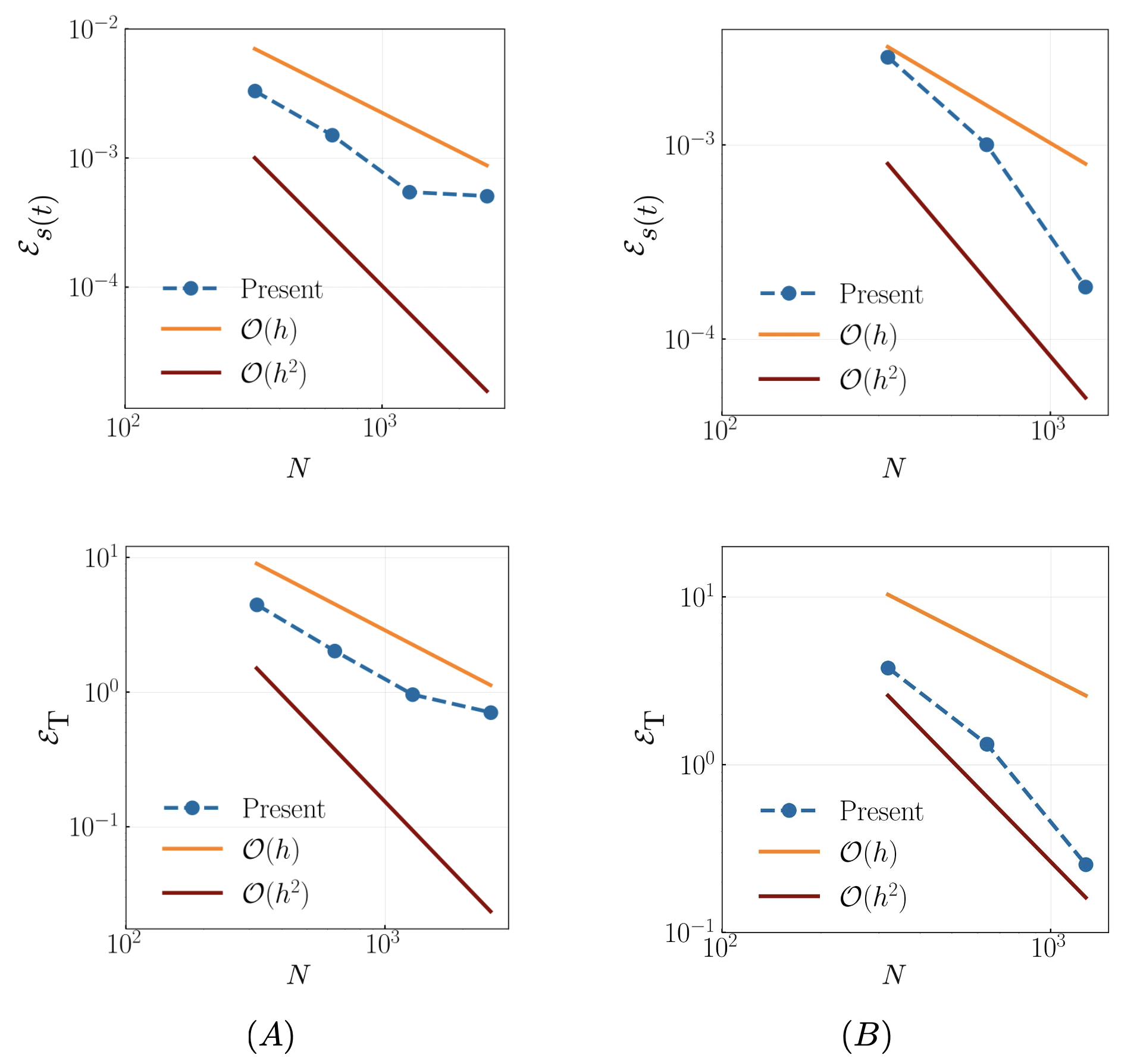}}
\caption{Convergence rates of the improved low-Mach enthalpy method for the one-dimensional Stefan problem with equal phase densities, $\rhos = \rhol = 2700$~kg/m$^3$, for $\DeltaT = 10$~K. The reference solutions are (A) the mushy Stefan similarity solution and (B) the numerical solution on the finest grid ($N_x \times N_y = 2560 \times 128$). The error in the interface position $x = \svarphit$ is computed over the time interval $0 \leq t \leq 10$~s, and the error in temperature $T(x,t)$ is computed over the entire domain ($0 \leq x \leq 1$~m) at time $t = 5$~s.}
\label{fig_ConvRateRatio1}
\end{figure}

Fig.~\ref{fig_ConvRateRatio1} illustrates the convergence behavior of the improved enthalpy method for the equal-density case, $\rhol/\rhos=1$, considered in Sec.~\ref{subsec_ratio1DiffDeltaT}. Convergence rates computed with respect to the mushy Stefan similarity solution are shown in Fig.~\ref{fig_ConvRateRatio1}(A). The method exhibits approximately first-order convergence for the temperature field and between first- and second-order convergence for the phase-front location. Convergence rates computed using the finest-grid numerical solution as the reference are shown in Fig.~\ref{fig_ConvRateRatio1}(B). In this case, both the temperature and phase-front location converge at rates between first and second order.

\subsubsection{Low-density-ratio volume-expansion case}\label{subsubsec_conv_low_density_ratio}

\begin{figure}[]
	\centering{\includegraphics[width=0.92\textwidth]{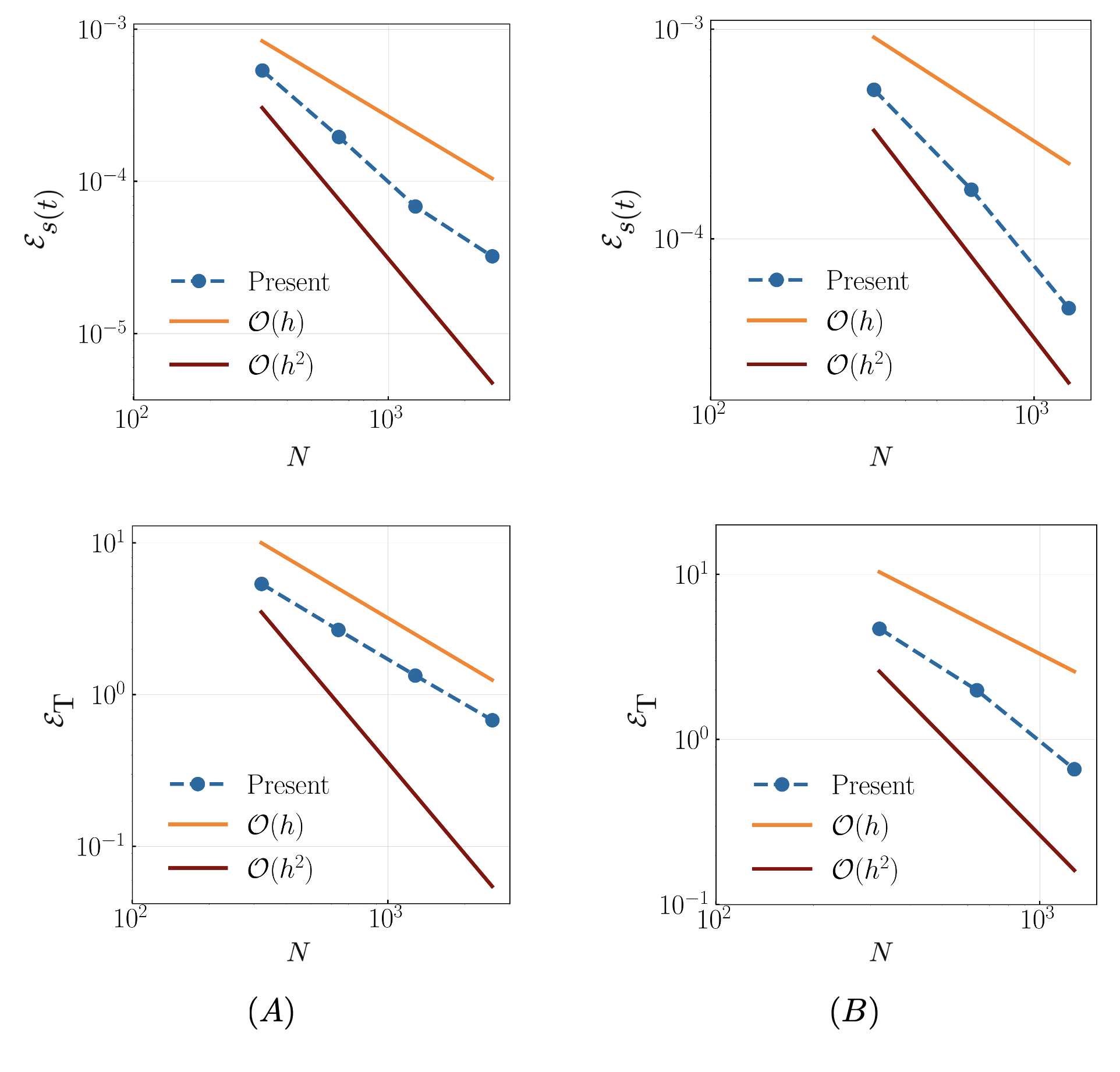}}
\caption{Convergence rates of the improved low-Mach enthalpy method for the one-dimensional Stefan problem with a low density ratio between the phases, $\rhos = 1350$~kg/m$^3$ and $\rhol = 2700$~kg/m$^3$, for $\DeltaT = 10$~K. The reference solutions are (A) the mushy Stefan similarity solution and (B) the numerical solution on the finest grid ($N_x \times N_y = 2560 \times 128$). The error in the interface position $x = \svarphit$ is computed over the time interval $0 \leq t \leq 10$~s, and the error in temperature $T(x,t)$ is computed over the entire domain ($0 \leq x \leq 1$~m) at time $t = 5$~s.}
\label{fig_ConvRateRatio2}
\end{figure}

Fig.~\ref{fig_ConvRateRatio2} presents the convergence study for the low-density-ratio case described in Sec.~\ref{subsec_ratio2DiffDeltaT}. The solid and liquid densities are $\rhos=1350$ kg/m$^3$ and $\rhol=2700$ kg/m$^3$, respectively. Convergence rates computed with respect to the mushy Stefan similarity solution are shown in Fig.~\ref{fig_ConvRateRatio2}(A). The method exhibits between first- and second-order convergence for the phase-front location and approximately first-order convergence for the temperature field. Convergence rates computed using the finest-grid numerical solution as the reference are shown in Fig.~\ref{fig_ConvRateRatio2}(B). In this case, the phase-front location exhibits nearly second-order convergence, whereas the temperature field converges at a rate between first and second order.

\subsubsection{Large-density-ratio volume-expansion case}\label{subsubsec_conv_large_ratio}

\begin{figure}[]
	\centering{\includegraphics[width=0.92\textwidth]{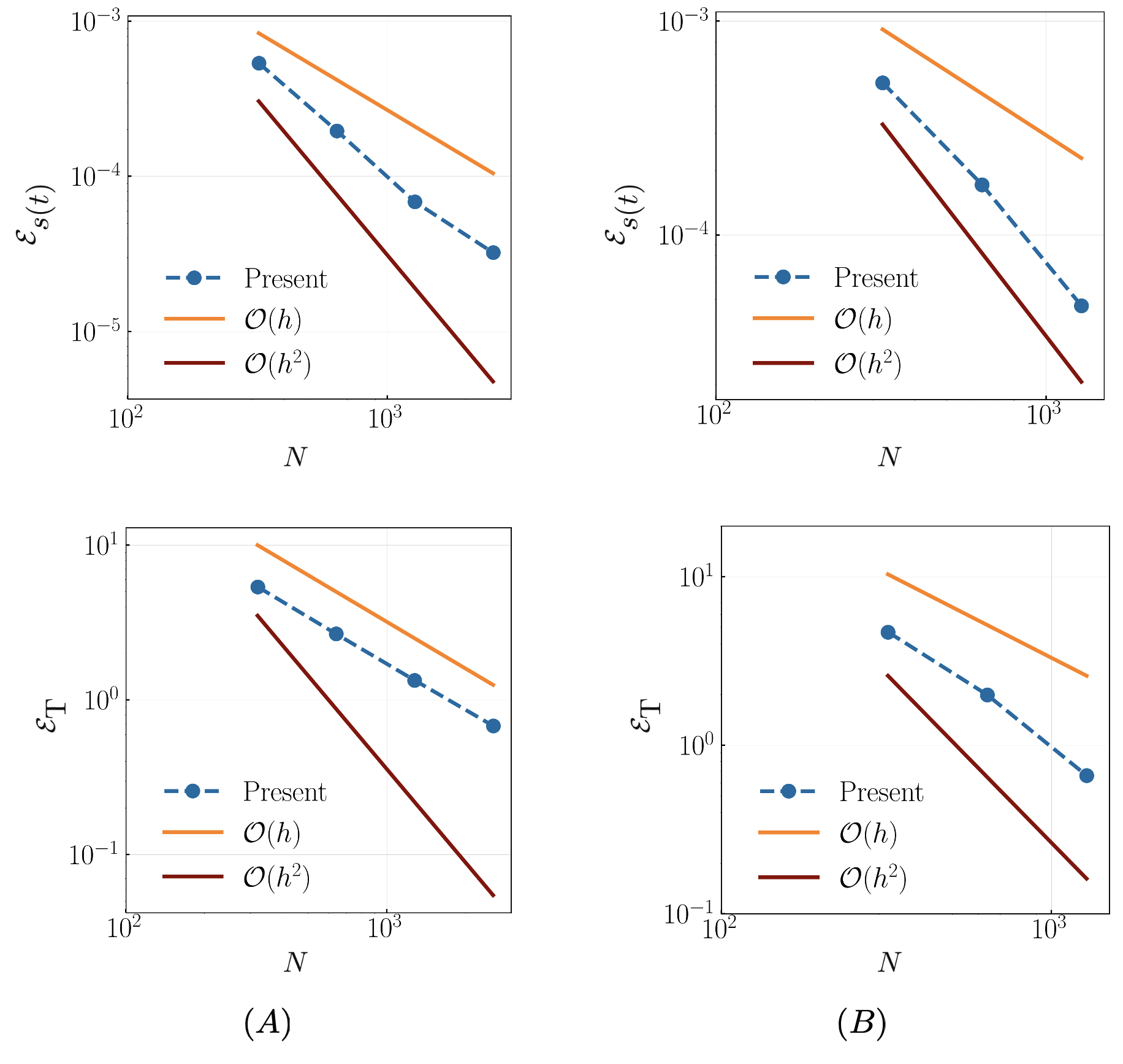}}
\caption{Spatio-temporal convergence rates of the improved low-Mach enthalpy method for the one-dimensional large-density-ratio Stefan problem with $\rhos=5$ kg/m$^3$ and $\rhol=2700$ kg/m$^3$, using $\DeltaT=51.2$ K. The reference solutions are (A) the mushy Stefan similarity solution and (B) the finest-grid numerical solution ($N_x\times N_y=2560\times128$). Errors in the phase-front location, $x=\svarphit$, are computed over the time interval $0.01\le t\le0.5$ s, while errors in the temperature field, $T(x,t)$, are computed over the entire domain ($0\le x\le0.25$ m) at $t=0.05$ s.}\label{fig_ConvRateRatio540}
\end{figure}

Fig.~\ref{fig_ConvRateRatio540} shows the spatio-temporal convergence of the improved enthalpy method for the large-density-ratio case considered in Sec.~\ref{subsec_ratio540DiffDeltaT}. The solid and liquid densities are $\rhos=5.0$ kg/m$^3$ and $\rhol=2700$ kg/m$^3$, respectively. Convergence rates computed with respect to the mushy Stefan similarity solution are shown in Fig.~\ref{fig_ConvRateRatio540}(A). The phase-front location exhibits between first- and second-order convergence, whereas the temperature field converges at approximately first order. Convergence rates computed using the finest-grid numerical solution as the reference are shown in Fig.~\ref{fig_ConvRateRatio540}(B). In this case, the phase-front location exhibits nearly second-order convergence, while the temperature field converges at a rate between first and second order. For this case, the errors are evaluated over the time interval $0.01\le t\le0.5$ s to exclude the initial transient associated with the sudden appearance of the solid region, which is not explicitly modeled by the enthalpy method.

The results of this section verify that the improved low-Mach enthalpy method (i) converges to the mushy Stefan solution and (ii) remains stable and accurate for phase-change problems with large density contrasts between the two phases.

Having validated the improved low-Mach enthalpy method against the one-dimensional similarity solution of the mushy Stefan problem, we next consider two-dimensional phase-change problems representing defect formation in metal-manufacturing processes involving large density contrasts between the solid and liquid phases. These simulations demonstrate that the improved formulation remains stable and accurately captures the underlying phase-change dynamics, whereas the original low-Mach enthalpy method exhibits numerical instabilities or produces physically inconsistent results.

\subsection{Pipe shrinkage defect}\label{subsec_pipe_shrinkage_defect}

Pipe-shrinkage defects occur during casting when the solid metal is denser than the liquid metal. If the inlet feed velocity is not properly maintained, the liquid--gas free surface recedes as a result of the density difference between the solid and liquid phases, leading to cavity formation in the resulting solid. This problem was investigated in our previous work~\cite{thirumalaisamy2023low} for a density ratio of $\rhos/\rhol = 1.09$, representative of an aluminum-like alloy. To demonstrate the capabilities of the improved low-Mach enthalpy method, we revisit this problem for density ratios of $\rhos/\rhol = 1.09$ and $3.0$. 

We consider a square computational domain $\Omega = [0,L]^2$ with side length
$L = 8\times10^{-3}$~m. Initially, the region below $y =
5\times10^{-3}$~m is filled with liquid metal, while the remainder of the
domain is occupied by the gas phase. The thermophysical properties of the
solid and liquid phases, except for their densities, are taken from
Table~\ref{tab_thermophys_properties_Al}. The solid and liquid phases are
assigned the same dynamic viscosity,
$\mul=\mus=1.4\times10^{-3}$~kg/(m$\cdot$s). The liquid--gas surface-tension
coefficient is set to $\sigma=0.87$~N/m. The gas phase is assigned the
properties $\rhog=0.4$~kg/m$^3$, $\kg=6.1\times10^{-2}$~W/(m$\cdot$K),
$\cpg=1100$~J/(kg$\cdot$K), and
$\mug=4\times10^{-5}$~kg/(m$\cdot$s). Gravity acts in the negative
$y$-direction, $\mathbf{g}=[0,-9.81]$~m/s$^2$. Outflow boundary conditions are imposed on the top boundary, while no-slip boundary conditions are imposed on the remaining boundaries. A fixed temperature condition, $T=0.5T_m$, is applied on all boundaries except the bottom wall, which is treated as adiabatic. Here, $T_m$ denotes the PCM melting temperature. The initial temperature throughout the domain is $T=2T_m$.

\begin{figure}[]
	\centering{\includegraphics[width=\textwidth, height = 0.91\textheight,  trim={40 35 30 0},clip]{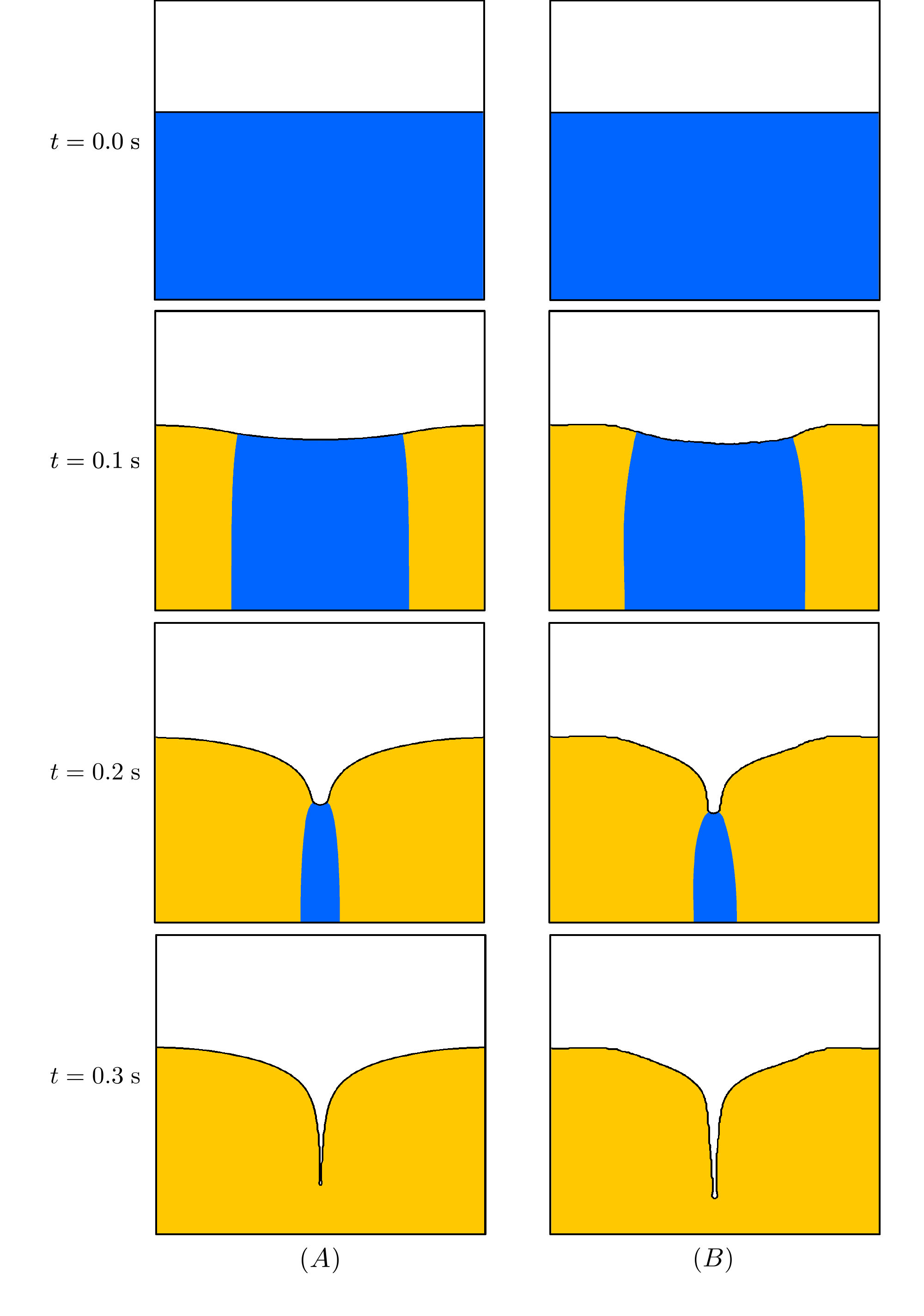}}
\caption{Evolution of metal solidification exhibiting pipe shrinkage defect for the low density ratio of $\rhos/\rhol = 1.09$ using (A) the original low-Mach enthalpy method and (B) the improved low-Mach enthalpy method. Solid, liquid, and gas regions are depicted in yellow, blue, and white colors, respectively. }\label{fig_Ratio1_1_Comparison_PipeShrinkage}
\end{figure}

\begin{figure}[]
	\centering{\includegraphics[width=1.0\textwidth, trim = {0 250 0 0}, clip]{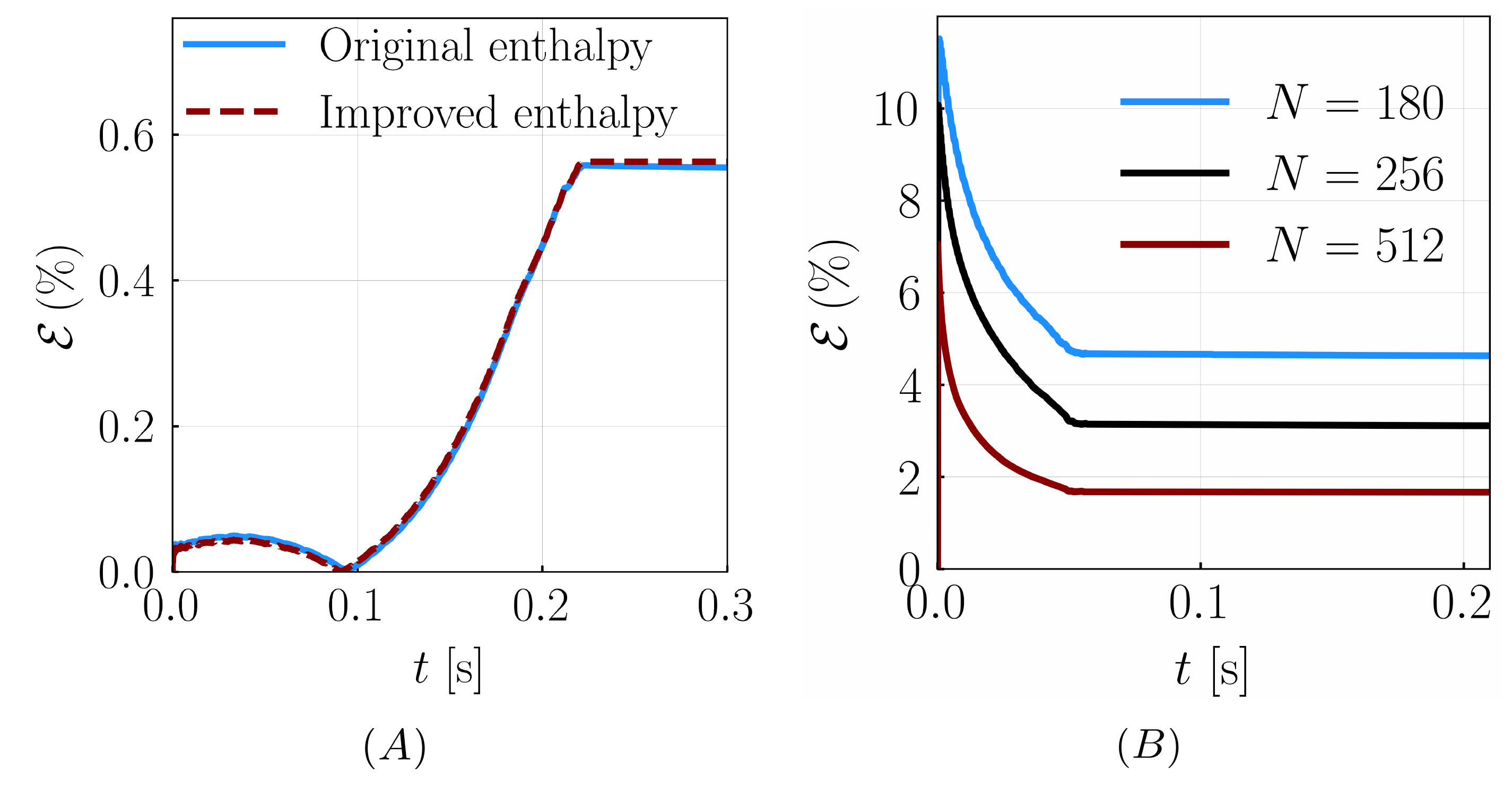}}
\caption {Comparison of the relative PCM mass error, $\mathcal{E}$, expressed as a
percentage, as a function of time for the pipe-shrinkage cases: (A) results
obtained using the original and improved formulations for
$\rhos/\rhol=1.09$; and (B) results obtained using the improved formulation
for $\rhos/\rhol=3.0$ at different grid resolutions.}\label{fig_PipeShrinkage3MassError}
\end{figure}

First, we consider the low-density-ratio case with
$\rhos/\rhol=1.09$, corresponding to solid and liquid densities of
$\rhos=2700$~kg/m$^3$ and $\rhol=2475$~kg/m$^3$, respectively. The objective
of this test is to verify that the improved formulation reproduces the
behavior of the original formulation in the low-density-ratio limit. The
domain is discretized using $N_x\times N_y=256\times256$ cells, with a
constant time step of $\Delta t=1.0\times10^{-5}$~s. The simulation is
continued until solidification is complete. Fig.~\ref{fig_Ratio1_1_Comparison_PipeShrinkage} compares the solidification dynamics, whereas Fig.~\ref{fig_PipeShrinkage3MassError}(A) compares the \% change in PCM mass over time with the original and improved formulations. Both formulations
produce nearly identical solidification dynamics. As solidification
progresses, the solid phase advances into the liquid region, while the
liquid--gas free surface gradually recedes, resulting in the formation of a
narrow central cavity. Both formulations predict a similar final change in
the PCM mass, approximately $0.6\%$. This close agreement demonstrates that
the improved formulation retains the behavior of the original formulation in
the low-density-ratio limit.

\begin{figure}[]
	\centering{\includegraphics[width=1.0\textwidth]{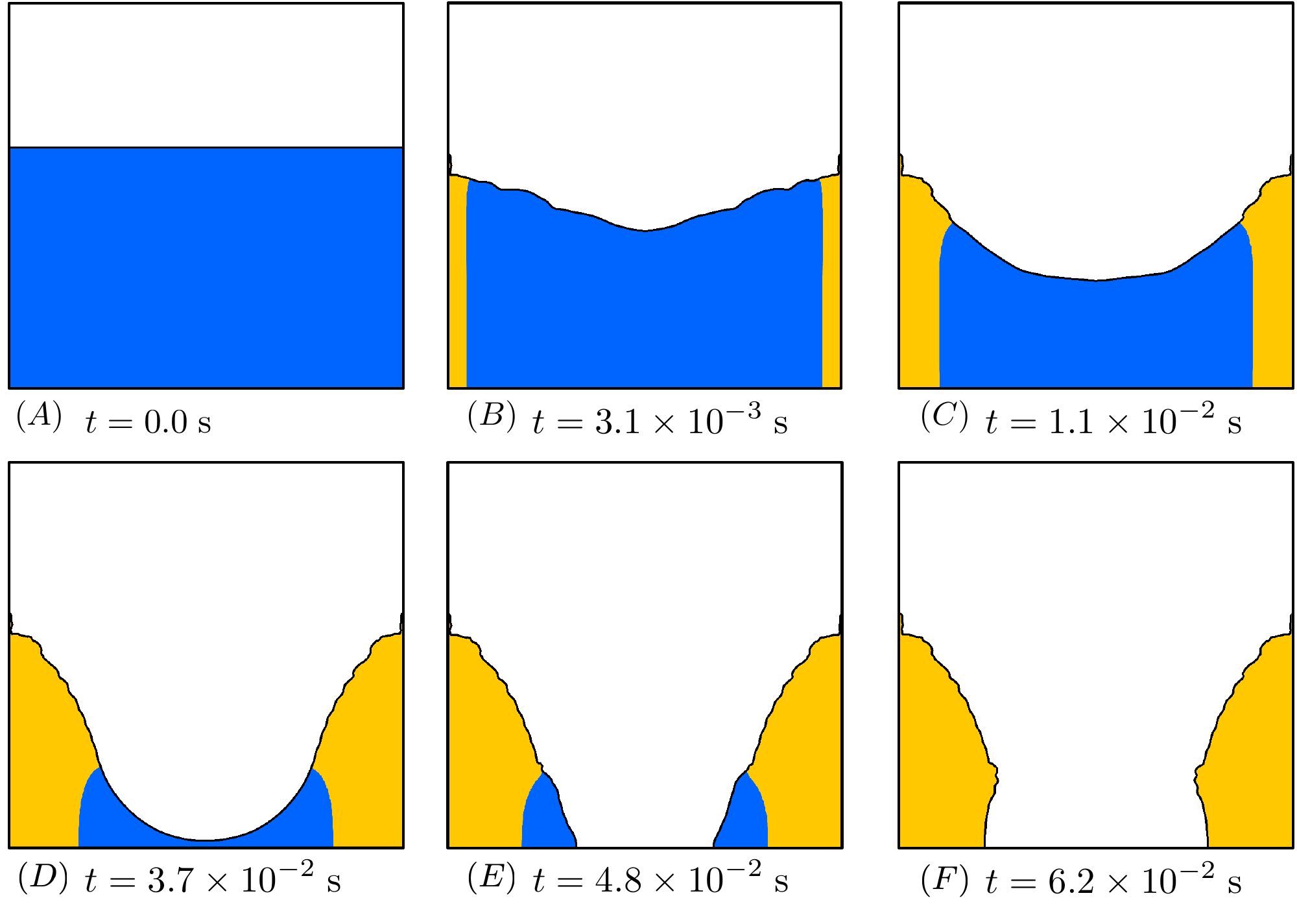}}
\caption{Evolution of metal solidification exhibiting a pipe shrinkage defect with the large density ratio of $\rhos/\rhol = 3.0$ using the improved low-Mach enthalpy method. Solid, liquid, and gas regions are depicted in yellow, blue, and white colors, respectively.}\label{fig_Ratio3PipeShrinkage}
\end{figure}

Next, we consider a large-density-ratio case with
$\rhos/\rhol=3.0$, corresponding to solid and liquid densities of
$\rhos=2700$~kg/m$^3$ and $\rhol=900$~kg/m$^3$, respectively. A time step of
$\Delta t=4.0\times10^{-6}$~s is employed, and the simulation is continued
until a steady state is reached. This case involves a substantially larger
volume change during solidification and, consequently, a greater recession of
the liquid--gas interface than in the low-density-ratio case. Fig.~\ref{fig_Ratio3PipeShrinkage} presents the simulation results. 
Solidification begins along the cooled sidewalls and progresses inward,
reducing the volume of the remaining liquid region. As solidification
continues, the free surface reaches the bottom boundary and pinches off,
forming a central cavity. The remaining liquid then solidifies into two
distinct solid regions separated by a central gas void.

Fig.~\ref{fig_PipeShrinkage3MassError}(B) shows the evolution of the PCM mass
change for $\rhos/\rhol=3.0$ on coarse, medium, and fine grids, corresponding
to $N_x\times N_y=\{180^2,256^2,512^2\}$ cells. The coarse- and medium-grid
simulations use $\Delta t=4.0\times10^{-6}$~s, whereas the fine-grid
simulation uses $\Delta t=1.0\times10^{-6}$~s. The magnitude of the PCM mass
change decreases with grid refinement, from approximately $4.5\%$ on the
coarsest grid to $1.8\%$ on the finest grid. When the density ratio of $\rhos/\rhol = 3.0$
is simulated using the original low-Mach enthalpy method, the solution
diverges rapidly (data not shown).

\subsection{Porosity defect}\label{subsec_porosity_defect}

Gas-bubble entrapment, commonly referred to as a porosity defect, occurs when
ambient gas or vapor generated by laser-induced metal vaporization becomes
trapped within the rapidly solidifying melt. Porosity reduces the mechanical
strength of the final part and can ultimately lead to failure. Numerical
simulations can help predict pore formation as a function of process
parameters, such as laser power and scan speed, as well as the flow of
ambient inert gas used to control oxidation of the molten metal.

In our previous work~\cite{thirumalaisamy2025consistent}, we investigated
porosity formation in an aluminum-based PCM with a solid-to-liquid density
ratio of $\rhos/\rhol=1.09$, representative of metals that undergo
solidification-induced volume shrinkage. Here, we revisit this problem using
the improved low-Mach enthalpy formulation and progressively increase the
solid-to-liquid density ratio to investigate pore-formation dynamics as a
function of the volumetric shrinkage of the material.

Consider a square computational domain $\Omega=[0,L]^2$ with side length
$L=10^{-2}$~m. Initially, the solid PCM occupies the region
$y\leq 0.1L$, while the liquid PCM extends from $y=0.1L$ to $y=0.75L$.
The remainder of the domain is occupied by the gas phase. Two gas bubbles are
initially placed within the liquid region; their locations and radii are
listed in Table~\ref{tab_bubbles_loc}.

The solid-phase density is fixed at $\rhos=2700$~kg/m$^3$, whereas the
liquid-phase density is varied as $\rhol=2475$, $500$, and $50$~kg/m$^3$,
corresponding to density ratios $\rhos/\rhol=1.09$, $5.4$, and $54$,
respectively. The solid and liquid phases are assigned the same dynamic
viscosity, $\mul=\mus=1.4\times10^{-3}$~kg/(m$\cdot$s). The liquid--gas
surface-tension coefficient is set to $\sigma=0.87$~N/m. The ambient gas is
assigned the properties $\rhog=0.4$~kg/m$^3$,
$\kg=6.1\times10^{-2}$~W/(m$\cdot$K),
$\cpg=1100$~J/(kg$\cdot$K), and
$\mug=4\times10^{-5}$~kg/(m$\cdot$s). Gravity acts in the negative
$y$-direction, $\mathbf{g}=[0,-9.81]$~m/s$^2$. 

No-slip velocity boundary conditions are imposed on the bottom and side
walls, while a zero-pressure outflow condition is applied at the top wall.
Adiabatic thermal boundary conditions are imposed on all walls except the
bottom wall, where the temperature is fixed at $T=0.5T_m$. Here, $T_m$
denotes the PCM melting temperature. Initially, all phases are at rest. The
initial temperature is $0.5T_m$ in the solid phase and $1.1T_m$ in the liquid
and gas phases. In all cases, the phase-change temperature interval is
$\Tliq-\Tsol=5$~K. The computational domain is discretized using a base grid of
$N_x\times N_y=128\times128$ cells with one additional level of adaptive mesh
refinement (AMR) and a refinement ratio of two. The PCM--gas and
solid--liquid interfaces are dynamically refined to maintain adequate
resolution of the evolving interfaces. Additional details regarding the AMR
tagging criteria are provided in~\cite{thirumalaisamy2025consistent}.

\begin{table}[]
\centering
\caption{Initial bubble locations and radii.}
\label{tab_bubbles_loc}
\begin{tabular}{clcl}
\midrule
Bubble	&Center	&Radius \\
\midrule
Left  	&[0.175L, 0.19L]	&0.05L     \\
Right 	&[0.8L, 0.21L]		&0.075L     \\
\bottomrule
\end{tabular}
\end{table}

\begin{figure}[]
	\centering{\includegraphics[width=\textwidth, height = 0.91\textheight,  trim={0 102 0 0},clip]{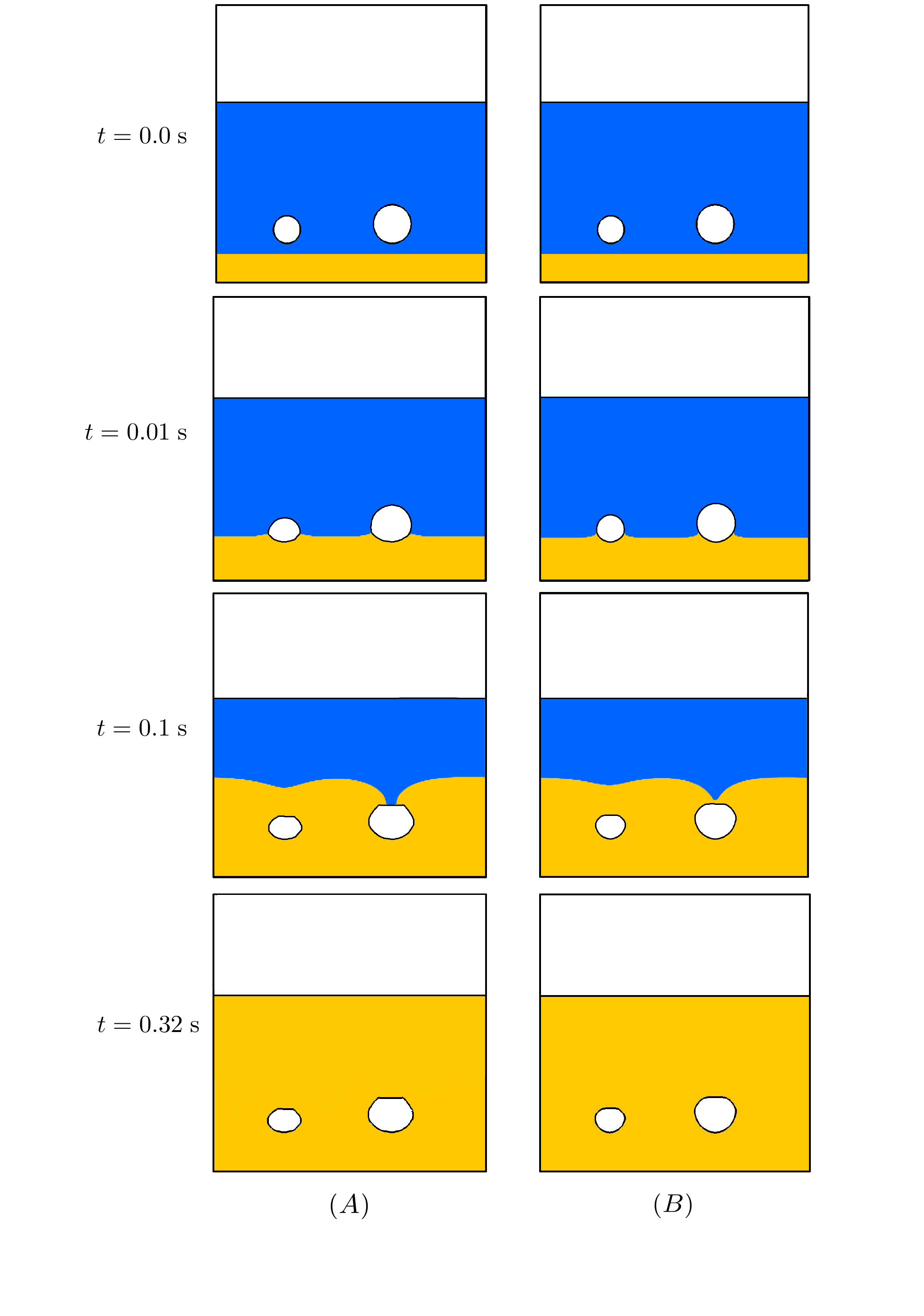}}
\caption{Evolution of metal solidification at the low density ratio
$\rhos/\rhol=1.09$ using (A) the original and (B) the improved low-Mach
enthalpy methods. The solid, liquid, and gas regions are shown in yellow,
blue, and white, respectively.}\label{fig_Ratio1_1_Comparison}
\end{figure}

We first consider the low-density-ratio case with
$\rhos/\rhol=1.09$, comparing the results obtained using the original and
improved low-Mach enthalpy formulations. The time step is fixed at
$\Delta t=1.0\times10^{-6}$~s. Fig.~\ref{fig_Ratio1_1_Comparison}
provides a direct comparison of the two methods. The gas bubbles initially
rise under buoyancy while solidification proceeds simultaneously from the
bottom of the liquid region. Both bubbles are subsequently intercepted by
the solidification front, after which the surrounding liquid solidifies, and
the bubbles become trapped within the solid PCM. As solidification continues,
the upper gas--liquid interface gradually recedes. This displacement is
relatively small for the present density ratio but becomes increasingly
pronounced as $\rhos/\rhol$ increases.

\begin{figure}[]
	\centering{\includegraphics[width=1.0\textwidth, trim = {0 250 0 0}, clip]{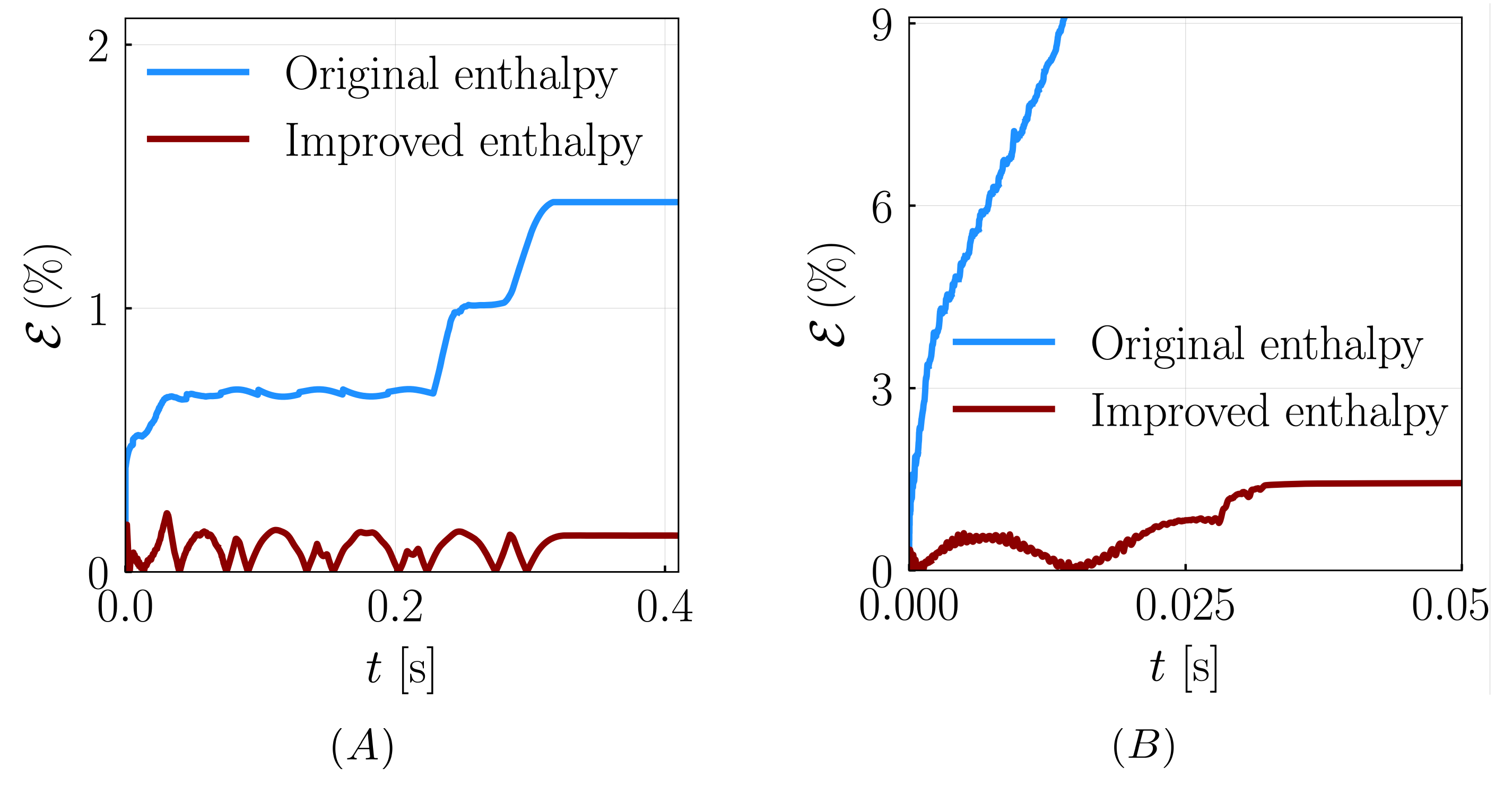}}
\caption{Comparison of the relative PCM mass error, $\mathcal{E}$, expressed as a
percentage, as a function of time for the metal-solidification case with two
gas bubbles, obtained using the original and improved low-Mach enthalpy
methods for density ratios of (A) $\rhos/\rhol=1.09$ and
(B) $\rhos/\rhol=5.4$.}\label{fig_Ratio1_1_MassComparison}
\end{figure}

The original and improved formulations predict similar solidification and
bubble dynamics. Fig.~\ref{fig_Ratio1_1_MassComparison}(A) compares the
relative PCM mass change during the simulations. The maximum mass change is
approximately $1.4\%$ for the original formulation and $0.3\%$ for the
improved formulation. Thus, both methods exhibit only a small change in PCM
mass. These results demonstrate that the improved low-Mach enthalpy
formulation retains the behavior of the original formulation in the
low-density-ratio limit, as expected.


\begin{figure}[]
	\centering{\includegraphics[width=1.0\textwidth, trim = {0 350 0 0}, clip]{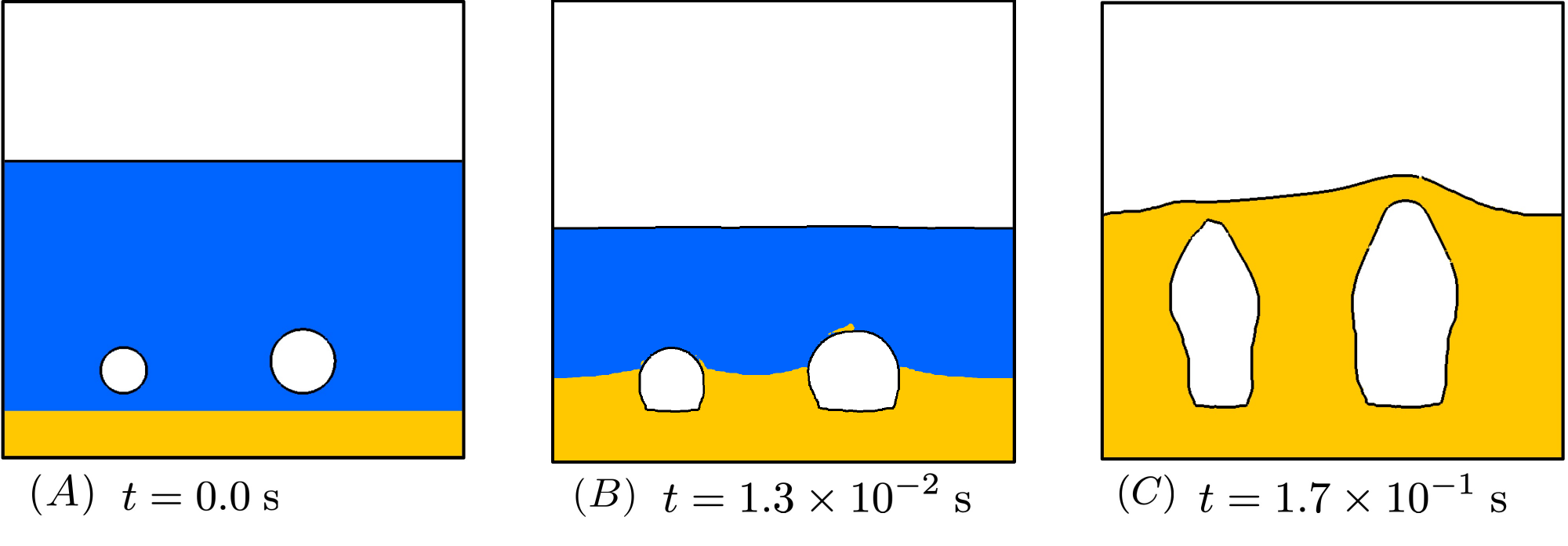}}
\caption{Evolution of metal solidification at the intermediate density ratio
$\rhos/\rhol=5.4$ using the original low-Mach enthalpy method. The solid,
liquid, and gas regions are shown in yellow, blue, and white, respectively.}\label{fig_Ratio5_4_Original}
\end{figure}

In the next test, we increase the density ratio to
$\rhos/\rhol=5.4$. The simulation is performed using a constant time step of
$\Delta t=1.0\times10^{-6}$~s. At this density ratio, the differences between
the original and improved low-Mach enthalpy formulations become more
pronounced. The original formulation exhibits a substantial change in PCM
mass and fails to predict a physically correct evolution of the gas--liquid
interface, rendering the results unreliable. Fig.~\ref{fig_Ratio5_4_Original}
shows the nonphysical results obtained with the original formulation.

\begin{figure}[]
	\centering{\includegraphics[width=1.0\textwidth]{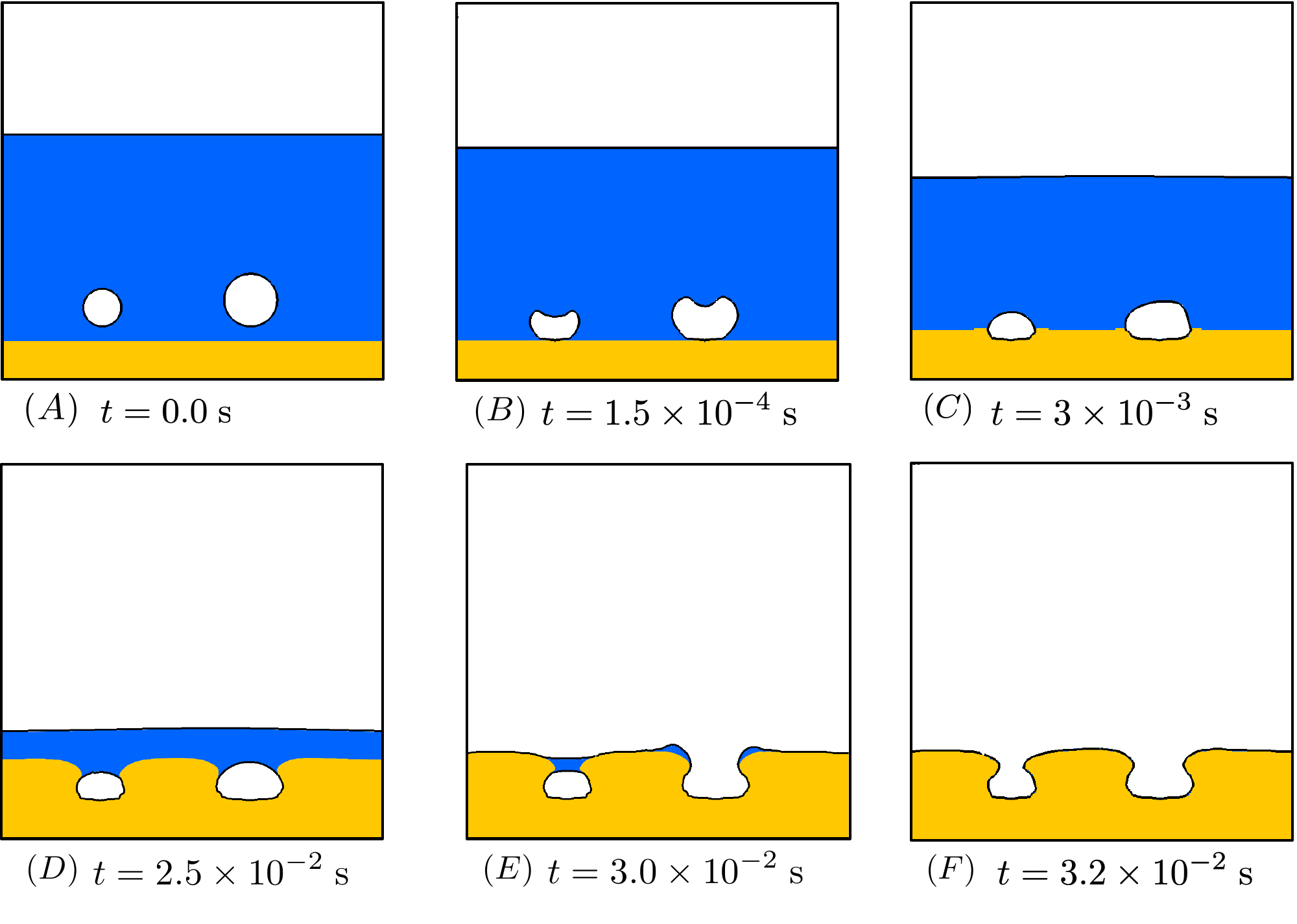}}
\caption{Evolution of metal solidification at the intermediate density ratio
$\rhos/\rhol=5.4$ using the improved low-Mach enthalpy method. The solid,
liquid, and gas regions are shown in yellow, blue, and white, respectively.}\label{fig_Ratio5_4_ProrosityImproved}
\end{figure}

Fig.~\ref{fig_Ratio5_4_ProrosityImproved} shows the solidification
evolution predicted by the improved formulation. Compared with the
low-density-ratio case, the larger density ratio produces a greater volume
change during solidification. Consequently, the upper gas--liquid free
surface recedes more rapidly. The gas bubbles exhibit only limited upward
displacement due to buoyancy. As the PCM solidifies around the bubbles and
the free surface approaches them, the free surface eventually intersects the
bubbles, allowing the trapped gas to escape into the ambient gas. This
process leaves two cavities in the solidified PCM that remain open to the
atmosphere.

Fig.~\ref{fig_Ratio1_1_MassComparison}(B) presents the PCM mass change
predicted by the improved formulation. The maximum relative mass change is
approximately $1.5\%$, which remains small despite the relatively large volume
change associated with solidification in this case.

\begin{figure}[]
	\centering{\includegraphics[width=1.0\textwidth]{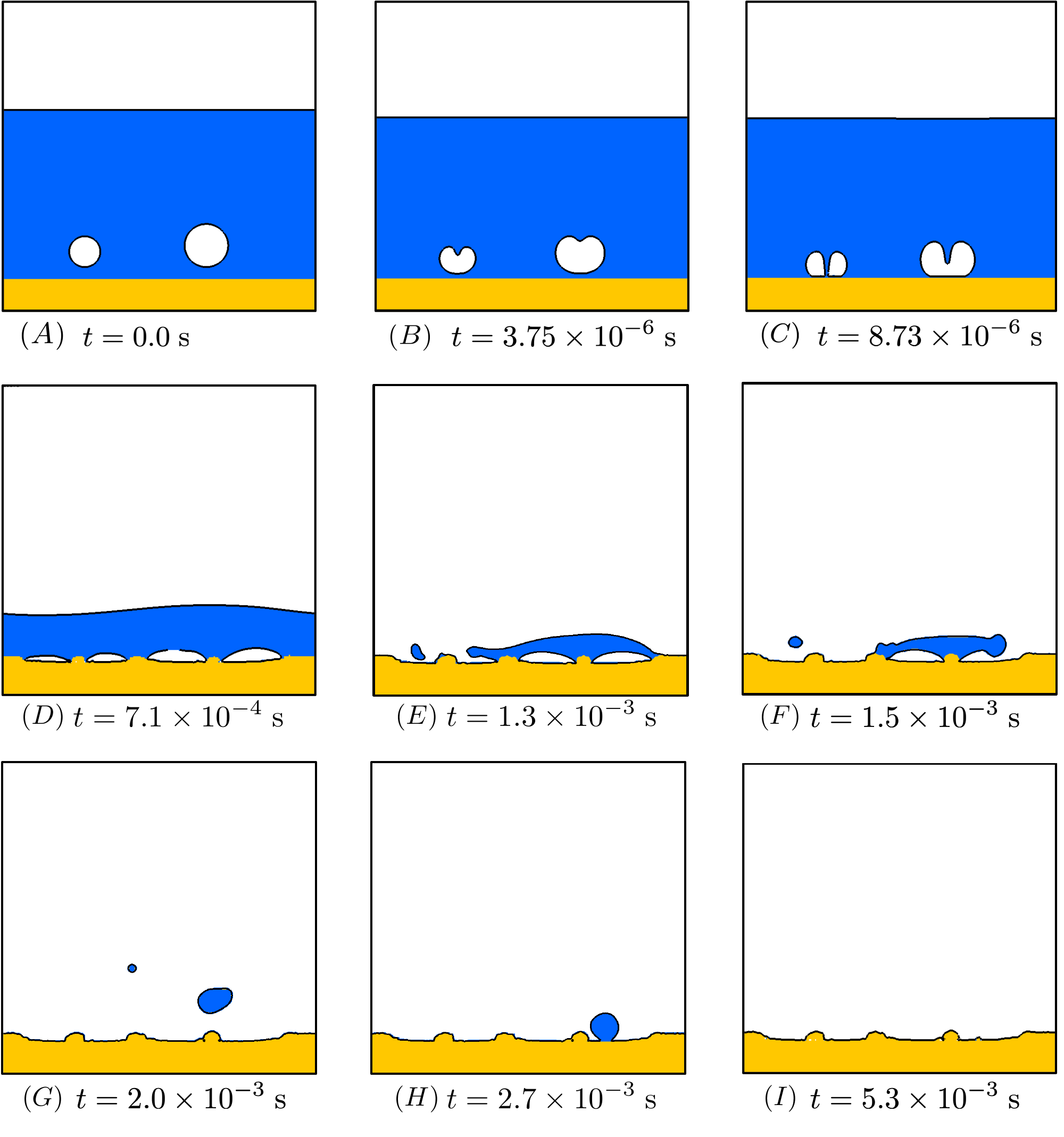}}
\caption{Evolution of metal solidification at the large density ratio
$\rhos/\rhol=54$ using the improved low-Mach enthalpy method. The solid,
liquid, and gas regions are shown in yellow, blue, and white, respectively.}\label{fig_Ratio54PorosityDefect}
\end{figure}

\begin{figure}[]
	\centering{\includegraphics[scale=0.415]{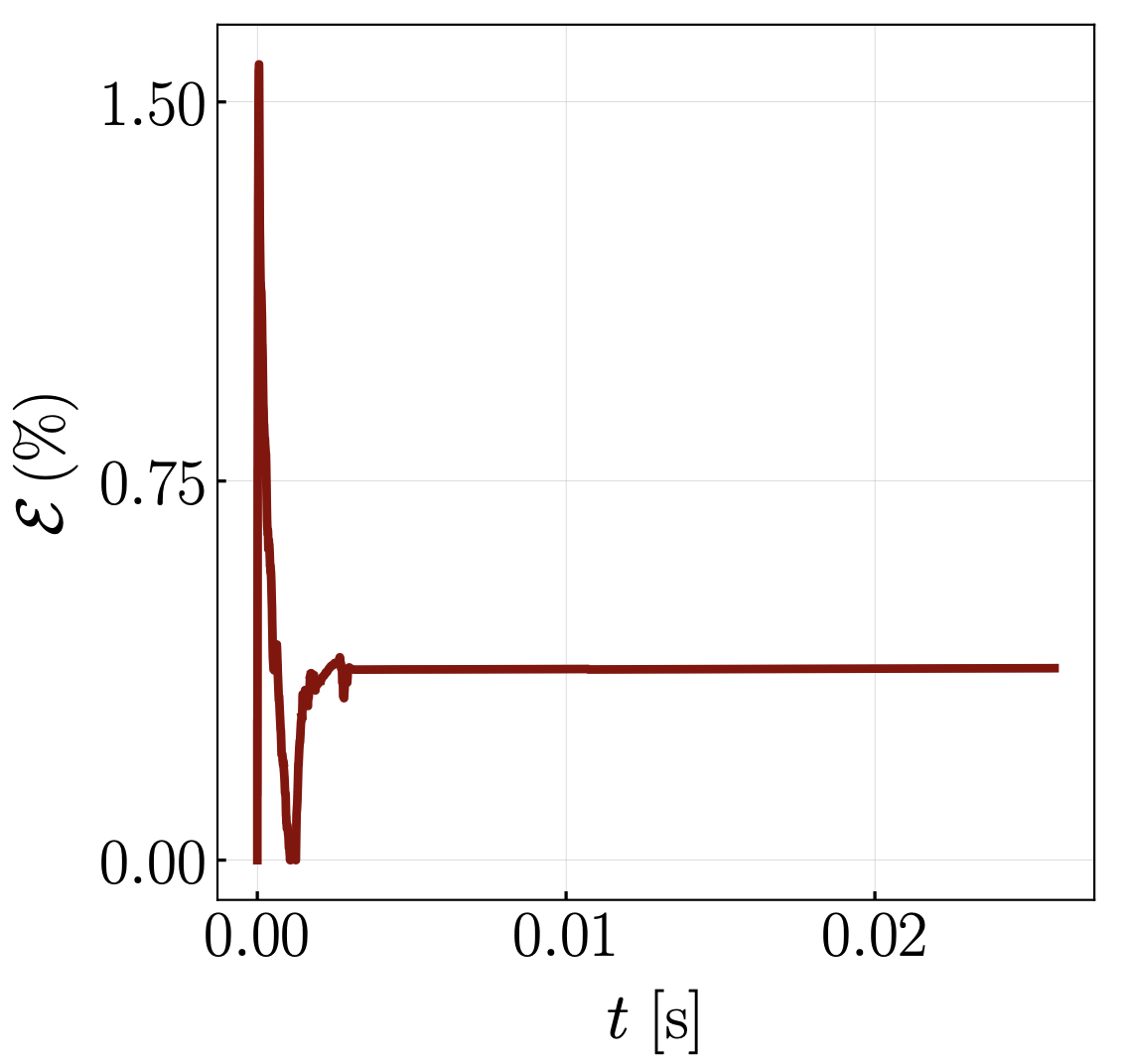}}
\caption{Relative PCM mass error, $\mathcal{E}$, expressed as a percentage, as a
function of time for the metal-solidification case with two gas bubbles at a
density ratio of $\rhos/\rhol=54$, obtained using the improved low-Mach
enthalpy method.}\label{fig_PorosityMassErrorRatio54}
\end{figure}

Finally, we consider an extreme hypothetical case with a large density ratio,
$\rhos/\rhol=54$. A reduced time step of
$\Delta t=8.0\times10^{-8}$~s is used. This simulation tests the robustness
of the formulation under the very large volume changes associated with
solidification. The original low-Mach enthalpy method becomes numerically
unstable shortly after the onset of solidification and therefore fails to
produce reliable results (data not shown). Fig.~\ref{fig_Ratio54PorosityDefect} shows the evolution of the solid phase
predicted by the improved method. Because of the rapid volume contraction,
the initial gas--liquid interface recedes downward at a much faster rate.
This motion suppresses the upward buoyant motion of the gas bubbles, which
remain below their initial positions. As the liquid layer contracts, the
gas--liquid interface eventually reaches the bubbles, allowing the trapped
gas to escape. The large volume change associated with solidification then
ejects the remaining liquid layer from the solidification front. The ejected
liquid forms a drop that travels through the ambient gas, becomes
approximately spherical under the action of surface tension, and subsequently
lands on the solid, where it completely solidifies. In this case, the bubble dynamics are governed primarily by the rapid recession of the gas--liquid PCM free surface rather than by the buoyant rise of the bubbles. After solidification is complete, the final PCM mass change is approximately $0.4\%$, as shown in Fig.~\ref{fig_PorosityMassErrorRatio54}.

These results demonstrate the ability of the improved low-Mach enthalpy
method to robustly solve increasingly challenging problems as the density
contrast between the solid and liquid metal phases increases.

\section{Conclusions and future work}\label{sec_conclusions}

This work fills a longstanding gap in the enthalpy method literature by deriving the similarity solution of the non-isothermal variable-density mushy Stefan problem. The resulting similarity solution provides a rigorous benchmark for validating enthalpy methods at finite values of the phase-change interval, $\Delta T$. We showed that, for equal- and low-density-ratio cases, the mushy and sharp-interface Stefan problems remain in close agreement over a broad range of $\Delta T$ values. In contrast, for large-density-ratio cases, $\Delta T$ must be sufficiently small for the mushy Stefan problem to closely approximate the sharp-interface Stefan problem. Thus, large-density-ratio cases provide a significantly more stringent validation benchmark for modern enthalpy methods that account for density-change-induced flow.

The similarity formulation also enabled improvements to the low-Mach enthalpy method, yielding a formulation that remains stable and accurate for phase-change problems involving large density contrasts between the two phases, e.g., density ratios of $\mathcal{O}(3)$. Specifically, at large density ratios, the original low-Mach enthalpy method exhibited unstable evolution of the phase-change front regardless of the magnitude of the Carman--Kozeny drag force. The primary difference between the original and improved formulations lies in the definition of the mushy region specific enthalpy, which gives rise to nonlinear and linear liquid-fraction--temperature ($\varphi$--$T$) relations, respectively. It was further demonstrated that the low-Mach enthalpy method solves the mushy Stefan problem rather than the classical Stefan problem, with the latter recovered only in the limiting case as $\Delta T\rightarrow0$.

We also used the improved low-Mach enthalpy method to simulate
metal-manufacturing defects, including pipe shrinkage and porosity. The
improved formulation enabled simulations at large density ratios and captured
the complex dynamics associated with substantial solidification-induced volume
shrinkage. In contrast, simulations at these density ratios using the original
formulation failed or produced nonphysical results.

Although enthalpy methods have traditionally been applied to melting and
solidification problems involving density ratios close to unity, the
theoretical and numerical developments presented in this work extend their
potential applicability to boiling and condensation problems, which involve
density changes of $\mathcal{O}(3)$ between the liquid and vapor phases.
Applying the improved enthalpy method to boiling and condensation problems
remains an important direction for future work. Overall, this study provides
a rigorous theoretical foundation and a robust numerical framework for
extending enthalpy methods to a broader class of variable-density
phase-change problems.

\section*{Acknowledgements}
A.P.S.B~acknowledges support of NSF award CBET CAREER 2234387. 

\section*{Supplementary material}
\noindent \textbf{Supplementary video S1}: Time evolution of the solid, liquid, and gas regions during metal
solidification with pipe-shrinkage defects for density ratios of
$\rhos/\rhol=1.09$ and $3.0$, obtained using the improved low-Mach enthalpy method (see Sec.~\ref{subsec_pipe_shrinkage_defect}). \\
\noindent \textbf{Supplementary video S2}: Time evolution of the solid, liquid, and gas regions during metal
solidification with porosity defects for density ratios of
$\rhos/\rhol=1.09$, $5.4$, and $54$, obtained using the improved low-Mach
enthalpy method (see Sec.~\ref{subsec_porosity_defect}).
 
\bibliography{Ref.bib}
\end{document}